\documentclass[a4paper,reqno, 11pt]{amsart}

\usepackage[margin=3cm]{geometry}
\usepackage{amssymb}
\usepackage{amsmath}
\usepackage{amsthm}
\usepackage{mathrsfs}
\usepackage{enumitem}
\usepackage{hyperref}

\usepackage{tikz}
\usepackage{tikz-cd}
\usetikzlibrary{matrix,calc,arrows,decorations.pathmorphing}
\tikzset{node distance=2cm, auto}

\tikzset{cd/.style=matrix of math nodes,row sep=2em,column sep=2em, text height=1.5ex, text depth=0.5ex}
\tikzset{cdar/.style=->,auto}
\tikzset{mid/.style={anchor=mid}} 
\tikzset{narrowfill/.style={inner sep=1pt, fill=white}}

\newtheorem{theorem}{Theorem}[section]
\newtheorem{lemma}[theorem]{Lemma}
\newtheorem{corollary}[theorem]{Corollary}

\newtheorem{proposition}[theorem]{Proposition}
\theoremstyle{definition}
\newtheorem{definition}[theorem]{Definition}
\newtheorem{example}[theorem]{Example}

\newcommand{\RR}{\mathbb{R}}
\newcommand{\CC}{\mathbb{C}}
\newcommand{\NN}{\mathbb{N}}
\newcommand{\ZZ}{\mathbb{Z}}

\newcommand{\Kk}{\mathcal{K}}
\newcommand{\TT}{\mathbb{T}}
\newcommand{\Aa}{\mathcal{A}}

\newcommand{\Ll}{\mathcal{L}}
\newcommand{\Ff}{\mathcal{F}}

\newcommand{\Cc}{\mathcal{C}}

\newcommand{\Ee}{\mathcal{E}}
\newcommand{\Gg}{\mathcal{G}}
\newcommand{\Ww}{\mathcal{W}}

\newcommand{\bB}{\mathfrak{B}}
\newcommand{\mM}{\mathfrak{M}}

\newcommand{\sspan}{\operatorname{span}}

\newcommand{\eps}{\varepsilon}

\newcommand{\Id}{\operatorname{id}}

\newcommand{\supp}{\operatorname{supp}}
\newcommand{\ess}{\mathrm{ess}}
\newcommand{\red}{\mathrm{r}}

\newcommand{\pt}{\mathrm{pt}}

\newcommand{\Cart}{\mathrm{Cartan}}
\newcommand{\EssCart}{\mathrm{EssCart}}
\newcommand{\Fell}{\mathrm{Fell}}
\newcommand{\Haus}{\mathrm{Haus}}

\newcommand{\Mloc}{M_{\mathrm{loc}}}

\newcommand{\Bis}{\mathrm{Bis}}

\newcommand{\dom}{\mathrm{dom}}

\newcommand{\llangle}{\langle\langle}
\newcommand{\rrangle}{\rangle\rangle}

\newcommand{\baltimes}[2]{
	\mathbin{_{#1}\times_{#2}}%
}

\theoremstyle{remark}
\newtheorem{remark}[theorem]{Remark}

\title{Functoriality for groupoid and Fell bundle $C^*$-algebras}

\author{Jonathan Taylor}
\address{Jonathan Taylor, Institut für Mathematik, University of Potsdam, Campus Golm, Haus 9, Karl-Liebknecht-Str. 24-25, 14476, Germany}
\email{jonathan.taylor@uni-potsdam.de}

\begin{document}
	
	\begin{abstract}
		We define a class of morphisms between \'etale groupoids and show that there is a functor from the category with these morphisms to the category of $C^*$-algebras. 
		We show that all homomorphisms between Cartan pairs of $C^*$-algebras that preserve the Cartan structure arise from such morphisms between the underlying Weyl groupoids and twists, and attain an equivalence of categories between Cartan pairs with structure preserving homomorphisms and a their associated twists.
		We define analogous morphisms for Fell bundles over $C^*$-algebras and show these functorially induce ${}^*$-homomorphisms between the Fell bundle $C^*$-algebras.
		
		We also construct colimit groupoids and Fell bundles for inductive systems of such morphisms, and show that the functor to $C^*$-algebras preserves these colimits.
	\end{abstract}
	
	\maketitle
	
	\section{Introduction}
	
	$C^*$-algebras of (twists over) \'etale groupoids have been the object of study for many authors for various reasons.
	They generalise discrete group $C^*$-algebras, crossed products of actions on topological spaces by discrete groups, graph algebras, as well as providing a wide menagerie of tractable objects to study within the category of $C^*$-algebras.
	When studying a groupoid $C^*$-algebra, one may of course use the structure of the underlying groupoid to determine qualities of the algebra.
	This is, of course, one of the main motivations for considering such algebras, and also a motivation for determining which $C^*$-algebras are isomorphic to groupoid $C^*$-algebras.
	
	Kumjian \cite{K1} showed that $C^*$-algebras with a $C^*$-diagonal must arise as $C^*$-algebras of twists over principal \'etale groupoids.
	Extending upon this, Renault \cite{R1} showed that one may relax the diagonal condition of Kumjian a little and recover a similar twist over an effective groupoid for $C^*$-algebras with a Cartan subalgebra.
	This gave a $C^*$-analogue to the known case for von Neumann algebras and their Cartan subalgebras, going back to the work of von Neumann.
	
	An extremely strong connection to the Elliott classification program was later established by Li \cite{L1}, where it was shown that inductive limits of Cartan pairs of $C^*$-algebras give a Cartan structure to their inductive limits, so long as the morphisms in the diagram preserve the relevant Cartan structure.
	Since all Elliot-classifiable $C^*$-algebras arise as such inductive limits, the immediate upshot is that all such $C^*$-algebras have twisted groupoid models.
	Moreover, Barlak and Li \cite{BL1} explicitly construct the inductive limit groupoids and twists from the $C^*$-algebraic data, giving a firmer handle for constructing and analysing such $C^*$-algebras.
	To do so, they rely heavily on the construction of Renault for the groupoid models of the building block Cartan pairs, rather than considering ${}^*$-homomorphisms between generic groupoid $C^*$-algebras.
	
	The construction of the groupoid $C^*$-algebra (and various flavours thereof) is not factorial for groupoid homomorphisms (i.e., continuous functors) in general.
	The obstruction to functoriality arises as a byproduct of groupoids being a generalisation of both groups and topological spaces, but rarely fitting into either category neatly.
	The construction of the group $C^*$-algebra is covariantly functorial for group homomorphisms, whereas for topological spaces, Gelfand duality is contravariant.
	Hence, for the groupoid $C^*$-algebra construction to be a functorial for groupoid homomorphisms, it would need to be contravariant on the $C^*$-algebras associated to the unit spaces of groupoids, while simultaneously working covariantly on source and range fibres within the groupoid.
	Barlak and Li \cite{L1} identify two classes of groupoid homomorphisms for which the groupoid $C^*$-algebra construction is functorial.
	They consider firstly open embeddings of groupoids into other groupoids, which covariantly induce inclusions of the $C^*$-algebras.
	Next, they consider surjective, proper, and fibrewise bijective groupoid homomorphisms, which induce ${}^*$-homomorphisms contravariantly between the groupoid $C^*$-algebras via pullback.
	The fibrewise bijectivity of the homomorphisms ensures that, on each fibre, one may invert the homomorphism to gain a homomorphism in the (now covariant) other direction, mimicking the case for groups.
	Li constructs his inductive limit groupoids by considering zig-zags of these two classes of groupoid homomorphisms, and taking projective and inductive limits along the resulting `grid' of these homomorphisms.
	
	To avoid the mud arising from groupoid homomorphisms in the $C^*$-algebra construction, one may ask to consider a different class of morphisms between groupoids.
	Buneci and Stachura \cite{BS1} detail a construction of \emph{algebraic morphisms} for groupoids with Haar system, which goes back to Zakrzewski.
	Instead of considering functions between groupoids, one considers actions of groupoids on other groupoids with certain compatibility conditions.
	Buneci and Stachura show that these morphisms functorially induce ${}^*$-homomorphisms between certain types of groupoid $C^*$-algebras.
	Meyer and Zhu \cite{MZ1} refer to such actions with this compatibility as \emph{actors}, and show that actors covariantly induce functors between categories of actions by groupoids.
	
	We define actors between \'etale groupoids following Buneci-Stachura and Meyer-Zhu, and show that the construction of the induced ${}^*$-homomorphisms continues to function without requiring the groupoids to be Hausdorff.
	Motivated by the results of Renault \cite{R1} as well as the author \cite{T1} regarding Cartan pairs and twists, we determine conditions under which the induced ${}^*$-homomorphism arising from an actor entwines the obvious Cartan-like structure of groupoid $C^*$-algebras.
	We are able to show that the category of groupoids (and twists) giving rise to Cartan pairs, considered with actors as morphisms, is covariantly equivalent to the category of Cartan pairs with non-degenerate structure-preserving ${}^*$-homomorphisms.
	We are also able to demonstrate that the construction of the inductive limits of Li generalises comfortably to inductive systems of actors, and the intermediate groupoids Li constructs from morphisms of Cartan pairs arise as transformation groupoids associated to actors.
	The construction of the inductive limit groupoid is then compatible with the $C^*$-functor, in that the $C^*$-algebra of the inductive limit groupoid is naturally isomorphic to the inductive limit of the $C^*$-algebras of the building blocks.
	
	We also consider Fell bundles over groupoids, and generalise the definition of actors to this setting.
	Such \emph{Fell actors} include an underlying actor between the base-space groupoids as part of their defining data.
	We are able to construct inductive limit Fell bundles for certain classes of inductive systems of Fell bundles, and gain a similar compatibility with the $C^*$-functor as in the case for groupoids.
	In particular, the inductive limit Fell bundle is a bundle over the inductive limit of the underlying groupoids in the system.
	
	Study of the relationship between Fell bundles and $C^*$-algebras has enjoyed recent advancements by Bice in \cite{B1}, where a duality between inclusions of $C^*$-algebras equipped with a kind of dynamical structure and certain classes of Fell bundles is established.
	Many of the results in this article have analogous counterparts in the article \cite{B1}, but the approach and perspective of the proofs in this article differ.
	There are also two main ways in which the results in this article cover topics not considered in \cite{B1}.
	The first is in the class of groupoids considered: in this article, we consider $C^*$-inclusions with generalised conditional expectations, which correspond to the underlying groupoid models being non-Hausdorff.
	Bice considers genuine expectations, which restricts the scope to Fell bundles over Hausdorff groupoids which allows for powerful duality results.
	Secondly, the perspective of this article allows us to naturally consider inductive limits of certain classes of actors, and show that the functor to $C^*$-algebras preserves them.
	
	\section*{Acknowledgements}
	The author would like to thank their doctoral supervisor and postdoctoral mentor, Ralf Meyer, for the invaluable advice and insight provided throughout this project.
	Further thanks go to Tristan Bice, who pointed out the article \cite{B1} has significant thematic overlap with the an early version of this article, and helped contextualise the results with established theory.
	Thanks to Ali I. Raad for revealing a number of small errata in the details of the article.
	
	The author was supported by the RTG2491 as an associated member.
	The majority of the work on this project was completed while the author was employed as a Postdoctoral Researcher at the University of Göttingen.
	
	\section{Groupoids and actors}
	
	We shall set our notation and recall some important facts regarding groupoids and their $C^*$-algebras here, and refer the reader to \cite{S1} and \cite[Section~7]{KM2} for further details.
	
	A \emph{groupoid} is a small category in which every arrow is invertible.
	Equivalently (following \cite{S1}), a groupoid is a set $G$ together with a distinguished subset $G^{(2)}\subseteq G\times G$ of \emph{composable pairs} together with a multiplication map $G^{(2)}\to G$ denoted by juxtaposition, and an inverse map $\gamma\to\gamma^{-1}$ from $G$ to $G$, such that
	\begin{enumerate}
		\item $(\gamma^{-1})^{-1}=\gamma$ for all $\gamma\in G$;
		\item if $(\alpha,\beta), (\beta,\gamma)\in G^{(2)}$ are composable pairs, then $(\alpha\beta,\gamma), (\alpha,\beta\gamma)\in G^{(2)}$ are also, and $(\alpha\beta)\gamma=\alpha(\beta\gamma)$;
		\item $(\gamma,\gamma^{-1})\in G^{(2)}$ is a composable pair for all $\gamma\in G$, and for all $\eta\in G$ with $(\gamma,\eta)\in G^{(2)}$ we have $\gamma^{-1}(\gamma\eta)=\eta$ and $(\gamma\eta)\eta^{-1}=\gamma$ for all $\gamma,\eta\in G$.
	\end{enumerate}
	From the categorical perspective, we identify all objects with their identity arrows, and then the range and source maps become $s(\gamma)=\gamma^{-1}\gamma$ and $r(\gamma)=\gamma\gamma^{-1}$ for all arrows $\gamma\in G$.
	We call $G^{(0)}:=\{\gamma^{-1}\gamma:\gamma\in G\}$ the space of \emph{units}, and note that it is the image of both the range and source maps.
	
	A \emph{topological groupoid} is a groupoid equipped with a topology for which the inverse map and the multiplication are continuous (where $G^{(2)}$ takes the subspace topology of the product topology on $G\times G$).
	Automatically the range and source maps are continuous.
	
	An \emph{\'etale groupoid} is a topological groupoid for which the range (equivalently the source) map is a local homeomorphism.
	For groupoids with only one unit (groups) the source and range maps are constant, so \'etale groups are exactly discrete groups.
	Similarly, transformation groupoids for actions of discrete groups on spaces (defined later in Definition~\ref{defn-grpActionTrsfmGrpd}) give rise to \'etale groupoids.
	
	We view the unit space of a groupoid as an underlying topological space, and the \'etale condition says that the entire groupoid locally has the topological structure of the unit space.
	In fact, a groupoid being \'etale is equivalent to it having a basis of open bisections: open subsets to which the range and source maps restrict to homeomorphisms onto their images.
	In particular, the unit space of a groupoid is an open bisection; it is open since it is the range of an open map (namely the range map).
	As we wish to construct $C^*$-algebras from \'etale groupoids, the unit spaces of such groupoids should lift to $C^*$-algebras modelling topological spaces.
	By Gelfand duality, such a $C^*$-algebra should be a commutative $C^*$-algebra that has the unit space of this groupoid as its spectrum.
	Hence we restrict our attention to \'etale groupoids with locally compact Hausdorff unit space.
	The whole groupoid is then locally compact and locally Hausdorff, but not necessarily globally Hausdorff.
	Unless otherwise stated, we shall consider \'etale groupoids with locally compact Hausdorff unit space throughout this article.
	
	The source and range fibres over a unit $x\in G^{(0)}$ are defined as $G_x:=s^{-1}(\{x\})$ and $G^x:=r^{-1}(\{x\})$.
	If $G$ is \'etale then these sets are discrete.
	
	Let $G$ be an \'etale groupoid with locally compact Hausdorff unit space $G^{(0)}$.
	For an open bisection $U\subseteq G$, we denote by $\Cc_c(U)$ the functions $G\to\CC$ with compact support contained in $U$ that have continuous restriction to $U$.
	That is, elements of $\Cc_c(U)$ are continuous functions $U\to\CC$ with compact support which are extended to functions on $G$ by zero.
	Note that if $G$ is Hausdorff, then these are exactly continuous functions $G\to\CC$ with compact support contained in $U$, but if $G$ is not Hausdorff, these functions need not be continuous on all of $G$.
	This particular phenomenon arises since the compact support of such a function may not be closed if $G$ is not Hausdorff.
	We define $\Cc_c(G)$ as the linear span of the spaces $\Cc_c(U)$ over all open bisections $U\subseteq G$ (considered in the ambient vector space of all functions $G\to\Cc$).
	Since the open bisections of $G$ form a basis, if $G$ is Hausdorff then $\Cc_c(G)$ is (by a partition of unity argument) exactly $C_c(G)$, the algebra of continuous compactly supported functions $G\to\CC$.
	
	We define a convolution product on $\Cc_c(G)$ as follows: for $f,g\in\Cc_c(G)$ and $\gamma\in G$ we define
	$$(fg)(\gamma):=\sum_{\eta\in G_{r(\gamma)}}f(\eta^{-1})g(\eta\gamma),$$
	that is, the sum over values $f(\alpha)g(\beta)$ for all composable pairs $(\alpha,\beta)\in G^{(2)}$ with $\alpha\beta=\gamma$.
	We define an involution ${}^*$ on $\Cc_c(G)$ via
	$$f^*(\gamma)=\overline{f(\gamma^{-1})}.$$
	The \emph{groupoid $C^*$-algebra} $C^*(G)$ of $G$ is then defined as the completion of $\Cc_c(G)$ in the maximal $C^*$-norm.
	Since $G^{(0)}$ is open, the algebra $C_0(G^{(0)})$ embeds naturally into $C^*(G)$ by continuously extending the inclusion $\Cc_c(G^{(0)})\subseteq\Cc_c(G)$.
		
	If $G$ is an \'etale groupoid, the collection $\Bis(G)$ of open bisections of $G$ forms an inverse semigroup under pointwise product $U\cdot V:=\{\gamma\cdot\eta:\gamma\in U,\eta\in V,s(\gamma)=r(\eta)\}$ and inverse $U^*=U^{-1}=\{\gamma^{-1}:\gamma\in U\}$ for all $U,V\in G$.
	
	For the following definition, we recall that the \emph{pullback} or \emph{fibred product} of two maps $f:X\to Z$ and $g:Y\to Z$ is the set $X\baltimes{f}{g}Y:=\{(x,y)\in X\times Y: f(x)=g(y)\}$.
	If $X$ and $Y$ are topological spaces, we endow $X\baltimes{f}{g}Y$ with the subspace topology of the product topology on $X\times Y$.
	
	\begin{definition}[{\cite[Definition and Lemma~4.1]{MZ1}}]\label{defn-action}
		Let $G$ be a topological groupoid and let $X$ be a topological space.
		A \emph{right action} $h$ of $G$ on $X$ consists of a continuous map $\sigma:X\to G^{(0)}$ (called the anchor) and a continuous multiplication $X\baltimes{\sigma}{r} G\to X$ (denoted by $\cdot_h$) satisfying
		\begin{enumerate}
			\item $\sigma(x\cdot_h g)=s(g)$ for all $(x,g)\in X\baltimes{\sigma}{r} G$;
			\item $(x\cdot_h g_1)\cdot_h g_2=x\cdot_h(g_1g_2)$ for all $(x,g_1)\in X\baltimes{\sigma}{r} G$ and $g_2\in G$ with $s(g_1)=r(g_2)$;
			\item $x\cdot_h \sigma(x)=x$ for all $x\in X$.
		\end{enumerate}
		
		A \emph{left action} $h$ of $G$ on $X$ consists of a continuous map $\rho_h:X\to G^{(0)}$ (called the anchor) and a continuous multiplication $G\baltimes{s}{\rho} X\to X$ (denoted by $\cdot_h$) satisfying
		\begin{enumerate}
			\item $\rho(g\cdot_h x)=r(g)$ for all $(g,x)\in G\baltimes{s}{\rho} X$;
			\item $g_1\cdot_h (g_2\cdot_\alpha x)=(g_1g_2)\cdot_h x$ for all $(g_2,x)\in G\baltimes{s}{\rho} X$ and $g_1\in G$ with $r(g_2)=s(g_1)$;
			\item $\rho(x)\cdot_h x=x$ for all $x\in X$.
		\end{enumerate}
		If the action $h$ is understood then we may omit the subscript from the multiplication and write $g\cdot x$ for $g\cdot_h x$.
	\end{definition}

	Any groupoid acts on itself from both the left and right by the groupoid multiplication.
	The left and right anchor maps are respectively the source and range maps of the groupoid.
	
	\begin{definition}[{\cite[Definition~2]{BS1}, \cite[Definition~4.15]{MZ1}}]\label{defn-actor}
		Let $G$ and $H$ be topological groupoids.
		An \emph{actor} from $G$ to $H$, denoted $G\curvearrowright H$, is a left action of $G$ on $H$ that commutes with the right multiplication of $H$ on itself.
		That is, $(\gamma,x)\in G\baltimes{s}{\rho}H$ and $(x,y)\in H^{(2)}$ we have $\rho(xy)=\rho(x)$, $s(\gamma\cdot x)=s(x)$, and $\gamma\cdot (xy)=(\gamma\cdot x)y$, where $\rho:H\to G^{(0)}$ is the anchor map of the action.
	\end{definition}
	
	\begin{remark}
		Given an actor $G\curvearrowright H$, since $s(\gamma\cdot x)=s(x)$ for any composable $(\gamma,x)\in G\baltimes{s}{\rho}H$, one must exercise a modicum of caution when constructing examples of actors, as some common examples of actions are not actors.
		For example, if $H$ consists only of units, every actor on $H$ is trivial, since $\gamma\cdot x\in H=H^{(0)}$ implies that $\gamma\cdot x=s(\gamma\cdot x)=s(x)=x$.
		In particular, apart from the trivial actions, actions of discrete groups on topological spaces by homeomorphisms are not actors, but the action of the group on the resulting transformation groupoid is (see Definition~\ref{eg-transGrpdActor}).
	\end{remark}
	
	Meyer and Zhu show that an actor $G\curvearrowright H$ induces a functor between the categories of $G$ and $H$ actions with action-equivariant maps as morphisms (see \cite[Proposition~4.18]{MZ1}).
	They also show that any functor with certain properties comes from a unique actor, and so the composition of such functors gives the composition of actors.
	More explicitly, for given actors $h:G\curvearrowright H$ and $k:H\curvearrowright K$ with anchor maps $\rho_h$ and $\rho_k$ respectively, the composition actor $kh:G\curvearrowright K$ has anchor map $\rho_{kh}:=\rho_h\circ\rho_k$ and action given by $\gamma\cdot_{kh}x:=(\gamma\cdot_h \rho_k(x))\cdot_k x$.
	In particular, \cite[Proposition~4.21]{MZ1} implies that this composition is associative, and that there is a category of groupoids with actors as morphisms, where the identity actor of a groupoid is given by the left multiplication of the groupoid on itself.
	One may think of multiplying or acting with $\rho_{k}(x)$ as `multiplying by a unit'.
	
	\begin{lemma}\label{lem-actorMultOpen}
		Let $G$ and $H$ be \'etale groupoids and let $G\curvearrowright H$ be an actor.
		The multiplication map $G\baltimes{s}{\rho}H\to H$ is open.
		\begin{proof}
			It suffices to show that the multiplication of the actor is open on (nonempty) basic open subsets of the form $U\baltimes{s}{\rho} V$, where $U\subseteq G$ and $V\subseteq H$ are open bisections.
			Fix a pair $(\gamma,x)\in U\baltimes{s}{\rho}V$, and let $(x_\lambda)_\lambda\subseteq H$ be a net converging to $\gamma\cdot x$.
			It suffices to show that $x_\lambda$ is eventually contained in $U\cdot V:=\{\gamma\cdot x:(\gamma,x)\in U\baltimes{s}{\rho}V\}$ for large enough $\lambda$. 
			Since the range map in $G$ is open, each $r(U)$ is an open neighbourhood of $r(\gamma)=\rho(\gamma\cdot x)$.
			Since $\rho$ is continuous, the net $(\rho(x_\lambda))_\lambda$ converges to $\rho(\gamma\cdot x)=r(x)$, so $\rho(x_\lambda)$ belongs to $r(U)$ eventually.
			For each such $\lambda$ large enough, set $\gamma_\lambda:=r|_{U}^{-1}(\rho(x_\lambda))$, so that $\gamma_\lambda\in U_j$ is the unique element satisfying $r(\gamma_\lambda)=\rho(x_\lambda)$.
			Then $(r(\gamma_\lambda))_\lambda$ is a net converging to $r(\gamma)$, whereby $\gamma_\lambda\to \gamma$ in $U$ since $r$ restricts to a homeomorphism on $U$.
			It follows that $(\gamma_\lambda^{-1}\cdot x_\lambda)_\lambda$ is a net converging to $\gamma^{-1}\cdot(\gamma\cdot x)=x$, therefore $\gamma_\lambda^{-1}\cdot x_\lambda$ eventually belongs to $V$ for large enough $\lambda$.
			Finally, $x_\lambda=\gamma_\lambda\cdot (\gamma_\lambda^{-1}\cdot x_\lambda)\in U\cdot V$ for large enough $\lambda$, as required.
		\end{proof}
	\end{lemma}
	
	\begin{lemma}\label{lem-actorBisectionProp}
		Let $G$ and $H$ be \'etale groupoids and let $G\curvearrowright H$ be an actor.
		For open bisections $U\subseteq G$ and $V\subseteq H$ the set $U\cdot V=\{\gamma\cdot x:\gamma\in U,x\in V, s(\gamma)=\rho(x)\}$ is an open bisection of $H$.
		In particular, for each $y\in U\cdot V$ there is a unique pair $(\gamma,x)\in U\baltimes{s}{\rho}V$ such that $y=\gamma\cdot x$.
		\begin{proof}
			We note first that $U\cdot V=(U\cdot r(V))\cdot V$, so it suffices to show that $U\cdot r(V)$ is an open bisection.
			
			First we show that the range and source maps in $H$ are injective on $U\cdot r(V)$.
			Let be elements $x_1,x_2\in U\cdot r(V)$.
			By definition there exist $\gamma_1,\gamma_2\in U$ and $y_1,y_2\in r(V)$ such that $x_i=\gamma_i\cdot y_i$.
			Suppose first that $s(x_1)=s(x_2)$.
			Then $y_1=s(\gamma_1\cdot y_1)=s(x_1)=s(x_2)=s(\gamma_2\cdot y_2)=y_2$, and $s(\gamma_1)=\rho(y_1)=\rho(y_2)=s(\gamma_2)$.
			Since $U$ is a bisection, we have $\gamma_1=\gamma_2$ whereby $x_1=\gamma_1\cdot y_1=\gamma_2\cdot y_2=x_2$.
			Now suppose $r(x_1)=r(x_2)$.
			Then $r(\gamma_1)=\rho(r(\gamma_1\cdot y_1))=\rho(r(\gamma_2\cdot y_2))=r(\gamma_2)$, and so $\gamma_1=\gamma_2$ since $U$ is a bisection.
			Then $(\gamma_1\cdot y_1)(\gamma_1^{-1}\cdot r(\gamma_1\cdot y_1))=r(\gamma_1\cdot y_1)=r(\gamma_1\cdot y_2)=(\gamma_1\cdot y_2)(\gamma_1^{-1}\cdot r(\gamma_1)\cdot y_1))$, and acting on the left by $\gamma_1^{-1}$ then yields $y_1(\gamma_1^{-1}\cdot r(\gamma_1\cdot y_1))=y_2(\gamma_1^{-1}\cdot r(\gamma_1)\cdot y_1))$.
			Taking the range of each side of the equation, we gain $y_1=y_2$ as required, so $U\cdot r(V)$ is a bisection.
			
			The set $U\cdot V$ is open by Lemma~\ref{lem-actorMultOpen}.
			Suppose two pairs $(\gamma_1,x_1),(\gamma_2,x_2)\in U\baltimes{s}{\rho}V$ satisfy $\gamma_1\cdot x_1=\gamma_2\cdot x_2$.
			Then $r(\gamma_1)=\rho(\gamma_1\cdot x_1)=\rho(\gamma_2\cdot x_2)=r(\gamma_2)$, whereby $\gamma_1=\gamma_2$ since $U$ is a bisection.
			Similarly, $s(x_1)=s(\gamma_1\cdot x_1)=s(\gamma_2\cdot x_2)=s(x_2)$, implying $x_1=x_2$ as $V$ is a bisection.
			This proves the final claim.
		\end{proof}
	\end{lemma}
	
	\begin{corollary}\label{cor-multRestrictsToHomeo}
		Let $G$ and $H$ be \'etale groupoids and let $G\curvearrowright H$ be an actor.
		For open bisections $U\subseteq G$ and $V\subseteq H$ the multiplication associated to the actor restricts to a homeomorphism $U\baltimes{s}{\rho}V\to U\cdot V$.
		\begin{proof}
			The subset $U\cdot V$ is a bisection by Lemma~\ref{lem-actorBisectionProp} so for any $y\in U\cdot V$ there are unique $\gamma\in U$ and $x\in V$ such that $\gamma\cdot x=y$.
			Hence the multiplication map of the actor is bijective on $U\baltimes{s}{\rho}V$, and is open by Lemma~\ref{lem-actorMultOpen}.
			Hence the multiplication map restricts to a homeomorphism $U\baltimes{s}{\rho}V\to U\cdot V$.
		\end{proof}
	\end{corollary}
	
	\subsection{Inducing morphisms between $C^*$-algebras}
	
	Buneci and Stachura \cite{BS1} construct what they call \emph{algebraic morphisms} between locally compact Hausdorff groupoids with Haar systems, and show that these induce ${}^*$-homomorphisms between $C^*$-algebras.
	While we wish to consider \'etale groupoids (which always carry a convenient Haar system of counting measures), we do not wish to restrict to strictly Hausdorff groupoids.
	Thus we shall reconstruct the results of Buneci and Stachura in the \'etale groupoid setting for groupoids with locally compact Hausdorff unit space (but not necessarily globally Hausdorff).
	
	Consider a continuous function $f:Y\to X$ between topological spaces.
	By Gelfand duality, the pullback along $f$ induces a ${}^*$-homomorphism $f^*:C_0(X)\to C_b(Y)=M(C_0(Y))$, and this ${}^*$-homomorphism has image contained in $C_0(Y)$ precisely when the map $f$ is proper.
	Continuous maps are an important subclass of actors, and we wish to consider actors between groupoids that induce ${}^*$-homomorphisms between the groupoid $C^*$-algebras (and not taking values in the local multiplier algebras).
	This motivates the following definition.
	
	\begin{definition}\label{defn-properActor}
		An actor $h:G\curvearrowright H$ is \emph{proper} if the anchor map $\rho:H\to G^{(0)}$ restricts to a proper map $\rho|_{H^{(0)}}:H^{(0)}\to G^{(0)}$ on the unit space $H^{(0)}$ of $H$.
	\end{definition}

	The term `proper' has its obvious heritage from that of proper continuous maps, but there is also the notion of a proper groupoid action: one where the action map $G\baltimes{s}{\rho}X\to X\times X$, $(\gamma,x)\mapsto (\gamma\cdot x,x)$ is a proper map.
	We would like to emphasise that a proper actor and a proper action are inequivalent notions.
	For example, the trivial actor $\{e\}\curvearrowright X$ on a topological space $X$ (viewed as a groupoid) is always a proper action, but it is only a proper actor if $X$ is compact.
	Throughout this article we adopt the usage of `proper' in the sense of Definition~\ref{defn-properActor}.
	
	\begin{proposition}\label{prop-actorsGiveStarHoms}
		Let $G$ and $H$ be \'etale groupoids with locally compact Hausdorff unit spaces.
		Let $h:G\curvearrowright H$ be a proper actor.
		There is a ${}^*$-homomorphism $\varphi_h:\Cc_c(G)\to\Cc_c(H)$ given by	
		$$[\varphi_h(f)](x)=\sum_{\substack{\gamma\in G_{\rho(s(x))}\\ x=\gamma\cdot s(x)}}f(\gamma).$$
		In particular if $f\in\Cc_c(U)$ for a bisection $U\subseteq G$ then $\varphi_h(f)\in\Cc_c(U\cdot H^{(0)})$, and for $\gamma\in U$, $x\in H^{(0)}$ with $s(\gamma)=\rho(x)$, the function $\varphi_h(f)$ satisfies
		$$[\varphi_h(f)](\gamma\cdot x)=f(\gamma).$$
		\begin{proof}
			First note that if $f\in\Cc_c(U)$ for some bisection $U\subseteq G$ then $[\varphi_h(f)](x)$ is non-zero only if $x=\gamma_x\cdot s(x)$ for some $\gamma_x\in U$ with $s(\gamma_x)=\rho(s(x))$, and so $\gamma_x$ is necessarily unique in $U$ as $U$ is a bisection.
			Thus
			$$[\varphi_h(f)](x)=\sum_{\substack{\gamma\in G_{\rho(s(x))}\\ x=\gamma\cdot s(x)}}f(\gamma)=f(\gamma_x).$$
			To see the function $\varphi_h(f)$ is continuous on $U\cdot H^{(0)}$, first recall that Lemma~\ref{lem-actorBisectionProp} implies $U\cdot H^{(0)}$ is an open bisection in $H$, and so the source map in $H$ restricts to a homeomorphism $U\cdot H^{(0)}\to s(U\cdot H^{(0)})=\rho^{-1}(s(U))\cap H^{(0)}$.
			Since $U$ is an open bisection, the source map in $G$ restricted to $U$ is a homeomorphism onto its image and so has an inverse $s|_U^{-1}:s(U)\to U$.
			We then see that on $U\cdot H^{(0)}$, the function $\varphi_h(f)$ is given by the composition of maps
			$$U\cdot H^{(0)}\xrightarrow{s}\rho^{-1}(s(U))\cap H^{(0)}\xrightarrow{\rho}s(U)\xrightarrow{s|_U^{-1}} U\xrightarrow{f}\CC,$$
			all of which are continuous on the relevant domains, and so $\varphi_h(f)$ is continuous on $U\cdot H^{(0)}$.
			To see that $\varphi_h(f)$ has compact support, we note that $[\varphi_h(f)](\gamma\cdot x)$ is non-zero for some $\gamma\cdot x\in U\cdot H^{(0)}$ if and only if $\gamma$ lies in the support of $f$.
			This implies that $x=s(\gamma)\cdot x$ belongs to $s(\supp(f))\cdot H^{(0)}$, which is exactly $\rho^{-1}(s(\supp(f)))\cap H^{(0)}$.
			This is compact since $s$ is continuous and the restriction of $\rho$ to $H^{(0)}$ is proper.
			
			The assignment $f\mapsto\varphi_h(f)$ is clearly a linear map $\Cc_c(G)\to\Cc_c(H)$.
			To see it is multiplicative, fix $f_1\in\Cc_c(U_1)$ and $f_2\in\Cc_c(U_2)$ where $U_i\subseteq G$ are open bisections.
			Then $\varphi_h(f_1f_2)$ and $\varphi_h(f_1)\varphi_h(f_2)$ both have support contained $(U_1U_2)\cdot H^{(0)}=(U_1\cdot H^{(0)})(U_2\cdot H^{(0)})$, and for each $(\gamma_1\gamma_2)\cdot x\in (U_1U_2)\cdot H^{(0)}$ (where $\gamma_i\in U_i$ and $x\in H^{(0)}$) we have
			\begin{align*}
				[\varphi_h(f_1f_2)]((\gamma_1\gamma_2)\cdot x)&=(f_1f_2)(\gamma_1\gamma_2)\\
				&=f_1(\gamma_1)f_2(\gamma_2)\\
				&=[\varphi_h(f_1)](\gamma_1\cdot r(\gamma_2\cdot x))[\varphi_h(f_2)](\gamma_2\cdot x)\\
				&=[\varphi_h(f_1)\varphi_h(f_2)]((\gamma_1\gamma_2)\cdot x).
			\end{align*}
			Thus $\varphi_h$ is multiplicative.
			For $x\in H^{(0)}$ and $\gamma\in G$ with $s(\gamma)=\rho(x)$ we have $r(\gamma)=\rho(r(\gamma\cdot x))$ so $(\gamma^{-1},r(\gamma\cdot x))$ is a composable pair in the actor.
			Moreover we see $(\gamma^{-1}\cdot r(\gamma\cdot x))(\gamma\cdot x)=(\gamma^{-1}\gamma)\cdot x=x$, so $\gamma^{-1}\cdot r(\gamma\cdot x)$ is the inverse of $\gamma\cdot x$.
			We also note that if $x=\gamma\cdot s(x)$ for some $x\in H$ and $\gamma\in G$, then $r(\gamma)=\rho(\gamma\cdot s(x))=\rho(x)=\rho(r(x))$.
			Thus for $f\in\Cc_c(G)$ and $x\in H$, we see that
			\begin{align*}
				[\varphi_h(f^*)](x)&=\sum_{\substack{\gamma\in G_{\rho(s(x))}\\ x=\gamma\cdot s(x)}}\overline{f(\gamma^{-1})}\\
				&=\sum_{\substack{\gamma\in G^{\rho(r(x))}_{\rho(s(x))}\\ x^{-1}=\gamma^{-1}\cdot s(x^{-1})}}\overline{f(\gamma^{-1})}\\
				&=\sum_{\substack{\eta\in G^{\rho(s(x))}_{\rho(r(x))}\\ x^{-1}=\eta\cdot s(x^{-1})}}\overline{f(\eta)}\\
				&=\sum_{\substack{\eta\in G_{\rho(s(x^{-1}))}\\ x^{-1}=\eta\cdot s(x^{-1})}}\overline{f(\eta)}\\
				&=\overline{[\varphi_h(f)](x^{-1})}=[\varphi_h(f)]^*(x).
			\end{align*}
			Thus $\varphi_h$ is a ${}^*$-homomorphism.
		\end{proof}
	\end{proposition}
	
	\begin{corollary}\label{cor-actorMapExtendsToFullAlg}
		The map $\varphi_h$ extends to a ${}^*$-homomorphism $C^*(G)\to C^*(H)$.
		\begin{proof}
			The norm $||f||_{\varphi_h}:=||\varphi_h(f)||$ is a $C^*$-norm on $\Cc_c(G)$.
			Hence it is bounded by the maximal $C^*$-norm in $C^*(G)$, whereby $\varphi_h$ extends continuously to $C^*(G)$.
		\end{proof}
	\end{corollary}
	
	\begin{corollary}\label{cor-actorHomPreservesSubalg}
		The image of $C_0(G^{(0)})$ under the ${}^*$-homomorphism $\varphi_h$ is contained in $C_0(H^{(0)})$.
		\begin{proof}
			Note that $G^{(0)}\subseteq G$ is an open bisection.
			By Proposition~\ref{prop-actorsGiveStarHoms} for $f\in C_0(G^{(0)})$ the open support of $\varphi_h(f)$ is contained in $G^{(0)}\cdot H^{(0)}=H^{(0)}$.
		\end{proof}
	\end{corollary}
	
	\begin{remark}
	If the actor $h$ is not proper, one may still construct a multiplier $\varphi_h(f)$ on $\Cc_c(H)$ given by
	$$[\varphi_h(f)g](x)=\sum_{\gamma\in G_{\rho(x)}}f(\gamma^{-1})g(\gamma\cdot x).$$
	This then extends to a multiplier of $C^*(H)$, and the assignment $f\mapsto\varphi_h(f)\in M(C^*(H))$ then extends to a ${}^*$-homomorphism $C^*(G)\to M(C^*(H))$.
	The subalgebras $C_0(G^{(0)})$ and $C_0(H^{(0)})$ embed non-degenerately into $C^*(G)$ and $C^*(H)$ respectively. 
	Between these two subalgebras, the ${}^*$-homomorphism maps to $C_0(H^{(0)})$ and not $C_b(H^{(0)})=M(C_0(H^{(0)}))$ if and only if the underlying Gelfand dual map is proper.
	The underlying Gelfand dual map is the restriction of the anchor map $\rho$ to the unit space $H^{(0)}$ in this case, and so taking values in $C^*(H)$ (rather than multipliers) arises directly as a result of properness of the actor.
	\end{remark}

	We shall now show that the assignment taking a groupoid to its $C^*$-algebra and an actor to its induced ${}^*$-homomorphism is a functor.
	
	\begin{lemma}\label{lem-actorToStarHomFunctorial}
		Let $h:G\curvearrowright H$ and $k:H\curvearrowright K$ be proper actors and let $\varphi_h:C^*(G)\to C^*(H)$ and $\varphi_k:C^*(H)\to C^*(K)$ be the respective ${}^*$-homomorphisms induced by $h$ and $k$ in Corollary~\textup{\ref{cor-actorMapExtendsToFullAlg}}.
		Let $kh:G\curvearrowright K$ be the composite actor.
		Then $\varphi_{kh}=\varphi_k\circ\varphi_h$.
		\begin{proof}
			Fix an open bisection $U\subseteq G$ and $f\in\Cc_c(U)$.
			For $\gamma\in U$ and $t\in K^{(0)}$ with $\rho_{kh}(t)=s(\gamma)$ we have $[\varphi_{kh}](f)(\gamma\cdot t)=f(\gamma)$ by Proposition~\ref{prop-actorsGiveStarHoms}.
			By the same Proposition, we know that $\varphi_{kh}(f)$ and $\varphi_k(\varphi_h(f))$ are both supported on the bisection $U\cdot_{hk}K^{(0)}=(U\cdot_h H^{(0)})\cdot_k K^{(0)}$.
			Observe then 
			$$[\varphi_k(\varphi_k(f))](\gamma\cdot t)=[\varphi_k(\varphi_h(f))]((\gamma\cdot 1_{\rho_k(t)})\cdot t)=[\varphi_h(f)](\gamma\cdot 1_{\rho_k(t)})=f(\gamma),$$
			whereby $\varphi_{hk}(f)=\varphi_h(\varphi_k(f))$.
			Both $\varphi_{hk}$ and $\varphi_h\circ\varphi_k$ are continuous and linear, so are equal as they coincide on a densely spanning subset of $C^*(G)$.
		\end{proof}
	\end{lemma}
	
	\subsection{Free actions and descension to quotients}
	The full groupoid $C^*$-algebra $C^*(G)$ has two quotients of particular interest: the reduced groupoid $C^*$-algebra $C^*_\red(G)$ and the essential groupoid $C^*$-algebra $C^*_\ess(G)$.
	For Hausdorff groupoids, these two quotients coincide, and for non-Hausdorff groupoids, we are mostly interested in the essential groupoid $C^*$-algebra (which itself is a quotient of the reduced).
	A natural question in the context of actors is then: what properties of an actor ensure that the induced ${}^*$-homomorphism descends to reduced or essential groupoid $C^*$-algebras?
	To answer this question, we must first recall the definition of the essential groupoid $C^*$-algebra.
	For this, we require some definitions.
	
	\begin{definition}[{\cite{AM1}}]
		Let $A$ be a $C^*$-algebra.
		An ideal $I\triangleleft A$ is \emph{essential} if $I^{\perp}:=\{a\in A: a\cdot b=0\text{ for all }b\in I\}$ is the zero ideal.
		The essential ideals of $A$ form a directed set ordered by containment (so that $I\leq J$ if and only if $I\supseteq J$).
		If $J\subseteq I$ is an essential ideal, then multipliers on $I$ have unique restrictions to multipliers on $J$, giving an injective ${}^*$-homomorphism $M(I)\to M(J)$.
		The \emph{local multiplier algebra} $\Mloc(A)$ of $A$ is the inductive limit of the system of essential ideals with restriction maps, and $A$ canonically embeds into $\Mloc(A)$ since $A$ embeds into $M(I)$ as multipliers for each essential ideal $I\subseteq A$.
		We will identify $A$ with its image in $\Mloc(A)$.
	\end{definition}
	
	\begin{definition}
		Let $X$ be a topological space.
		A subset $D\subseteq X$ is \emph{nowhere dense} if the closure of $D$ has empty interior (that is, $\bar{D}^\circ=\emptyset$).
		A subset $C\subseteq X$ is \emph{meagre} if it is a countable union of nowhere dense sets.
	\end{definition}
	
	\begin{definition}[{\cite[Definition~3.3]{KM3}}]
		Let $A\subseteq B$ be an inclusion of $C^*$-algebras, and let $\tilde{A}\supseteq A$ be a $C^*$-algebra containing $A$.
		A \emph{generalised conditional expectation} taking values in $\tilde{A}$ is a completely positive contractive linear map $E:B\to\tilde{A}$ restricting to the identity on $A$.
		If $\tilde{A}=\Mloc(A)$, the generalised expectation is called a \emph{local expectation}.
	\end{definition}
	
	By a result of Frank, Hamana's injective hull (see \cite{H1}) and the local multiplier algebra agree for commutative $C^*$-algebras \cite[Theorem~1]{F1}, so commutative inclusions $A\subseteq B$ always have a local multiplier algebra valued conditional expectation.
	By results of Dixmier \cite{D1} and later Gonshor \cite{G1}, the local multiplier algebra of $C_0(X)$ (and hence the injective hull of $C_0(X)$) is isomorphic to
	$$\Mloc(C_0(X))\cong\bB(X)/\mM(X),$$
	where $\bB(X)$ is the algebra of bounded Borel-measurable functions $X\to\CC$, and $\mM(X)$ is the ideal of such functions that vanish on a comeagre subset of $X$ (i.e. functions with `meagre support').
	This isomorphism maps the subalgebra $C_0(X)\subseteq\Mloc(C_0(X))$ identically to its image in the quotient by $\mM(X)$, where $C_0(X)$ embeds identically as a subspace of $\bB(X)$.
	The quotient map is of course injective on $C_0(X)$ in $\bB(X)$, as any non-zero continuous function on a Baire space has non-meagre support.
	
	Given a function $f\in\Cc_c(G)$, one may consider its restriction to the unit space $G^{(0)}$.
	This is a function on $G^{(0)}$, but need not be continuous.
	This restriction $f|_{G^{(0)}}$ will however be an element of $\bB(G^{(0)})$.
	Mapping $f$ to the equivalence class of its restriction $f|_{G^{(0)}}$ in $\bB(G^{(0)})/\mM(G^{(0)})=\Mloc(C_0(G^{(0)}))$ extends to a generalised conditional expectation $C^*(G)\to\Mloc(C_0(G^{(0)}))$.
	
	The essential groupoid $C^*$-algebra $C^*_\ess(G)$ is defined as the quotient of $C^*(G)$ by the largest ideal contained in $\ker(EL)$ (see \cite[Definition~7.12]{KM2}).
	By \cite[Proposition~4.3]{KM2}, the conditional expectation defining the reduced groupoid $C^*$-algebra dominates $EL$, implying that $C^*_\ess(G)$ is a quotient of $C^*_{\red}(G)$.
	Hence elements of $C^*_\ess(G)$ may be viewed as equivalence classes of Borel-measurable functions $G\to\CC$, by embedding $C^*_{\red}(G)$ in $\bB(G)$ using \cite[Proposition~7.10]{KM2}.
	If $G$ is covered by countably many bisections, then \cite[Proposition~7.18]{KM2} implies that a function $f\in C^*_\red(G)\subseteq\bB(G)$ represents the zero element of $C^*_\ess(G)$ if and only if its support is meagre.
	
	Since the essential groupoid $C^*$-algebra is defined as a quotient by an ideal contained in the kernel of a local expectation, it is enough to show that $\varphi_h$ entwines conditional expectations in an appropriate sense.
	Of course this is not the only way to ensure that a ${}^*$-homomorphism descends to a quotient, for instance, if $G$ is Hausdorff and satisfies weak containment (in particular, if $G$ is amenable) then $C^*(G)=C^*_\ess(G)$, so all ${}^*$-homomorphisms automatically descend independent of how they interact with conditional expectations.
	We shall however restrict our attention to the case where these ${}^*$-homomorphisms entwine conditional expectations, since we are interested in morphisms that preserve the structure of the inclusion $C_0(G^{(0)})\subseteq C^*(G)$.
	
	Since we do not assume all groupoids are Hausdorff in this article, we need to refine our notion of `entwining conditional expectations', since the expression $\varphi_h\circ EL$ is not defined if $EL:C^*(G)\to\Mloc(C_0(G^{(0)}))$ is not a genuine expectation.
	
	\begin{lemma}\label{lem-liftToLocMultAlg}
		Let $\varphi:A\to B$ be a ${}^*$-homomorphism.
		For each ideal $I\triangleleft A$, let $B_I:=\overline\sspan B\varphi(I)B$ be the ideal of $B$ generated by the image of $I$ in $B$.
		If $\varphi(I)$ contains an approximate unit for $B_I$ for all ideals $I\triangleleft A$, then $\varphi$ extends to a ${}^*$-homomorphism $\tilde\varphi:\Mloc(A)\to\Mloc(B)$.
		If $\varphi$ is injective, then so is $\tilde\varphi$.
		\begin{proof}
			The map $\varphi$ restricts to a non-degenerate ${}^*$-homomorphism $\varphi|_I:I\to B_I$ for each ideal $I\triangleleft A$ by assumption.
			Hence there is an extension $\tilde\varphi_I:M(I)\to M(B_I)$ satisfying $\tilde\varphi_I(\tau)ab=\varphi(\tau a)b$ for all $a\in I$ and $b\in B_I$.
			If $J\subseteq I$ is an essential ideal, the restriction map $M(I)\to M(J)$ is injective.
			For any $a\in J$ and $b\in B_J$ we have
			$$\tilde\varphi_J(\tau|_J)ab=\varphi(\tau a)b=\tilde\varphi(\tau)ab,$$
			whereby $\tilde\varphi_J(\tau|_J)=\tilde\varphi_I(\tau)|_J$, since $\varphi|_J:J\to B_J$ is non-degenerate.
			The system $(\tilde\varphi_I)_{I\triangleleft A}$ ranging over essential ideals $I\triangleleft A$ is compatible with the inductive limit structure of $\Mloc(A)$ and so extends to a ${}^*$-homomorphism $\tilde\varphi:\Mloc(A)\to\Mloc(B)$.
			
			Now suppose that $\varphi$ is injective. 
			The maps $\tilde\varphi_I$ are then also injective, hence isometric.
			Thus $\tilde\varphi$ is isometric on the dense subspace $\sspan M(I)\subseteq\Mloc(A)$, and hence is globally isometric.
		\end{proof}
	\end{lemma}
	
	\begin{corollary}\label{cor-alwaysLiftCommutative}
		Let $\varphi:A\to B$ be a ${}^*$-homomorphism and suppose that $B$ is commutative.
		Then $\varphi$ extends to a ${}^*$-homomorphism $\tilde\varphi:\Mloc(A)\to\Mloc(B)$, and $\tilde\varphi$ is injective if and only if $\varphi$ is.
		\begin{proof}
			Using Lemma~\ref{lem-liftToLocMultAlg} it suffices to show that $\varphi(I)$ contains an approximate identity for $B_I:=\overline\sspan B\varphi(I)B$.
			Note that $\varphi(I)$ contains an approximate unit for the hereditary subalgebra $\overline\sspan\varphi(I)B\varphi(I)$, which is an ideal in $B$ as $B$ is commutative.
			This ideal is the ideal generated by $\varphi(I)$, and so is equal to $B_I$.
		\end{proof}
	\end{corollary}

	\begin{remark}[Hooptedoodle]
		At first glance, the proof of Corollary~\ref{cor-alwaysLiftCommutative} uses specifically that all hereditary subalgebras of the codomain algebra $B$ are ideals, and no other properties of being a commutative $C^*$-algebra.
		This is however an equivalent characterisation of commutativity for $C^*$-algebras.
		In general there is a bijection between closed left-ideals of a $C^*$-algebra $B$ and its hereditary subalgebras, given by $L\mapsto L\cap L^*$.
		Since (two-sided) ideals are both left-ideals and are ${}^*$-closed, we see that any ideal $I$ is mapped to itself under this map.
		Thus, if all hereditary subalgebras are ideals, then all closed left-ideals are two-sided.
		Given a pure state $\psi$ on $B$, the set $\{b\in B:\psi(b^*b)=0\}$ is a maximal closed left-ideal, and so is a maximal ideal under our assumptions.
		The quotient $B/L$ is both a $C^*$-algebra and a division ring, whereby $B/L$ is isomorphic to $\CC$.
		Hence all irreducible representations of $B$ are one-dimensional, whereby $B$ is commutative.
	\end{remark}
	
	\begin{definition}
		Let $G$ and $H$ be \'etale groupoids and let $\varphi:C^*(G)\to C^*(H)$ be a ${}^*$-homomorphism such that $\varphi(C_0(G^{(0)}))\subseteq C_0(H^{(0)})$.
		Let $\tilde\varphi|_{C_0(G^{(0)})}:\Mloc(C_0(G^{(0)}))\to\Mloc(C_0(H^{(0)}))$ be the extension of $\varphi|_{C_0(G^{(0)})}$ to a map between local multiplier algebras given by Lemma~\ref{lem-liftToLocMultAlg}.
		Let $EL_G$ and $EL_H$ be the canonical local expectations associated to $G$ and $H$ respectively.
		We say that $\varphi$ \emph{entwines expectations} if $EL_H\circ\varphi=\tilde\varphi|_{C_0(G^{(0)})}\circ EL_G$.
	\end{definition}	
	
	We wish to determine which actors induce ${}^*$-homomorphisms between groupoid $C^*$-algebras that entwine the (local) conditional expectations.
	For an actor $G\curvearrowright H$ and a function $f\in\Cc_c(G)$, Proposition~\ref{prop-actorsGiveStarHoms} gives a formula for the image of the function under the induced ${}^*$-homomorphism $\varphi:\Cc_c(G)\to \Cc_c(H)$.
	For points $t\in H^{(0)}$, this formula reads
	$$[\varphi(f)](t)=\sum_{\gamma\cdot t=t}f(\gamma).$$
	One such summand is $f(\rho(t))$, since $\rho(t)\cdot t=t$ by definition.
	If $\varphi(f)$ entwines expectations, we see that this should be the only such possibly non-zero summand, since $EL(f)$ will have support contained in $G^{(0)}$.
	Hence one could expect that the collection of arrows $\gamma$ in $G$ satisfying $\gamma\cdot t=t$ should contain only the unit $\rho(t)$.
	This motivates the following definition and deductions thereafter.
	
	\begin{definition}
		Let $h:G\curvearrowright H$ be an actor with anchor map $\rho$.
		For $t\in H^{(0)}$, we say the actor $h$ is \emph{free at $t$} if for all $\gamma\in G_{\rho(t)}$ we have $\gamma\cdot t=t$ if and only if $\gamma=s(\gamma)$.
		We say $h$ is \emph{free} if $h$ is free at all $t\in H^{(0)}$.
	\end{definition}
	
	\begin{lemma}
		Let $h:G_1\curvearrowright G_2$ be an actor.
		If $h$ is free and the anchor map $\rho_h:G_2\to G_1^{(0)}$ is surjective, then $h$ is a monomorphism (that is, left-cancellative).
		\begin{proof}
			Suppose that $h$ is free and $\rho_h$ is surjective, and fix actors $k_i:H\curvearrowright G_1$  with anchor maps $\rho_i:H\to G_1$ for $i=1,2$ such that $hk_1=hk_2$.
			For each $t\in G_2^{(0)}$ we have $\rho_1\circ\rho_h(t)=\rho_2\circ\rho_h(t)$, whereby $\rho_1=\rho_2$ as $\rho_h$ is surjective.
			For $\eta\in H$ with $s(\eta)=\rho_1(\rho_h(t))$ we have
			$$(\eta\cdot_{k_1}\rho_h(t))\cdot_h t=\eta\cdot_{hk_1}t=\eta\cdot_{hk_2}t=(\eta\cdot_{k_2}\rho_h(t))\cdot_h t,$$
			whereby $\eta\cdot_{k_1}\rho_h(t)=\eta\cdot_{k_2}\rho_h(t)$ as $h$ is free.
			Every unit in $G_1$ is of the form $\rho_h(t)$ for some $t\in G_2^{(0)}$ as $\rho_h$ is surjective, and actors are uniquely determined by how they act on the unit space of the target groupoid, so $k_1=k_2$ holds.
			Thus $h$ is left-cancellative.
		\end{proof}
	\end{lemma}
	
	For an actor $h:G\curvearrowright H$, $t\in H^{(0)}$, $x\in H^t$, and $\gamma\in G_{\rho(t)}$, it is clear that $\gamma\cdot x=x$ if and only if $\gamma\cdot t=t$.
	For a bisection $U\subseteq G$, define the set $S^h_U\subseteq H^{(0)}$ by
	\begin{align*}
		S^h_U&:=\{t\in H^{(0)}:\exists \gamma\in U\setminus G^{(0)},s(\gamma)=\rho(t),\gamma\cdot t=t\}\\
		&=(U\setminus G^{(0)})\cdot_h H^{(0)}\cap H^{(0)}.
	\end{align*}
	We may write $S_U$ for $S^h_U$ if the actor $h$ is clear from context.
	A point $t\in H^{(0)}$ belongs to $S_U$ for some bisection $U\subseteq G$ exactly when the actor is not free at $t$.
	We call $S_U$ the \emph{stable space} of $U$.
	If $G$ is Hausdorff then $S_U$ is an open set, but this does not generally hold for non-Hausdorff groupoids.
	
	\begin{proposition}\label{prop-anctorHomEntwines}
		Let $h:G\curvearrowright H$ be a proper actor with anchor map $\rho:H\to G^{(0)}$.
		Let $\varphi_h:C^*(G)\to C^*(H)$ be the induced ${}^*$-homomorphism, and let $\tilde\varphi_h:\Mloc(C_0(G^{(0)}))\to \Mloc(C_0(H^{(0)}))$ be the ${}^*$-homomorphism from Corollary~\textup{\ref{cor-alwaysLiftCommutative}} extending the restriction $\varphi_h|_{C_0(G^{(0)})}:C_0(G^{(0)})\to C_0(H^{(0)})$.
		The map $\varphi_h$ entwines the canonical local expectations $EL_G$ and $EL_H$ if and only if for each open bisection $U\subseteq G$, the stable space $S_U$ of $U$ has meagre intersection with the preimage under $\rho$ of any precompact open set of $G^{(0)}$.
		\begin{proof}
			On the image of $\Cc_c(G)$ in $C^*(G)$ under the canonical embedding, the local expectation $EL_G$ acts by restricting a function in $\Cc_c(G)$ to the unit space $G^{(0)}$ and taking the class of this restriction in the quotient by meagrely supported functions.
			Since $EL_G$ and $EL_H$ are continuous and linear, it suffices to consider functions $f\in\Cc_c(U)$ for arbitrary open bisections $U\subseteq G$.
			Suppose $\varphi_h$ entwines conditional expectations.
			Then for each $f\in\Cc_c(U)$ and a comeagre set $C_f\subseteq H^{(0)}$ we have 
			$$f(\rho(t))=[\tilde\varphi_h(EL_G(f))](t)=[EL_H(\varphi_h(f))](t)=\begin{cases}
			f(\gamma_t),&\exists\gamma\in U,s(\gamma)=\rho(t),\gamma\cdot t=t\\
			0,&\text{otherwise},
			\end{cases}$$
			for all $t\in C_f$.
			Thus, for each $t\in C_f$, it follows that $f(\rho(t))\neq 0$ if and only if $\rho(t)\in U$, so the set of points where this is not satisfied is contained in $H^{(0)}\setminus C_f$, which is meagre.
			Therefore the set
			$$S_f:=\{t\in H^{(0)}:\exists\gamma\in U\setminus G^{(0)},s(\gamma)=\rho(t),\gamma\cdot t=t,f(\gamma)\neq 0\}\subseteq S_U$$ 
			is contained in the complement of $C_f$ in $H^{(0)}$, and so is meagre.
			Let $\supp^\circ(f)$ be the open support of $f$ in $U$.
			Let $V:=\rho^{-1}(s(\supp^\circ(f)))$ be the preimage of the source of the support of $f$, which is open and precompact since $\rho$ is proper and continuous.
			The intersection of $V$ with $S_U$ consists of exactly the points $t\in H^{(0)}$, where there is an element $\gamma_t\in U\setminus G^{(0)}$ with $s(\gamma_t)=\rho(t)$, and $f(\gamma)\neq 0$.
			Thus $\rho^{-1}(s(\supp^\circ(f)))\cap S_U=S_f$.
			Every precompact open subset of $s(U)$ can be written as the source of the open support of a function in $\Cc_c(U)$, so we see that $S_U$ has meagre intersection with the preimage of every such set, as required.\\
			
			Now suppose that $S_U$ has meagre intersection with the preimge under $\rho$ of every precompact open subset in $H^{(0)}$.
			For a function $f\in\Cc_c(U)$, the source $V_f:=s(\supp^\circ(f))$ of the open support of $f$ is a precompact open subset of $G^{(0)}$.
			Since $\rho$ is proper, the preimage $\rho^{-1}(V_f)$ is also precompact.
			Note now that a point $t\in H^{(0)}$ that belongs to the set $S_f$ defined above will satisfy $\rho(t)\in V_f$, and will also be a subset of $S_U$ as noted above.
			Thus $S_f\subseteq \rho^{-1}(V_f)\cap S_U$, which is meagre by hypothesis.
			As noted above, $S_f$ contains all points in $H^{(0)}$ that do not satisfy $[\varphi_h(f)](t)= f(\rho(t))$.
			Since this set is meagre, the complement $C_f:=H^{(0)}\setminus S_f$ is comeagre, and all points $t\in C_f$ satisfy $[\varphi_h(f)](t)= f(\rho(t))$.
			Hence $\tilde\varphi_h(EL_G(f))=EL_H(\varphi_h(f))$, where we identify $\Mloc(C_0(H^{(0)}))$ with $\bB(H^{(0)})/\mM(H^{(0)})$.
		\end{proof}
	\end{proposition}
	
	\begin{corollary}\label{cor-freeActEntwines}
		Let $h:G\curvearrowright H$ be a free and proper actor.
		Then $\varphi_h:C^*(G)\to C^*(H)$ entwines conditional expectations and descends to a ${}^*$-homomorphism of essential groupoid $C^*$-algebras $C^*_\ess(G)\to C^*_\ess(H)$.
		\begin{proof}
			If $h$ is free then $S_U$ is empty for all bisections $U\subseteq G$.
			By Proposition~\ref{prop-anctorHomEntwines}, the resulting $\varphi_h$ entwines canonical local expectations, and so $\varphi_h$ descends to a map between the essential groupoid $C^*$-algebras.
		\end{proof}
	\end{corollary}

	Consider a unit $t\in G^{(0)}$ lying in the image of the anchor map of a free actor.
	The source fibre over $t$ is then represented faithfully on the $\ell^2$ spaces of the range fibres over units in $H$ that map to $t$ under the anchor map.
	However, if $t$ does not lie in the image of the anchor map, then the source fibre over $t$ does not act in the actor at all, and so is not detected by the induced ${}^*$-homomorphism.
	This faithfulness over the image of the anchor map detects the faithfulness over the induced ${}^*$-homomorphism.
	
	\begin{proposition}\label{prop-cartanMapInjIffInjOnSubalg}
		Let $h:G_1\curvearrowright G_2$ be a free and proper actor.
		The induced ${}^*$-homomorphism $\varphi_h:C^*_\ess(G_1)\to C^*_\ess(G_2)$ is injective if it is injective on $C_0(G_1^{(0)})$.
		That is, the anchor map $\rho:G_2\to G_1^{(0)}$ is surjective.
		\begin{proof}
			Identify $C^*_\red(G_1)$ as a subalgebra of $\bB(G_1)$ by \cite[Proposition~7.10]{KM2}.
			Fix $f\in C^*_\ess(G_1)$ and let $\tilde{f}\in C^*_\red(G_1)$ be a function representing $f$ in the quotient.
			Let $EL_i:C^*_\ess(G_i)\to \Mloc(C_0(G_i^{(0)}))$ be the canonical faithful local expectations associated to $C^*_\ess(G_i)$ for $i=1,2$.
			Let $\tilde\varphi_h:\Mloc(C_0(G_1^{(0)}))\to\Mloc(C_0(G_2^{(0)}))$ be the extension of $\varphi_h|_{C_0(G_1^{(0)})}:C_0(G_1^{(0)})\to C_0(G_2^{(0)})$ from Lemma~\ref{lem-liftToLocMultAlg} and Corollary~\ref{cor-alwaysLiftCommutative}.
			Note that this is injective since $\rho$ is surjective.
			Since $EL_2$ is faithful, we have $\varphi_h(f)=0$ if and only if $EL_2(\varphi_h(f^*f))=0$.
			The map $\varphi_h$ entwines conditional expectations by Corollary~\ref{cor-freeActEntwines}, and so we see $0=EL_2(\varphi_h(f^*f))=\tilde\varphi_h(EL_1(f^*f))$, whereby $EL_1(f^*f)\in\ker(\tilde\varphi_h)$.
		\end{proof}
	\end{proposition}	
	
	Sometimes actors are automatically free.
	
	\begin{lemma}
		Let $h:G\curvearrowright H$ be an actor with anchor map $\rho$.
		Suppose that $G$ is effective and the anchor map $\rho$ is open.
		Then $h$ is free.
		\begin{proof}
			Fix $t\in H^{(0)}$ and $\gamma\in G$ with $s(\gamma)=\rho(t)$ and $\gamma\cdot t=t$.
			Then $r(\gamma)=\rho(\gamma\cdot t)=\rho(t)=s(\gamma)$, whereby $\gamma$ lies in the isotropy $G'$ of $G$.
			Let $U\subseteq G$ be an open bisection neighbourhood of $\gamma$.
			Let $m:G\baltimes{s}{\rho}H\to H$ be the multiplication map of the actor $h$, and let $V\subseteq H^{(0)}$ be an open neighbourhood of $t$.
			The preimage $m^{-1}(V)$ is then a neighbourhood of $(\gamma,t)$ in $G\baltimes{s}{\rho}H$, and so $W:=U\baltimes{s}{\rho}H^{(0)}\cap m^{-1}(V)$ is as well.
			For every $(\eta,t')\in W$ we have $\eta\cdot t'=m(\eta,t')\in V\subseteq H^{(0)}$, whereby $\eta\cdot t'=t'$ and $\eta$ must also belong to the isotropy of $G$.
			The projection $p:G\baltimes{s}{\rho}H\to G$ is an open map since both $s$ and $\rho$ are open, and so $p(W)$ is an open neighbourhood of $\gamma$ contained in isotropy.
			Thus $\gamma$ belongs to the interior of the isotropy of $G$, which is the unit space since $G$ is effective.
		\end{proof}
	\end{lemma}	
	
	Proposition~\ref{prop-anctorHomEntwines} shows that the preimages of stable spaces $S_U$ of an actor $h$ must be somewhat `small' in order for the induced ${}^*$-homomorphism to entwine conditional expectations.
	If the sets $S_U$ are a priori `small' (for example, if they are meagre), then this becomes a condition on the anchor map preserving this notion of size.
	In particular, we wish to consider maps which map only meagre sets to meagre sets.
	
	\begin{definition}
		A map  $f:X\to Y$ between topological spaces is \emph{skeletal} if the preimage of any nowhere dense set is again nowhere dense.
	\end{definition}	
	
	\begin{lemma}
		The preimage of a meagre set under a skeletal map is meagre.
		\begin{proof}
			Taking preimages preserves all unions, in particular countable unions.
			Thus the preimage of a union of countably many nowhere dense sets is a countable union of nowhere dense sets; it is meagre.
		\end{proof}
	\end{lemma}
	
	Skeletal maps are a class of functions which only map meagre sets to meagre sets.
	This will allow us to give a more concise and pragmatic version of Proposition~\ref{prop-anctorHomEntwines} if the anchor map associated to an actor between groupoids is skeletal.
	
	\begin{lemma}\label{lem-skeletalMeagreDiffProperty}
		Let $h:G\curvearrowright H$ be a proper actor with anchor map $\rho:H\to G^{(0)}$.
		Suppose that the restriction of $\rho$ to $H^{(0)}$ is skeletal.
		For any bisection $U\subseteq G$, the set $(U\setminus \overline{G^{(0)}})\cdot_h H^{(0)}$ is open and dense in $(U\setminus G^{(0)})\cdot_h H^{(0)}$.
		In particular, it is comeagre.
		\begin{proof}
			Let $m:G\baltimes{s}{\rho}H\to H$ be the multiplication associated to the actor $h$.
			Suppose that $V\subseteq H$ is an open subset such that $V\cap (U\setminus G^{(0)})\cdot H^{(0)}\neq\emptyset$.
			Then the preimage of this set under $m$ is also non-empty, since $m$ is surjective, and so
			$$\emptyset\neq m^{-1}(V\cap (U\setminus G^{(0)})\cdot H^{(0)})=m^{-1}(V)\cap[(U\setminus G^{(0)})\baltimes{s}{\rho}H^{(0)}].$$
			The set $U\setminus\overline{G^{(0)}}$ is dense in $U\setminus G^{(0)}$, and so $(U\setminus\overline{G^{(0)}})\baltimes{s}{\rho}H^{(0)}$ is dense in $(U\setminus G^{(0)})\baltimes{s}{\rho}H^{(0)}$ in the product topology.
			The preimage $m^{-1}(V)$ is open since $m$ is continuous, and so $m^{-1}(V)\cap (U\setminus\overline{G^{(0)}})\baltimes{s}{\rho}H^{(0)}$ is nonempty.
			Thus $V\cap (U\setminus\overline{G^{(0)}})\cdot H^{(0)}\neq\emptyset$, whereby $(U\setminus\overline{G^{(0)}})\cdot H^{(0)}$ is dense in $(U\setminus G^{(0)})\cdot H^{(0)}$.
			
			Lastly note that $U\setminus\overline{G^{(0)}}$ is an open bisection in $G$, so Lemma~\ref{lem-actorBisectionProp} implies that $(U\setminus\overline{G^{(0)}})\cdot H^{(0)}$ is open.
		\end{proof}
	\end{lemma}	
	
	\begin{proposition}\label{prop-starHomEntwinesSkeletal}
		Let $h:G\curvearrowright H$ be a proper actor with anchor map $\rho:H\to G^{(0)}$.
		Suppose that the resctriction of $\rho$ to $H^{(0)}$ is skeletal.
		For an open bisection $U\subseteq G$, the stable set $S_U$ has meagre intersection with every precompact open subset of $H^{(0)}$ if and only if it is empty.
		\begin{proof}
			Lemma~\ref{lem-skeletalMeagreDiffProperty} implies that the set $S_U':=(U\setminus\overline{G^{(0)}})\cdot H^{(0)}\cap H^{(0)}$ is a dense open subset of $S_U=(U\setminus G^{(0)})\cdot H^{(0)}\cap H^{(0)}$.
			$S_U'$ is also an open bisection of $H$ by Lemma~\ref{lem-actorBisectionProp}, since $U\setminus \overline{G^{(0}}$ is an open bisection in $G$.
			For any precompact open $V\subseteq H^{(0)}$, the set $S_U\cap V$ is meagre if and only if $S_U'\cap V$ is.
			However, $S_U'\cap V$ is an open subset of a Baire space, and so can only be meagre if it is empty.
			Thus $S_U'\cap V$ is empty for any precompact open $V\subseteq H^{(0)}$.
			Since $H$ is locally compact, there exists a basis of such precompact open sets, and so $S_U'$ is empty.
			Thus $S_U$ is empty, as it is contained in the closure of $S_U'=\emptyset$.
		\end{proof}
	\end{proposition}
	
	\begin{corollary}
		Let $h:G\curvearrowright H$ be a proper actor with anchor map $\rho:H\to G^{(0)}$.
		Suppose that the resctriction of $\rho$ to $H^{(0)}$ is skeletal.
		Then the induced ${}^*$-homomorphism $\varphi_h:C^*(G)\to C^*(H)$ entwines canonical local expectations if and only if the stabiliser set $S_U$ is empty for all bisections $U\subseteq G$, if and only $h$ is free.
		\begin{proof}
			We note that $h$ is free precisely when $S_U$ is empty for all bisections $U\subseteq G$.
			Applying Propositions~\ref{prop-anctorHomEntwines} and \ref{prop-starHomEntwinesSkeletal} gives the desired results.
		\end{proof}
	\end{corollary}
	
	If $G$ is Hausdorff, we do not require that the anchor map is skeletal.
	
	\begin{corollary}
		Let $h:G\curvearrowright H$ be a proper actor with anchor map $\rho:H\to G^{(0)}$.
		Suppose $G$ is Hausdorff.
		For an open bisection $U\subseteq G$, the stable set $S_U$ has meagre intersection with every precompact open subset of $H^{(0)}$ if and only if it is empty.
		In particular the induced ${}^*$-homomorphism $\varphi_h:C^*(G)\to C^*(H)$ entwines canonical local expectations if and only if the stabiliser set $S_U$ is empty for all bisections $U\subseteq G$, if and only $h$ is free.
		\begin{proof}
			The unit space $G^{(0)}$ of $G$ is clopen since $G$ is \'etale and Hausdorff.
			Hence the sets $S_U$ and $S_U'$ in the proof of Proposition~\ref{prop-starHomEntwinesSkeletal} coincide, as $G^{(0)}=\overline{G^{(0)}}$.
			The same proof shows that $S_U$ is empty if and only if its intersection with every precompact open subset of $H^{(0)}$ is meagre.
		\end{proof}
	\end{corollary}
	
	\begin{example}
		The class of actors with skeletal anchor maps is not without examples.
		Open maps are in particular skeletal, so if the anchor map $\rho$ associated to an actor is open, then the induced ${}^*$-homomorphism between groupoid $C^*$-algebras entwines canonical local expectations if and only if the actor is free.
		If $G$ is a group, then every map $H\to G^{(0)}=\{\text{pt}\}$ is skeletal, since the only meagre subset of a single point space is empty.
		Thus every actor of a group on a groupoid with compact unit space is skeletal and proper (this is however a moot point, since discrete groups are Hausdorff).
	\end{example}

	\section{Examples}
	
	Here we give some examples of actors arising from common groupoid constructions.
	We shall recall two definitions of transformation groupoid, one associated to group actions by homeomorphisms, which is elementary to state, and a more general version for inverse semigroup actions by partial homeomorphisms.
	The two definitions coincide up to isomorphism for group actions, so there is no ambiguity when referring to the `transformation groupoid' in general.
	
	\begin{definition}\label{defn-grpActionTrsfmGrpd}
		Let $\Gamma$ be a discrete group and suppose $\Gamma$ acts on a locally compact Hausdorff space $X$ by homeomorphisms.
		As a topological space, define $\Gamma\ltimes X:=\Gamma\times X$ with the product topology.
		There is a multiplication on $\Gamma\ltimes X$ given by
		$$(\gamma_1,x_1)\cdot(\gamma_2,x_2):=(\gamma_1\gamma_2,x_2),$$
		defined whenever $x_1=\gamma_2\cdot x_2$ turning it into a groupoid.
		Inversion is given by $(\gamma,x)^{-1}=(\gamma^{-1},\gamma\cdot x)$, and the range and source maps are given by
		$$r(\gamma,x)=(1_\Gamma,\gamma\cdot x),\qquad s(\gamma,x)=(1_\Gamma,x).$$
		Hence the unit space of $\Gamma\ltimes X$ is $\{1_\Gamma\}\times X$, which we identify with $X$.
		This groupoid is \'etale, since there is a basis of bisections given by $\{\gamma\}\times U$ for $\gamma\in\Gamma$ and $U\subseteq X$ open.
	\end{definition}
	
	\begin{example}\label{eg-transGrpdActor}
		Let $\Gamma$ be a discrete group and let $\Gamma\curvearrowright X,Y$ be actions of $\Gamma$ on locally compact Hausdorff topological spaces $X$ and $Y$.
		Let $\rho:X\to Y$ be a continuous $\Gamma$-equivariant map.
		We shall define an actor of the transformation groupoids $\Gamma\ltimes Y\curvearrowright\Gamma\ltimes X$ with anchor map $(g,x)\mapsto \rho(g\cdot x)$.
		We shall overload $\rho$ to denote this anchor map.
		For $((g_1,y),(g_2,x))\in (\Gamma\ltimes Y)\baltimes{s}{\rho}(\Gamma\ltimes X)$ define $(\gamma,y)\cdot(\eta,x):=(\gamma\eta,y)$.
		To see the multiplication of the actor is continuous, fix $V\subseteq X$ open, and $g\in\Gamma$.
		Sets of the form $\{g\}\times V$ form a basis for $\Gamma\ltimes X$, and the preimage of $\{g\}\times V$ under the multiplication is the union
		$$\bigcup_{g_1g_2=g} (\{g_1\}\times Y)\baltimes{s}{\rho}(\{g_2\}\times V),$$
		which is open in $(\Gamma\ltimes Y)\baltimes{s}{\rho}(\Gamma\ltimes X)$.
		The $\Gamma$-equivariance of $\rho$ ensures that the proposed actor satisfies condition (1) of Definition~\ref{defn-action} and associativity of the group $\Gamma$ gives condition (2) as well as compatibility with the multiplication in $\Gamma\ltimes X$.
		Condition (3) is satisfied since units in $\Gamma\ltimes Y$ act by the unit of $\Gamma$.
		
		The actor described in this example is free precisely when the action of $\Gamma$ on $X$ is, and the actor is proper exactly when $\rho$ is proper between $X$ and $Y$.
		
		If $Y$ is just a point so that $\Gamma\ltimes Y\cong\Gamma$, then the actor described above arises from the constant map $X\to \{\pt\}$, which is always $\Gamma$-equivariant.
		The resulting ${}^*$-homomorphism $C^*(\Gamma)\to M(C^*(\Gamma\ltimes X))\cong M(C_0(X)\rtimes\Gamma)$ then coincides with the ${}^*$-homomorphism induced by the canonical representation of $\Gamma$ in $\mathcal{U}M(C_0(X)\rtimes\Gamma)$.
	\end{example}

	\begin{example}\label{eg-fibrewiseBijectiveMapsExample}
		Let $G$ and $H$ be \'etale groupoids and let $\phi:H\to G$ be a continuous  homomorphism that restricts to a bijection $\phi|_{H_t}:H_t\to G_{\phi(t)}$ for each $t\in H^{(0)}$.
		Set $\rho:=\phi|_{H^{(0)}}\circ r:H\to G^{(0)}$ and define a map $\cdot:G\baltimes{s}{\rho}H\to H$ by
		$$(\gamma,x)\mapsto \gamma\cdot x:=\phi|_{H_{r(x)}}^{-1}(\gamma)x.$$
		First we note that if $(\gamma,x)$ belongs to $G\baltimes{s}{\rho}H$ then $s(\gamma)=\rho(x)=\phi(r(x))$, so $\gamma$ is an element of $G_{\phi(r(x))}=\phi(H_{r(x)})$, hence $\phi|_{H_{r(x)}}^{-1}(\gamma)$ is composable with $x$.
		
		The product $\cdot$ defines an actor of $G$ on $H$ with anchor map $\rho$.
		Indeed, the algebraic axioms follow from $\phi$ being a homomorphism, it remains only to show that $\cdot$ is continuous.
		Let $(\gamma_\lambda,x_\lambda)\subseteq G\baltimes{s}{\rho}H$ be a net converging to some $(\gamma,x)$.
		In particular, $\gamma_\lambda\to\gamma$ and $x_\lambda\to x$.
		Since multiplication in $H$ is continuous and $\cdot$ commutes with it on the right, it suffices to consider the case where $x=xx^{-1}$ is a unit in $H^{(0)}$, and we may assume that $x_\lambda$ is a unit for all $\lambda$ since $H^{(0)}$ is open.
		We must show that $\gamma_\lambda\cdot x_\lambda=\phi|_{H_{x_\lambda}}^{-1}(\gamma_\lambda)$ converges to $\gamma\cdot x=\phi|_{H_x}^{-1}(\gamma)$.
		Fix an open bisection $V\subseteq H$ with $\phi|_{H_x}^{-1}(\gamma)$, and define $y_\lambda:=s|_{V}^{-1}(x_\lambda)$.
		That is, $y_\lambda$ is the unique element of $V$ with source equal to $x_\lambda$.
		The net $(y_\lambda)$ converges to $\phi|_{H_x}^{-1}(\gamma)$ since $s(y_\lambda)\to x=s(\phi|_{H_x}(\gamma))$ and $s$ restricts to a homeomorphism between $V$ and $s(V)$.
		By construction, $(\phi|_{H_{x_\lambda}}^{-1}(\gamma_\lambda),y_\lambda^{-1})$ is a composable pair in $H$ for all $\lambda$, and we have
		$$\phi(\phi|_{H_{x_\lambda}}^{-1}(\gamma_\lambda)y_\lambda^{-1})=\gamma_\lambda\phi(y_\lambda)^{-1}\to \gamma\phi(\phi|_{H_x}^{-1}(\gamma))=\gamma\gamma^{-1}=r(\gamma),$$
		as $\phi$ is continuous.
		Since $G^{(0)}$ is open, the net $\gamma_\lambda\phi(y_\lambda^{-1})$ eventually consists of units.
		Hence $\gamma_\lambda$ is eventually equal to $\phi(y_\lambda^{-1})^{-1}=\phi(y_\lambda)$.
		Hence $\phi|_{H_{x_\lambda}}^{-1}(\gamma_\lambda)=\phi|_{H_{x_\lambda}}^{-1}(\phi(y_\lambda))=y_\lambda\to \phi|_{H_x}^{-1}(\gamma)$, since $s(y_\lambda)=x_\lambda$ and $\phi$ is fibrewise bijective.
	\end{example}

	\begin{remark}
		Continuous fibrewise-bijective homomorphisms as in Example~\ref{eg-fibrewiseBijectiveMapsExample} are considered by a number of authors, going back to Brown \cite{Brown1}, where they are called \emph{star-bijective}.
		Barlak and Li \cite[Lemma~3.2]{BL1} showed that continuous proper fibrewise bijective maps induce homomorphisms between reduced groupoid $C^*$-algebras, and further between the twisted groupoid $C^*$-algebras if the homomorphisms respect the twist.
		In \cite{B1} continuous proper fibrewise bijective homomorphisms are named \emph{\'etale morphisms}, and it is shown that pullbacks of Fell bundles along these morphisms exist and induce ${}^*$-homomorphisms between the resulting section algebras.
	\end{remark}
	
	A more general version of the above example can be constructed using inverse semigroup actions by differing inverse semigroups.
	To this end, we need our second definition of transformation groupoid.
	
	\begin{definition}
		Let $S$ be a unital inverse semigroup and suppose $S$ acts on a locally compact Hausdorff space $X$ by local homeomorphisms.
		Consider $S\times X$ with the product topology, where $S$ is considered a discrete space.
		We define an equivalence relation on $S\times X$ by $(t,x)\sim (u,y)$ if and only if $x=y$ and there exists an idempotent $e\in E(S)$ such that $x$ is in the domain of $e$, and $te=ue$.
		This relation is called the \emph{germ relation}.
		Define $S\ltimes X$ as the quotient of $S\times X$ by $\sim$.
		There is a multiplication on $S\ltimes X$ given by
		$$[t,x]\cdot [u,y]=[tu,y],$$
		which is defined whenever $x=u\cdot y$.
		The inverse map is given by $[t,x]^{-1}=[t^*,t\cdot x]$, and the range and source maps are given by $s[t,x]=[1_S,x]$, $r[t,x]=[1_S,t\cdot x]$.
		The unit space is canonically homeomorphic to $X$ via the map $x\mapsto [1_S,x]$, and we often identify the two.
		One checks that the sets $\{[t,x]:x\in U\}$ for $t\in S$ and $U\subseteq X$ open form a basis of open bisections for $S\ltimes X$, whereby $S\ltimes X$ is \'etale.
	\end{definition}
	
	This definition agrees with Definition~\ref{defn-grpActionTrsfmGrpd}, since for a group $S=\Gamma$ the equivalence relation is trivial: there are no idempotents in a group save the identity.
	
	\begin{example}\label{eg-invSemigrpGrpdActor}
		Let $S$ and $T$ be unital inverse semigroups acting on locally compact Hausdorff spaces $X$ and $Y$ respectively by partial homeomorphisms.
		Call a pair $(\rho,\psi)$ an \emph{entertwining morphism} between the actions $S\curvearrowright X$ and $T\curvearrowright Y$ if $\rho:X\to Y$ is a continuous map, and $\psi:T\to S$ is a semigroup homomorphism, such that for any $t\in T$ we have $x\in\dom(\psi(t))$ if and only if $\rho(x)\in\dom(t)$, and
		$$t\cdot\rho(x)=\rho(\psi(t)\cdot x)$$
		for all $x\in X$ and $t\in T$ with $\rho(x)\in\dom(t)$.
		We define an anchor $\rho\circ r:S\ltimes X\to Y\cong (T\ltimes Y)^{(0)}$, which we shall also denote by $\rho$.
		For $[t,y]\in T\ltimes Y$ and $[u,x]\in S\ltimes X$ with $s[t,y]=y=\rho(u\cdot x)=\rho[u,x]$, define the multiplication
		$$[t,y]\cdot[u,x]:=[\psi(t)u,x].$$
		To see that this is well defined, suppose that $[t,y]=[t',y]$ in $T\ltimes Y$ and $[u,x]=[u',x]$ in $S\ltimes X$.
		That is, there exist idempotents $e\in E(T)$ and $f\in E(S)$ such that $te=t'e$, $uf=u'f$, and $y\in\dom(e)$, $x\in\dom(f)$.
		Since $y=\rho(u\cdot x)$ lies in the domain of $e$, we have $u\cdot x\in\dom(\psi(e))$.
		Then $\psi(t)u\cdot x=\psi(t)\psi(e)uf\cdot x=\psi(t'e)u'f\cdot x=\psi(t')u'\cdot x$, so the multiplication is well-defined.
		
		This defines an actor $T\ltimes Y\curvearrowright S\ltimes X$.
	\end{example}
	The above Example~\ref{eg-invSemigrpGrpdActor} in fact covers all actors of \'etale groupoids.
	
	\begin{proposition}
		Let $G$ be an \'etale groupoid.
		Then $\Bis(G)$ acts canonically on $G^{(0)}$ by partial homeomorphisms and $G$ is isomorphic to $\Bis(G)\ltimes G^{(0)}$.
		Moreover, any actor $G\curvearrowright H$ comes from a semigroup homomorphism $\psi:\Bis(G)\to\Bis(H)$ and a continuous map $\rho:H^{(0)}\to G^{(0)}$ that entwines the actions of $\Bis(G)$ and $\Bis(H)$.
		\begin{proof}
			An open bisection $U\subseteq G$ acts on $G^{(0)}$ via the partial homeomorphism $r\circ s|_U^{-1}$, and these form the action.
			The isomorphism $G\to\Bis(G)\ltimes G^{(0)}$ maps $\gamma\in G$ to the germ $[U_\gamma,s(\gamma)]$, where $U_\gamma$ is an open bisection neighbourhood of $\gamma$.
			One readily checks that this gives an isomorphism of groupoids.
			
			Let $h:\Bis(G)\ltimes G^{(0)}\curvearrowright\Bis(H)\ltimes H^{(0)}$ be an actor.
			Define $\psi:\Bis(G)\to\Bis(H)$ by $\psi(U):=U\cdot_h H^{(0)}$.
			Each $\psi(U)$ is a bisection in $H$ by Lemma~\ref{lem-actorBisectionProp}, and $\psi$ is a semigroup homomorphism, since 
			\begin{align*}
				\psi(UV)&=(UV)\cdot_hH^{(0)}\\
				&=U\cdot_h(V\cdot_hH^{(0)})\\
				&=U\cdot_h(H^{(0)}(V\cdot_hH^{(0)}))\\
				&=(U\cdot_hH^{(0)})(V\cdot_hH^{(0)})\\
				&=\psi(U)\psi(V).
			\end{align*}
			The anchor map for the actor $h$ restricts to a continuous map $\rho:H^{(0)}\to G^{(0)}$, and we claim this map entwines the actions of $\Bis(G)$ and $\Bis(H)$ via $\psi$.
			Firstly, we note that $t\in H^{(0)}$ belongs to the domain of $\psi(U)$ exactly when $x\in s(\psi(U))=\psi(s(U))=s(U)\cdot_h H^{(0)}$, which is exactly the preimage of $s(U)$ under $\rho$.
			For an open bisection $U\in\Bis(G)$ and a point $x\in\dom(\psi(U))$, the point $U\cdot\rho(x)$ is the range of the unique arrow $\gamma_x\in U\cap s^{-1}(\{\rho(x)\})$.
			Since $\gamma_x$ is the unique arrow of $U$ with $s(\gamma_x)=\rho(x)$, the unique arrow in $\psi(U)=U\cdot_h H^{(0)}$ with source $x$ is $\gamma_x\cdot x$.
			We then have
			$$\rho(\psi(U)\cdot x)=\rho(r(\gamma_x\cdot x))=r(\gamma_x)=U\cdot \rho(x),$$
			whereby $\rho$ entwines the actions of $\Bis(G)$ and $\Bis(H)$.
			
			Lastly, we must show that $h$ is exactly the actor described in Example~\ref{eg-invSemigrpGrpdActor}.
			We identify an open bisection $U\subseteq G$ with the set $\{[U,y]\in\Bis(G)\ltimes G^{(0)}: y\in s(U)\}$, noting that this is exactly the image of $U$ in $\Bis(G)\ltimes G^{(0)}$ under the isomorphism described above.
			For $[U,y]\in\Bis(G)\ltimes G^{(0)}$ and $[V,x]\in\Bis(H)\ltimes H^{(0)}$ with $\rho[V,x]=s[U,y]$, we have
			$$[\psi(U)V,x]=[(U\cdot_h H^{(0)})V,x]=[U\cdot_h V,x].$$
			Let $\gamma\in G$ and $\eta\in H$ be the arrows corresponding to $[U,y]$ and $[V,x]$ under the canonical isomorphisms $\Bis(G)\ltimes G^{(0)}\cong G$ and $\Bis(H)\ltimes H^{(0)}\cong H$ given above.
			That is, $\gamma\in U$ is the unique element with $s(\gamma)=y$ and $\eta\in V$ the unique element with $s(\eta)=x$.
			Then $\gamma\cdot_h\eta$ is contained in $U\cdot_h V$, and has source $s(\gamma\cdot_h\eta)=s(\eta)=x$, whereby it is represented by the class $[U\cdot_h V,x]$.
			Hence $h$ arises exactly from the pair $(\psi,\rho)$.
		\end{proof}
	\end{proposition}
	
	\section{Actors of Fell bundles over groupoids}
	
	Given a topological space $X$, one may consider a bundle of $C^*$-algebras over $X$ and construct an algebra of sections with pointwise operations.
	If the fibres of this bundle are noncommutative, then the resulting section algebra will also be noncommutative, so is never isomorphic to $C_0(G^{(0)})$ for a groupoid $G$.
	We can, however, consider arrows in a groupoid $G$ with unit space $X$ acting on the bundle over $X$ by considering partial equivalences between source and range fibre algebras of a given arrow.
	Kumjian \cite{K2} defines Fell bundles over groupoids in analogue to Fell's $C^*$-algebra bundles.
	This definition allows one to consider the arrows of a groupoid acting by Hilbert bimodules (or `partial Morita equivalences'), allowing for a richer class of noncommutative dynamical systems to be considered.
	In \cite{KM1}, Kwa\'sniewski and Meyer show that noncommutative Cartan inclusions as defined by Exel \cite{E1} arise as Fell bundle $C^*$-algebras for bundles over inverse semigroups, many examples of which come from the Fell bundles over groupoids.
	
	In this section, we define actors between Fell bundles over groupoids that give rise to ${}^*$-homomorphisms between the Fell bundle section $C^*$-algebras, and show for corical Fell bundles (as defined by Bice in \cite{B1}), such actors form a category and there is a functor from this category to the category of $C^*$-algebras with ${}^*$-homomorphisms.
	We mimic a number of results from the previous section on actors of groupoids, as well as giving particular attention to one-dimensional Fell bundles over groupoids, which may be identified with twists.
	Twists over \'etale groupoids are of particular interest, since all (essential) Cartan pairs have twisted groupoid models (see \cite{R1}, \cite{Raad1}, \cite{T1}).
	
	We recall first the definition of a Banach bundle, in particular Banach bundles not over a Hausdorff spaces.
	
	\begin{definition}[{\cite[Definition~2.1]{BE1}}]
		Let $X$ be a topological space.
		An \emph{upper-semicontinuous Banach bundle over $X$} consists of a topological bundle $p:A\to X$, where $p$ is an open continuous surjection, and the fibres $A_x:=p^{-1}(\{x\})$ are Banach spaces for each $x\in X$ satisfying
		\begin{enumerate}
			\item the map $A\to\RR^+$ given by $v\mapsto ||v||$ is upper-semicontinuous;
			\item the map $(a,b)\mapsto a+b$ is continuous on the subset $\{(a,b)\in A\times A: p(a)=p(b)\}$ (viewed as a subspace of $A\times A$);
			\item for each $\lambda\in\CC$, the map $a\mapsto\lambda a$ is continuous;
			\item if $(v_\lambda)\subseteq A$ is a net such that $p(v_\lambda)$ converges to $x\in X$ and $||v_\lambda||\to 0$, then $v_\lambda$ converges to $0_x$: the zero element in the fibre $A_x$.
		\end{enumerate}
	\end{definition}
	
	\begin{definition}[{see \cite[2.1]{K2}}]
		A \emph{Fell bundle} over an \'etale groupoid $G$ with locally compact Hausdorff unit space is an upper-semicontinuous Banach bundle $\Ee=(\Ee_\gamma)_{\gamma\in G}$ over $G$ together with a partially defined multiplication
		$$\cdot:\Ee^{(2)}:=\{(a,b)\in\Ee\times\Ee: (p(a),p(b))\in G^{(2)}\}\to \Ee,\qquad (a,b)\mapsto a\cdot b,$$
		where $p:\Ee\to G$ is the bundle map, and an involution
		$${}^*:\Ee\to\Ee,\qquad a\mapsto a^*,$$
		satisfying the following
		\begin{enumerate}
			\item $p(a\cdot b)=p(a)p(b)$ for $(a,b)\in\Ee^{(2)}$, so that $\Ee_\gamma\cdot\Ee_\eta\subseteq\Ee_{\gamma\eta}$ whenever $(\gamma,\eta)\in G^{(2)}$, and the multiplication restricts to a bilinear map $\Ee_\gamma\times\Ee_\eta\to\Ee_{\gamma\eta}$;
			\item the multiplication is associative, that is, $a\cdot(b\cdot c)=(a\cdot b)\cdot c$ for such composable $a,b,c\in\Ee$;
			\item the multiplication map $\Ee^{(2)}\to\Ee$ is continuous, where $\Ee^{(2)}$ is equipped with the subspace topology of the product topology from $\Ee\times\Ee$;
			\item $||a\cdot b||\leq ||a||\cdot||b||$ for all $(a,b)\in\Ee^{(2)}$;
			\item the involution ${}^*$ is conjugate linear, and $p(a^*)=p(a)^{-1}$ for all $a\in\Ee$, that is, $\Ee_\gamma^*\subseteq\Ee_{\gamma^{-1}}$ for all $\gamma\in G$;
			\item $(a^*)^*=a$ and $||a^*||=||a||$ for all $a\in \Ee$, and if $(a,b)\in\Ee^{(2)}$ then $(a\cdot b)^*=b^*\cdot a^*$;
			\item the involution ${}^*:\Ee\to\Ee$ is continuous;
			\item $||a^*a||=||a||^2$ for all $a\in\Ee$;
			\item $a^*a$ is a positive element of the $C^*$-algebra $\Ee_{s(p(a))}$ for all $a\in\Ee$.
		\end{enumerate}
		We say the Fell bundle $\Ee$ is \emph{saturated} if $\Ee_\gamma\cdot\Ee_\eta$ spans a dense subalgebra of $\Ee_{\gamma\eta}$ for all $(\gamma,\eta)\in G^{(2)}$.
	\end{definition}

	If necessary we shall write $\mu_{\gamma,\eta}:\Ee_\gamma\times\Ee_\eta\to\Ee_{\gamma\eta}$ for the multiplication in the fibres over $\gamma$ and $\eta$ in a Fell bundle, so that certain statements and diagrams later in this section are clearer.
	
	We wish to consider algebras of sections of Fell bundles, so we are in general only interested in elements of the Fell bundle that lie in the ranges of the particular sections in which we are interested.
	For a Banach bundle $p:\Ee\to X$ over a locally compact Hausdorff space $X$ a result of Douady and Soglio-H\'erault \cite[Appendix~C]{FD1} implies the bundle has `enough sections', in the sense that each element $\xi\in\Ee$ lies in the range of some continuous (partially defined) section $U\subseteq X\to\Ee$.
	Although we shall be considering Fell bundles over groupoids $G$ that are not globally Hausdorff, the sections we consider are constructed from local sections over locally compact Hausdorff sets in $G$.
	These particular open sets form a basis for the topology on $G$, and so we have `enough sections' in this sense.
	
	For a Fell bundle $\Ee$ over $G$, we write $\Ee^{(0)}$ for the restriction $\Ee|_{G^{(0)}}$ of the Fell bundle to the unit space $G^{(0)}$.
	The axioms in the definition of Fell bundles imply that $\Ee^{(0)}$ is a $C^*$-bundle (see \cite{BE1}).
	
	Throughout this section we let $G$ and $H$ be \'etale groupoids with locally compact Hausdorff unit spaces.
	We also let $\Ee=(\Ee_\gamma)_{\gamma\in G}$ and $\Ff=(\Ff_x)_{x\in H}$ be Fell bundles over $G$ and $H$ respectively.
	Let $p_\Ee:\Ee\to G$ and $p_\Ff:\Ff\to H$ be the covering maps, i.e. $p_\Ee(\xi)=\gamma$ exactly when $\xi\in\Ee_\gamma$, and similarly for $p_\Ff$.
	For an actor $h:G\curvearrowright H$ with anchor map $\rho$, we may overload the notation of $s:G\to G^{(0)}$ and $\rho:H\to G^{(0)}$ to also denote the pullbacks $s\circ p_\Ee:\Ee\to G^{(0)}$ and $\rho\circ p_\Ff:\Ff\to G^{(0)}$.
	This is to avoid cumbersome notation.
	
	\begin{definition}\label{defn-FellActor}
		A \emph{Fell actor} $A$ of $\Ee$ on $\Ff$, denoted $A:\Ee\curvearrowright\Ff$ consists of an actor $h:G\curvearrowright H$, together with a continuous map $M^A:\Ee\baltimes{s}{\rho}\Ff\to\Ff$ satisfying 
		\begin{enumerate}
			\item for each $(\gamma,x)\in G\baltimes{s}{\rho}H$ the restriction $M^A_{\gamma,x}:=M^A|_{\Ee_\gamma\times\Ff_x}$ is a bilinear map $\Ee_{\gamma}\times\Ff_x\to\Ff_{\gamma\cdot x}$;
			\item for each $\gamma,\eta\in G$ and $x\in H$ with $s(\gamma)=r(\eta)$, $s(\eta)=\rho(x)$, the following diagram commutes:
			\[
   				\begin{tikzpicture}[baseline=(current bounding box.west)]
			    	\node (1) at (0,1) {$(\Ee_\gamma\times \Ee_\eta) \times \Ff_x$};
      				\node (1a) at (0,0) {$\Ee_\gamma\times (\Ee_\eta \times \Ff_x)$};
      				\node (2) at (5,1) {$\Ee_{\gamma\eta} \times \Ff_x$};
      				\node (3) at (5,0) {$\Ee_\gamma\times \Ff_{\eta\cdot x}$};
      				\node (4) at (7,.5) {$\Ff_{\gamma\eta\cdot x}$};
      				\draw[<->] (1) -- node[swap] {ass.} (1a);
      				\draw[cdar] (1) -- node {$\mu_{\gamma,\eta}\times \Id_{\Ff_x}$} (2);
      				\draw[cdar] (1a) -- node[swap] {$\Id_{\Ee_\gamma}\times M^A_{\eta,x}$} (3);
     				\draw[cdar] (3.east) -- node[swap] {$M^A_{\gamma,\eta\cdot x}$} (4);
      				\draw[cdar] (2.east) -- node {$M^A_{\gamma\eta,x}$} (4);
    			\end{tikzpicture}
    		\]
    		\item for each $\gamma\in G$ and $x,y\in H$ with $s(\gamma)=\rho(x)$, $s(x)=r(y)$, the following diagram commutes:
    		\[
   				\begin{tikzpicture}[baseline=(current bounding box.west)]
			    	\node (1) at (0,1) {$(\Ee_\gamma\times \Ff_x) \times \Ff_y$};
      				\node (1a) at (0,0) {$\Ee_\gamma\times (\Ff_x \times \Ff_y)$};
      				\node (2) at (5,1) {$\Ee_{\gamma\cdot x} \times \Ff_y$};
      				\node (3) at (5,0) {$\Ee_\gamma\times \Ff_{xy}$};
      				\node (4) at (7,.5) {$\Ff_{\gamma\cdot xy}$};
      				\draw[<->] (1) -- node[swap] {ass.} (1a);
      				\draw[cdar] (1) -- node {$M^A_{\gamma,x}\times \Id_{\Ff_y}$} (2);
      				\draw[cdar] (1a) -- node[swap] {$\Id_{\Ee_\gamma}\times \mu_{x,y}$} (3);
     				\draw[cdar] (3.east) -- node[swap] {$M^A_{\gamma, xy}$} (4);
      				\draw[cdar] (2.east) -- node {$\mu_{\gamma\cdot x,y}$} (4);
    			\end{tikzpicture}
    		\]
    		\item for each $\gamma\in G$, $x,y\in H$ with $s(\gamma)=\rho(x)$, $r(y)=r(\gamma\cdot x)$, and each $\xi_\gamma\in\Ee_\gamma,f_x\in\Ff_x,g_y\in\Ff_y$, we have
    		$$M^A_{\gamma,x}(\xi_\gamma,f_x)^*g_y=f_x^*M^A_{\gamma^{-1},y}(\xi_\gamma^*,g_y).$$
    		
			We say the Fell actor $A$ is \emph{saturated} if for each $(\gamma,x)\in G\baltimes{s}{\rho}H$ the image of the multiplication map $M_{\gamma,x}:\Ee_\gamma\times\Ee_x\to\Ee_{\gamma\cdot x}$ spans a dense subspace of $\Ee_{\gamma\cdot x}$.
			
			We may write $M=M^A$ when the Fell actor $A$ is understood.
		\end{enumerate}
	\end{definition}
	
	The multiplication of a Fell bundle on itself gives rise to an actor, which we shall later see is the identity actor in the appropriate category.
	A Fell bundle is saturated precisely when the multiplication actor of the bundle on itself is saturated.
	We now show that such Fell actors can be characterised in terms of maps from fibres of $\Ee$ to the adjointable operators between fibres of $\Ff$.	
	For the following results we treat fibres of a Fell bundle over units in the groupoid as $C^*$-algebras.
	
	\begin{lemma}\label{lem-multDefinesStarHomsOnUnitSpaceFibres}
		Let $A$ be a Fell actor of $\Ee$ on $\Ff$.
		For each $t\in H^{(0)}$, the multiplication $M_{\rho(t),t}:\Ee_{\rho(t)}\times\Ff_t\to\Ff_t$ comes from a ${}^*$-homomorphism.
		That is, there is a ${}^*$-homomorphism $\nabla^A_t:\Ee_{\rho(t)}\to M(\Ff_t)$ satisfying
		$$M_{\rho(t),t}(a,b)=\nabla_t(a)b,$$
		for all $a\in\Ee_{\rho(t)}$ and $b\in\Ff_t$.
		If $A$ is saturated, then the maps $\nabla^A_t$ are non-degenerate.
		We may write $\nabla_t=\nabla^A_t$ when $A$ is understood.
		\begin{proof}
			For $a\in\Ee_{\rho(t)}$, we claim that the map $b\mapsto M_{\rho(t),t}(a,b)$ is an adjointable map on $\Ff_t$ (here adjointability is with respect to the canonical Hilbert module structure of $\Ff_t$ over itself).
			To see this, fix $b,c\in\Ff_t$ and compute:
			\begin{align*}
				\langle b,M_{\rho(t),t}(a,c)\rangle&=b^*M_{\rho(t),t}(a,c)\\
				&=M_{\rho(t),t}(a^*,b)^*c,&\text{(4) in Definition~\ref{defn-FellActor},}\\
				&=\langle M_{\rho(t),t}(a^*,b),c\rangle.
			\end{align*}
			We have shown that $\nabla_t(a):=M_{\rho(t),t}(a,\cdot)$ is adjointable with adjoint $\nabla_t(a^*)$.
			Thus $a\mapsto\nabla_t(a)$ defines a map $\Ee_t\to M(\Ff_t)$, and we have shown that it is ${}^*$-preserving.
			We note also that $\nabla_t$ is linear since $M_{\rho(t),t}$ is bilinear, and it is multiplicative by condition (2) in Definition~\ref{defn-FellActor}.
			Thus, $\nabla_t$ is a ${}^*$-homomorphism.
			
			By construction we have $\nabla_t(\Ee_{\rho(t)})\Ff_t=M_{\rho(t),t}(\Ee_t,\Ff_t)$.
			Hence if $A$ is saturated, we see that $\nabla_t(\Ee_{\rho(t)})\Ff_t$ spans a dense subspace of $\Ff_t$. 
			Thus $\nabla_t$ is non-degenerate.
		\end{proof}
	\end{lemma}

	We shall restrict our attention to saturated Fell bundles and actors, as these are the actors for which the multiplication maps $M_{\gamma,x}$ are surjective.
	This will be a convenient property to have later, when we wish to compose actors.
	
	\begin{lemma}\label{lem-saturatedActorsMultnSurjective}
		Let $A:\Ee\curvearrowright\Ff$ be a saturated Fell actor with underlying actor $h:G\curvearrowright H$.
		Suppose $\Ee$ and $\Ff$ are both saturated.
		Then $\Ff_{\gamma\cdot x}= M_{\gamma,x}(\Ee_\gamma,\Ff_x)$ for all $(\gamma,x)\in G\baltimes{s}{\rho}H$, that is, $M_{\gamma, x}$ is surjective.
		\begin{proof}
			Since $A$ is saturated, Lemma~\ref{lem-multDefinesStarHomsOnUnitSpaceFibres} gives that $\nabla_{r(\gamma\cdot x)}:\Ee_{r(\gamma)}\to M(\Ff_{r(\gamma\cdot x)})$ is non-degenerate.
			Since $\Ff$ is saturated, the Cohen-Hewitt factorisation theorem implies $\Ff_{r(\gamma\cdot x)}\cdot\Ff_{\gamma\cdot x}=\Ff_{\gamma\cdot x}$ giving 
			$$\Ff_{\gamma\cdot x}=\nabla_{r(\gamma\cdot x)}(\Ee_{r(\gamma)})\Ff_{\gamma\cdot x}=\{\nabla_{r(\gamma\cdot x}(a)f:a\in\Ee_{r(\gamma)}, f\in\Ff_{\gamma\cdot x}\}.$$
			Hence we see
			\begin{align*}
				\Ff_{\gamma\cdot x}&=\nabla_{r(\gamma\cdot x)}(\Ee_{r(\gamma)})\cdot\Ff_{\gamma\cdot x}\\
				&=M_{r(\gamma),\gamma\cdot x}(\Ee_{r(\gamma)},\Ff_{\gamma\cdot x})\\
				&=M_{\gamma,x}((\Ee_{r(\gamma)}\cdot\Ee_{\gamma}),\Ff_{x})\\
				&=M_{\gamma,x}(\Ee_\gamma,\Ff_x),
			\end{align*}
			where the final equality holds since $\Ee$ is saturated.
		\end{proof}
	\end{lemma}
	
	\begin{lemma}\label{lem-connectingMapsOnFibres}
		For $(\gamma,x)\in G\baltimes{s}{\rho}H$, the multiplication map $M_{\gamma,x}:\Ee_\gamma\times\Ff_x\to\Ff_{\gamma\cdot x}$ comes from a bounded, linear, left $\nabla^A_{r(\gamma\cdot x)}$-equivariant map $\nabla^A_{\gamma,x}:\Ee_\gamma\to\Ll(\Ff_x,\Ff_{\gamma\cdot x})$.
		That is, $\nabla^A_{\gamma, x}(a\cdot \xi_\gamma)=\nabla^A_{r(\gamma\cdot x)}(a)\nabla^A_{\gamma,x}(\xi_\gamma)$.
		Moreover, the adjoint of $\nabla^A_{\gamma,x}(\xi_\gamma)$ is given by $\nabla^A_{\gamma^{-1},\gamma\cdot x}(\xi_\gamma^*)$.
		We may write $\nabla_{\gamma,x}=\nabla^A_{\gamma,x}$ when $A$ is understood.
		\begin{proof}
			For $\xi_\gamma\in\Ee_\gamma$ we claim the map $\nabla_{\gamma,x}(\xi_\gamma):\Ff_x\to\Ff_{\gamma\cdot x}$, given by
			$$[\nabla_{\gamma,x}(\xi_\gamma)](f_x):=M_{\gamma,x}(\xi_\gamma,f_x)$$
			is adjointable.
			To see this, we note that $\nabla_{\gamma^{-1},\gamma\cdot x}(\xi_\gamma^*)$ satisfies
			\begin{align*}
				\langle g_{\gamma\cdot x},\nabla_{\gamma,x}(\xi_\gamma)f_x\rangle&=g_{\gamma\cdot x}^*\cdot M_{\gamma,x}(\xi_\gamma,f_x)\\
				&=M_{\gamma^{-1},\gamma\cdot x}(\xi_\gamma^*,g_{\gamma\cdot x})^*\cdot f_x\\
				&=\langle \nabla_{\gamma^{-1},\gamma\cdot x}(\xi_\gamma^*)g_{\gamma\cdot x},f_x\rangle
			\end{align*}
			for all $f\in\Ff_x$ and $g_{\gamma\cdot x}\in\Ff_{\gamma\cdot x}$, where $\nabla_{\gamma,x}(\xi_\gamma)$ is adjointable with adjoint $\nabla_{\gamma^{-1},\gamma\cdot x}(\xi_\gamma^*)$.
			Hence $\nabla_{\gamma,x}$ is also bounded and linear.
			To see that $\nabla_{\gamma,x}$ is $\nabla_{r(\gamma\cdot x)}$-equivariant, we first recall that $\Ll(\Ff_x,\Ff_{\gamma\cdot x})$ has a left Hilbert $M(\Ff_{r(\gamma\cdot x)})$-module structure (see \cite[Theorem~2.4]{Lance1}).
			Thus for $a\in\Ee_{r(\gamma)}$, $\xi_\gamma\in\Ee_\gamma$, and $f_x\in\Ff_x$ we note that
			$$\nabla_{\gamma\cdot x}(a\cdot \xi_\gamma)f_x=M_{\gamma,x}(a\cdot\xi_\gamma,f_x)=M_{r(\gamma),\gamma\cdot x}(a,M_{\gamma, x}(\xi_\gamma,f_x))=\nabla_{r(\gamma\cdot x)}(a)\nabla_{\gamma,x}(\xi_\gamma)f_x,$$
			whereby $\nabla_{\gamma,x}$ is left $\nabla_{r(\gamma\cdot x)}$-equivariant.
		\end{proof}
	\end{lemma}
	
	\begin{lemma}\label{lem-connectingMapsKindaMultiplicative}
		Let $\gamma,\eta\in G$ and $x\in H$ with $s(\gamma)=r(\eta)$ and $\rho(x)=s(\gamma)$.
		For all $\xi_\gamma\in\Ee_\gamma$ and $\zeta_\eta\in\Ee_\eta$ we have
		$$\nabla_{\gamma\eta,x}(\xi_\gamma\zeta_\eta)=\nabla_{\gamma,\eta\cdot x}(\xi_\gamma)\nabla_{\eta,x}(\zeta_\eta).$$
		In particular, if $x\in H^{(0)}$ and $\Ff_x$ is unital as a $C^*$-algebra, then $\nabla_{\eta,x}(\zeta_\eta)=\nabla_{\eta,x}(\zeta_\eta)1_x=M_{\eta,x}(\zeta_\eta,1_x)$.
		\begin{proof}
			For all $f_x\in\Ff_x$ we have
			\begin{align*}
				\nabla_{\gamma\eta,x}(\xi_\gamma\zeta_\eta)f_x&=M_{\gamma\eta,x}(\xi_\gamma\zeta_\eta,f_x)\\
				&=M_{\gamma,\eta\cdot x}(\xi_\gamma,M_{\eta,x}(\zeta_\eta,f_x))\\
				&=\nabla_{\gamma,\eta\cdot x}(\xi_\gamma)\nabla_{\eta,x}(\eta)f_x,
			\end{align*}
			where the second equality follows from the second axiom of Definition~\ref{defn-FellActor}.
		\end{proof}
	\end{lemma}
	
	If $A:\Ee\curvearrowright\Ff$ is a saturated Fell actor, then for each $(\gamma,t)\in G\baltimes{s}{\rho} H$, the products $\nabla_{\gamma,t}(\Ee_\gamma)\cdot\Ff_t$ span a dense subspace of $\Ff_{\gamma\cdot t}$. 
	
	We now aim to show that Fell actors induce ${}^*$-homomorphisms between Fell bundle $C^*$-algebras, and to this end we recall the construction of the Fell bundle $C^*$-algebra as in \cite{BE1}.
	Let $\Ee$ be a Fell bundle over $G$, and let $U\subseteq G$ be an open bisection.
	Since $G$ is \'etale and has locally compact Hausdorff unit space, the space $U$ is locally compact and Hausdorff.
	Write $\Cc_c(U,\Ee)$ for the collection of continuous compactly supported partial sections $U\to\Ee$, which have been extended by zero to global sections $G\to\Ee$.
	We warn that if $G$ is not Hausdorff, elements of $\Cc_c(U,\Ee)$ may not be continuous.
	Write $\Cc_c(G,\Ee)$ for the finite linear span of such functions over all open bisections $U\subseteq G$.
	If $G$ is Hausdorff, then $\Cc_c(G,\Ee)$ coincides with the continuous compactly supported sections of $\Ee$ via a partition of unity argument.
	We equip $\Cc_c(G,\Ee)$ with a convolution product
	$$f\ast g(\gamma)=\sum_{\eta\in G_{r(\gamma)}}f(\eta^{-1})\cdot g(\eta\gamma),$$
	and involution $f^*(\gamma)=f(\gamma^{-1})^*$, for $f,g\in\Cc_c(G,\Ee)$.
	We define $C^*(\Ee)$ as the completion of $\Cc_c(G,\Ee)$ under the maximal $C^*$-norm.
	As noted in \cite{BE1}, the same argument as \cite[Proposition~3.14]{E3} implies that such a maximal $C^*$-norm exists.
	
	By construction, $C^*(\Ee)$ has a dense subalgebra spanned by sections of $\Ee$ supported on bisections, allowing us to give prototypical constructions over bisections as in the case of \'etale groupoid $C^*$-algebras.
	
	\begin{lemma}\label{lem-fellMorphismPrimitiveForBisectionsDefn}
		Let $A:\Ee\curvearrowright\Ff$ be a Fell actor with underlying actor $h:G\curvearrowright H$.
		Let $\rho:H\to G^{(0)}$ be the anchor map of $h$.
		For open bisections $U\subseteq G$ and $V\subseteq H$, and sections $f\in\Cc_c(U,\Ee)$, $g\in\Cc_c(V,\Ff)$, there is a (restricted) section $\Phi^0(f)g:U\cdot V\to\Ff$ of $\Ff$ given by
		$$[\Phi^0(f)g](\gamma\cdot x)=M_{\gamma,x}(f(\gamma),g(x)).$$
		Moreover, $\Phi^0(f)g$ is continuous and compactly supported, so extends (by zero) to an element of $\Cc_c(U\cdot V,\Ff)$.
		\begin{proof}
			By Lemma~\ref{lem-actorBisectionProp}, $U\cdot V$ is an open bisection and any point $y\in U\cdot V$ can be expressed as $\gamma\cdot x$ for unique $\gamma\in U$ and $x\in V$.
			Hence the function $h$ is well defined.
			It is a section of $\Ff$ since $M_{\gamma,x}$ takes values in $\Ff_{\gamma\cdot x}$ by definition.
			
			To see that $\Phi^0(f)g$ is continuous, recall that the multiplication $m:G\baltimes{s}{\rho}H\to H$ associated to the actor $h$ restricts to a homeomorphism on $U\baltimes{s}{\rho}V$ by Corollary~\ref{cor-multRestrictsToHomeo}.
			Let $w:U\cdot V\to U\baltimes{s}{\rho}V$ be the inverse of $m$ restricted to $U\baltimes{s}{\rho}V$.
			We then see that $\Phi^0(f)g=M_{\gamma,x}\circ(f\times g)\circ w$, whereby $\Phi^0(f)g$ is continuous since each of the $M_{\gamma,x}$, $f$, $g$, and $w$ are.
			
			To see that $\Phi^0(f)g$ has compact support, let $K_f\subseteq U$ and $K_g\subseteq V$ be the compact supports of $f$ and $g$ respectively.
			Then $\Phi^0(f)g(\gamma\cdot x)\neq 0$ implies that $\gamma\in K_f$ and $x\in K_g$, since the multiplication map $M_{\gamma,x}$ is bilinear.
			Thus the support of $\Phi^0(f)g$ is contained in the image of $K_f\baltimes{s}{\rho} K_g$ under $m$, which is compact since $m$ is continuous.
		\end{proof}
	\end{lemma}
	
	The phrase `restricted section' in the above lemma means that we consider a function $\Phi^0(f)g:U\cdot V\to\Ff$ with the property that $p_\Ff(\Phi^0(f)g(y))=y$ for all $y\in U\cdot V$.
	This is restricted in the sense that we only consider fibres of $\Ff$ over points in the subset $U\cdot V\subseteq H$.
	
	\begin{lemma}\label{lem-multiplierOnFellSections}
		Let $A:\Ee\curvearrowright\Ff$ be a Fell actor with underlying actor $h:G\curvearrowright H$.
		Let $\rho:H\to G^{(0)}$ be the anchor map of $h$.
		Fix $f\in\Cc_c(G,\Ee)$.
		There is a map $\Phi_A^0(f):\Cc_c(H,\Ff)\to\Cc_c(H,\Ff)$ given by
		$$[\Phi^0_A(f)g](x):=\sum_{\gamma\in G_{\rho(x)}}M_{\gamma^{-1},\gamma\cdot x}(f(\gamma^{-1}),g(\gamma\cdot x)).$$
		\begin{proof}
			If $f\in\Cc_c(U,\Ee)$ and $g\in\Cc_c(V,\Ff)$ for open bisections $U\subseteq G$ and $V\subseteq H$, then $\Phi^0_A(f)g$ is exactly the section described in Lemma~\ref{lem-fellMorphismPrimitiveForBisectionsDefn}, so is an element of $\Cc_c(U\cdot V,\Ff)$.
			If it exists, the map $\Phi^0_A(f)$ is linear on $\Cc_c(H,\Ff)$ since $M_{\gamma,x}$ is bilinear for each $(\gamma,x)\in G\baltimes{s}{\rho}H$.
			The assignment $f\mapsto\Phi^0_A(f)$ is also linear for the same reason, thus $\Phi^0_A(f)$ is a linear map $\Cc_c(H,\Ff)\to\Cc_c(H,\Ff)$.
		\end{proof}
	\end{lemma}
	
	The linear map $\Phi^0_A(f)$ in Lemma~\ref{lem-multiplierOnFellSections} defines something reminiscent of a multiplier on $\Cc_c(H,\Ff)$, which one may expect to extend to a multiplier of $C^*(\Ff)$.
	If we wish for $\Phi^0_A(f)$ to be represented by a section in $\Cc_c(H,\Ff)$, we require (as in the case of groupoids) that the actor $h$ is proper.
	We moreover require that the maps $\nabla_t:\Ee_{\rho(t)}\to M(\Ff_t)$ for $t\in H^{(0)}$ take values in $\Ff_t$.
	
	\begin{definition}
		Let $A:\Ee\curvearrowright\Ff$ be a Fell actor with underlying actor $h:G\curvearrowright H$.
		Let $\rho:H\to G^{(0)}$ be the anchor map associated to $h$.
		We say $A$ is \emph{proper} if $h$ is proper and the maps $\nabla_{\gamma,x}:\Ee_{\gamma}\to \Ll(\Ff_x,\Ff_{\gamma\cdot x})$ each take values in $\Kk(\Ff_x,\Ff_{\gamma\cdot x})\subseteq\Ll(\Ff_x,\Ff_{\gamma\cdot x})$, the subspace of compact operators.
	\end{definition}
	
	\begin{proposition}\label{prop-properFellActorGivesDirectMap}
		Let $A:\Ee\curvearrowright\Ff$ be a proper Fell actor with underlying actor $h:G\curvearrowright H$ and associated anchor map $\rho$.
		For each $f\in\Cc_c(G,\Ee)$, the map $\Phi^0_A(f):\Cc_c(H,\Ff)\to\Cc_c(H,\Ff)$ is given by left multiplication by an element of $\Cc_c(H,\Ff)$.
		That is, we may consider $\Phi^0_A(f)\in\Cc_c(H,\Ff)$.
		\begin{proof}
			Since $\Phi^0_A$ is linear, it suffices to consider $f\in\Cc_c(U,\Ee)$ where $U\subseteq G$ is an open bisection.
			The Hilbert module completion $\Cc_0(U,\Ee)$ of $\Cc_c(U,\Ee)$ is a right Hilbert module over $C_0(s(U),\Ee)$, and the Cohen-Hewitt factorisation theorem gives sections $f_1\in\Cc_0(U,\Ee)$ and $f_2\in C_0(s(U),\Ee)$ with $f=f_1f_2$.
			In fact, both $f_1$ and $f_2$ can be taken to have compact support by multiplying by an Urysohn function for the support of $f$.
			Assuming we have done this, we define $L_f:U\cdot H^{(0)}\to\Ff$ by
			$$L_f(\gamma\cdot t):=M_{\gamma,t}(f_1(\gamma),\nabla_t(f_2(\rho(t)))),\qquad\text{for }(\gamma,t)\in U\baltimes{s}{\rho}H^{(0)}.$$
			
			First we shall show that $L_f$ is a well defined function $U\cdot H^{(0)}\to\Ff$.
			Let $f_1f_2=f_1'f_2'$ for $f_1,f_1'\in\Cc_c(U,\Ee)$ and $f_2,f_2'\in C_c(s(U),\Ee)$.
			Then for all $\eta\in\Ff_{\gamma\cdot t}$ we have
			\begin{align*}
				\llangle M_{\gamma,t}(f_1(\gamma),\nabla_t(f_2(\rho(t)))),\eta\rrangle&=M_{\gamma,t}(f_1(\gamma),\nabla_t(f_2(\rho(t)))))\eta^*\\
				&=M_{\gamma,(\gamma\cdot t)^{-1}}(f_1(\gamma),\nabla_t(f_2(\rho(t))))\eta^*)\\
				&=M_{\gamma,(\gamma\cdot t)^{-1}}(f_1(\gamma),M_{s(\gamma),(\gamma\cdot t)^{-1}}(f_2(\rho(t))),\eta^*))\\
				&=M_{\gamma,(\gamma\cdot t)^{-1}}(f_1(\gamma)\cdot f_2(\rho(t)),\eta^*)\\
				&=M_{\gamma,(\gamma\cdot t)^{-1}}(f_1'(\gamma)\cdot f_2'(\rho(t)),\eta^*)\\
				&=M_{\gamma,(\gamma\cdot t)^{-1}}(f_1'(\gamma),M_{s(\gamma),(\gamma\cdot t)^{-1}}(f_2'(\rho(t)),\eta^*))\\
				&=M_{\gamma,(\gamma\cdot t)^{-1}}(f_1'(\gamma),\nabla_t(f_2'(\rho(t)))\eta^*)\\
				&=\llangle M_{\gamma,t}(f_1'(\gamma),\nabla_t(f_2'(\rho(t)))),\eta\rrangle.
			\end{align*}
			Here we have used the associativity condition (2) in Definition~\ref{defn-FellActor}.
			Since the left inner product $\llangle\cdot,\cdot\rrangle$ on $\Ff_{\gamma\cdot t}$ is positive definite, we see that $M_{\gamma,t}(f_1(\gamma),\nabla_t(f_2(\rho(t))))=M_{\gamma,t}(f_1'(\gamma),\nabla_t(f_2'(\rho(t))))$ for all $(\gamma,t)\in U\baltimes{s}{\rho}H^{(0)}$, so $L_f$ is well defined.
			
			That $L_f$ is a section follows since $M_{\gamma,t}$ has range contained in $\Ff_{\gamma\cdot t}$ by definition.
			We also see that $L_f$ has compact support, since the support of $L_f$ is exactly the image $\supp(f)\baltimes{s}{\rho} H^{(0)}$ under the multiplication $M$, which is compact since $f$ has compact support and $\rho$ is proper.
			
			To see that $L_f$ is continuous on $U\cdot H^{(0)}$, we shall show that the product of $L_f$ with any $g\in\Cc_c(H^{(0)},\Ff)$ is continuous.
			Let $w:U\cdot H^{(0)}\to U\baltimes{s}{\rho}H^{(0)}$ be the inverse to the multiplication of the actor $h$, which exsists by Corollary~\ref{cor-multRestrictsToHomeo}. 
			We compute
			$$(L_fg)(\gamma\cdot t)=M_{\gamma,t}(\xi_{f,\gamma},\nabla_t(a_{f,\gamma}))g(t)=M_{\gamma,t}(\xi_{f,\gamma}\cdot a_{f,\gamma},g(t))=M_{\gamma,t}(f(\gamma),g(t)),$$
			whereby $L_fg=M\circ f\times g\circ w$, and so $L_fg$ is continuous.
			Expressing $f=f_1f_2$ as above then yields $L_{f_1f_2}=L_{f_1}L_{f_2}$, which is continuous since $L_{f_1}=M(f_1\circ\rho,\cdot)$ is an element of $C_0(H^{(0)},\Ff)$.
		\end{proof}
	\end{proposition}
	
	\begin{corollary}\label{cor-StarHomFellBundles}[{cf. \cite[Theorem~4.6]{B1}}]
		Under the conditions of Proposition~\textup{\ref{prop-properFellActorGivesDirectMap}}, the assignment $\Cc_c(G,\Ee)\to\Cc_c(H,\Ff)$, $f\mapsto \Phi^0_A(f)$ is a ${}^*$-homomorphism, so extends so a homomorphism $\Phi_A:C^*(\Ee)\to C^*(\Ff)$ of the full $C^*$-algebras of the Fell bundles $\Ee$ and $\Ff$.
		\begin{proof}
			Linearity of $\Phi^0_A$ is clear, and the properties in Definition~\ref{defn-FellActor} give that $\Phi^0_A$ is multiplicative and ${}^*$-preserving.
		\end{proof}
	\end{corollary}

	\begin{remark}
		If the underlying actor of a proper Fell actor arises from a continuous proper fibrewise bijective homomorphism of groupoids as in Example~\ref{eg-fibrewiseBijectiveMapsExample}, then the induced ${}^*$-homomorphism between Fell bundle $C^*$-algebras in Corollary~\ref{cor-StarHomFellBundles} coincides with the ${}^*$-homomorphism induced by a pullback along a continuous fibrewise bijective homomorphism as in \cite[Theorem~4.6]{B1}.
		The perspective of the proofs of Proposition~\ref{prop-properFellActorGivesDirectMap} and Corollary~\ref{cor-StarHomFellBundles} is quite different to the approach in \cite{B1}, and moreover does not require the hypothesis that the underlying groupoids be Hausdorff.
	\end{remark}

	The $C^*$-algebra associated to a Fell bundle over an \'etale groupoid has a canonical subalgebra given by the restriction of the Fell bundle to the unit space of the groupoid.
	If the Fell bundle is the trivial bundle $G\times\CC$ over an \'etale groupoid $G$, the $C^*$-algebra of the Fell bundle is simply the groupoid $C^*$-algebra $C^*(G)$, and the canonical subalgebra is given by the continuous functions on the unit space.
	Just as in the case for groupoid $C^*$-algebras and actors between groupoids, sections of a Fell bundle supported on open bisections are mapped to sections supported on open bisections by a ${}^*$-homomorphism induced from a Fell actor.
	In particular, such ${}^*$-homomorphisms map sections supported on the unit space to sections supported on the unit space of the codomain bundle.
	
	\begin{corollary}\label{cor-FellActorHomPreservesSubalg}
		Let $A:\Ee\curvearrowright\Ff$ be proper Fell actor and let $U\subseteq G$ be an open bisection.
		Fix $f\in\Cc_c(U,\Ee)$.
		The image of $f$ under the ${}^*$-homomorphism $\Phi_A$ in Corollary~\textup{\ref{cor-StarHomFellBundles}} has support contained in $U\cdot H^{(0)}$.
		In particular, the image of the canonical subalgebra $C^*(\Ee^{(0)})\subseteq C^*(\Ee)$ under $\Phi_A$ is contained in the canonical subalgebra $C^*(\Ff^{(0)})\subseteq C^*(\Ff)$.
	\end{corollary}

	We also wish to show that this construction taking a Fell actor to a ${}^*$-homomorphism is functorial, and to do that we need to define the composition of Fell actors $\Ee\curvearrowright\Ff$ and $\Ff\curvearrowright\Gg$ for Fell bundles over groupoids $G$, $H$, and $K$.
	If we choose to model this on the composition of groupoid actors $h:G\curvearrowright H$ and $k:H\curvearrowright K$, then for $\gamma\in G$, $t\in K$ with $\rho_{hk}(t)=s(\gamma)$, the multiplication $M_{\gamma,t}:\Ee_\gamma\times\Gg_t\to\Gg_{(\gamma\cdot\rho_k(t))\cdot t}$ should factor through $M_{\gamma,\rho_k(t)}:\Ee_\gamma\times\Ff_{\rho_k(t)}\to\Ff_{\gamma\cdot\rho_k(t)}$ and $M_{\gamma\cdot\rho_k(t),t}:\Ff_{\gamma\cdot\rho_k(t)}\times\Gg_t\to\Gg_{(\gamma\cdot\rho_k(t))\cdot t}$.
	Ideally we would simply substitute a `unit' in the $\Ff_{\rho_k(t)}$ and $\Ff_{\gamma\cdot\rho_k(t)}$ places and call it a day, but the $C^*$-fibre $\Ff_{\rho_k(t)}$ is not in general unital, and the bimodule $\Ff_{\gamma\cdot \rho_k(t)}$ does not immediately have such a notion.
	Moreover, to show that the composition of Fell actors is well-defined, we shall need continuous local sections of units.
	Hence we consider saturated Fell bundles with local unital sections, together with saturated Fell actors.
	This motivates the following definition.
	
	\begin{definition}[{see \cite[Definitions 1.10, 1.21]{B1}}]
		A \emph{topological category} is a category carrying a topology making the source, range, and product maps continuous.
		A Fell bundle $p:\Ee\to G$ is \emph{categorical} if it forms a topological category.
	\end{definition}
	
	Categorical Fell bundles are considered by Bice in \cite{B1}, where it is shown that the bundle map $p:\Ee\to G$ restricts to a homeomorphism between the spaces of units of each $\Ee$ and $G$ when viewed as categories (see \cite[Proposition~1.22]{B1}).
	In fact, \cite[Proposition~1.22]{B1} shows that existence of the section of units $G^{(0)}\to\Ee$ is equivalent to $\Ee$ being a topological Fell bundle, making it the right notion to consider.
	
	We shall adapt the proof of \cite[Proposition~1.22]{B1} to local sections of units, which shall be a useful tool later.
	
	\begin{lemma}[{cf. \cite[Proposition~1.22]{B1}}]\label{lem-localUnitSections}
		Let $\Ee$ be a categorical Fell bundle and let $U\subseteq G^{(0)}$ be a precompact open set of units.
		There exists a continuous compactly supported section $u\in C_c(G^{(0)},\Ee)$ with $u(x)=1_x$ for all $x\in U$.
		\begin{proof}
			Since $\Ee$ is a topological category, the source map $\mathbf{s}:\Ee\to\Ee$ is continuous.
			The restricted zero-section $G^{(0)}\to\Ee$ given by $x\mapsto 0_x$ is continuous since $\Ee$ is a Banach bundle, and so composing with the source map $\mathbf s:\Ee\to\Ee$ in the topological category $\Ee$ gives a global section of units, and this is continuous.
			
			For each $x\in\bar{U}$, let $V_x$ be a precompact open neighbourhood of $x$.
			Then $(V_x)_{x\in U}$ is an open cover of $\bar{U}$, so has a finite subcover $V_1,\dots, V_n$.
			Then $V:=\bigcup_{j=1}^n V_j$ is a precompact set that contains $\bar{U}$ in its interior. 
			Pick an Urysohn function $f:G^{(0)}\to [0,1]$ separating $\bar{U}$ and $X\setminus V$.
			Then $u:= f1_\Ee$ is continuous, and has support contained in $\bar{V}$ which is compact.
		\end{proof}
	\end{lemma}

	Akemann, Pederson, and Tomiyama prove that the the unit of sections is continuous in the strict topology of the fibres \cite[Theorem~3.3]{APT1}, but this may fail to be continuous even when all fibres are unital.
	
	\begin{example}[{\cite[Section~1.3]{B1}}]
		Let $\Ee$ be the subbundle of the trivial bundle $M_2(\CC)\times[0,1]$ given by
		$$\Ee:=\{(M,t)\in M_2(\CC)\times[0,1]: M\in\CC\Theta_{11}\text{ if }t=1\},$$
		that is, the fibre at $t=1$ consists of matrices with only top-left entry non-zero.
		Then each fibre of this $C^*$-algebra is unital, but the unit at $t=1$ is the same matrix as at all other fibres.
		Thus the section of units is the function $f:[0,1]\to \Ee$ given by
		$$f(t)=\begin{pmatrix}
			1&0\\
			0&1-\delta_1(t)
		\end{pmatrix},$$
		where $\delta_1:[0,1]\to\{0,1\}$ is the Kronecker delta at 1.
		The function $f$ is certainly not continuous, since the limit $\lim\limits_{t\to 1}f(t)$ does not exist in $\Ee$ (but does exist in $M_2(\CC)$).
	\end{example}

	We can now define the composition of saturated Fell actors between categorical Fell bundles.
	
	\begin{proposition}\label{prop-defineFellActComp}
		Let $\Ee$, $\Ff$, and $\Gg$ be saturated Fell bundles over $G$, $H$, and $K$ respectively, and let $A:\Ee\curvearrowright \Ff$ and $B:\Ff\curvearrowright\Gg$ be saturated Fell actors with respective underlying actors $h:G\curvearrowright H$ and $k:H\curvearrowright K$.
		Suppose that $\Ff$ is categorical.
		There is an actor $BA:\Ee\curvearrowright\Gg$ with underlying actor $kh:G\curvearrowright H$ defining an associative composition of Fell actors.
		In particular, categorical Fell bundles over \'etale groupoids together with Fell actors form a category.
		\begin{proof}
			For $(\gamma,t)\in G\baltimes{s}{\rho_{kh}}K$, let $\nabla^B_{r(t)}:\Ff_{\rho_k(t)}\to\Gg_{r(t)}$ be the connecting ${}^*$-homomorphisms from Lemma~\ref{lem-multDefinesStarHomsOnUnitSpaceFibres} for the actor $B$.
			Fix $(\gamma,t)\in G\baltimes{s}{\rho_{kh}}K$, $\xi_\gamma\in\Ee_\gamma$, and $\zeta_t\in\Gg_{t}$.
			By Lemma~\ref{lem-multDefinesStarHomsOnUnitSpaceFibres} the map $\nabla^B_{r(t)}$ is unital, whereby we have $\zeta_t=\nabla^B_{r(t)}(1_{\rho_k(t)})\zeta_t$, where $1_{\rho_k(t)}$ is the unit in the fibre $\Ff_{\rho_k(t)}$.
			We define $M^{BA}_{\gamma,t}:\Ee_\gamma\times\Ff_t\to\Ff_{\gamma\cdot t}$ by
			$$M^{BA}_{\gamma,t}(\xi_\gamma,\zeta_t):=M^{B}_{\gamma\cdot\rho_k(t),t}(M^{A}_{\gamma,\rho_k(t)}(\xi_\gamma,1_{\rho_k(t)}),\zeta_t)$$
			
			To see that $M^{BA}$ is continuous, fix a net $(\xi_\lambda,\zeta_\lambda)\in\Ee\baltimes{s}{\rho_{kh}}\Gg$ converging to $(\xi,\zeta)$.
			Since $\Ee\baltimes{s}{\rho_{kh}}\Gg$ inherits the subspace topology from $\Ee\times\Gg$, the nets $\xi_\lambda$ and $\zeta_\lambda$ each converge individually to $\xi$ and $\zeta$ respectively.
			Without loss of generality, we may assume there is a precompact open set $V\subseteq K^{(0)}$ with  $r(\zeta_\lambda),r(\zeta)\in V$, so that there is a section $f_V\in C_0(H^{(0)},\Ff)$ with $f_V(\rho(t))=1_{\rho_k(t)}\in\Ff_{\rho_k(t)}$ for each $t\in V$ by Lemma~\ref{lem-localUnitSections}.
			We then have
			$$M^{BA}(\xi_\lambda,\zeta_\lambda)=M^B(M^A(\xi_\lambda,f_V(\rho_k(r(\zeta_\lambda)),\zeta_\lambda)))\xrightarrow{\lambda\to\infty}M^B(M^A(\xi,f_V(\rho_k(r(\zeta)))),\zeta),$$
			since the maps $M^A$, $M^B$, $f_V$, $\rho_k$, and $r$ are all continuous.
			Hence $M^{BA}$ is continuous.
			
			To show that $M^{BA}$ satisfies axiom (2) of Definition~\ref{defn-FellActor}, we first note that by Lemma~\ref{lem-connectingMapsKindaMultiplicative} for $\gamma,\eta\in G$ and $x\in H^{(0)}$ with $s(\gamma)=r(\eta)$ and $s(\eta)=\rho_h(x)$, and for $\xi_\gamma\in\Ee_\gamma$ and $\xi_\eta\in\Ee_\eta$, we have
			$$M^A_{\gamma,r(\eta\cdot x)}(\xi_\gamma,1_{r(\eta\cdot x)})M^A_{\eta,x}(\xi_\eta,1_x)=\nabla^B_{\gamma\eta,x}(\xi_\gamma\xi_\eta)=M^A_{\gamma\eta,x}(\xi_\gamma\xi_\eta,1_x).$$
			Hence for $\gamma,\eta\in G$ as above and $t\in K$ with $\rho_{kh}(t)=s(\eta)$, and for $\xi_\gamma\in\Ee_\gamma$, $\xi_\eta\in\Ee_\eta$, and $\zeta_t\in\Gg_t$, we have 
			\begin{align*}
				M^{BA}_{\gamma,\eta\cdot t}(\xi_\gamma,&M^{BA}_{\eta,t}(\xi_\eta,\zeta_t))\\
				&=M^{BA}_{\gamma,\eta\cdot t}(\xi_\gamma,M^B_{\eta\cdot\rho_k(t),t}(M^A_{\eta,\rho_k(t)}(\xi_\eta,1_{\rho_k(t)}),\zeta_t))\\
				&=M^B_{\gamma\cdot \rho_k(\eta\cdot t),\rho_k(\eta\cdot t)}(M^A_{\gamma,\rho_k(\eta\cdot t)}(\xi_\gamma,1_{\rho_k(\eta\cdot t)}),M^B_{\eta\cdot\rho_k(t)}(M^A_{\eta,\rho_k(t)}(\xi_\eta,1_{\rho_k(t)}),\zeta_t))\\
				&=M^B_{\gamma\eta,\rho_k(t)}(M^A(\xi_\gamma,1_{\rho_k(\eta\cdot t)})M^A_{\eta,\rho_k(t)}(\xi_\eta,1_{\rho_k(t)}),\zeta_t)\\
				&=M^B_{\gamma\eta,\rho_k(t)}(M^A(\xi_\gamma\xi_\eta,1_{\rho_k(t)}),\zeta_t)\\
				&=M^{BA}_{\gamma\eta,t}(\xi_\gamma\xi_\eta,\zeta_t).
			\end{align*}
			
			Axiom (3) of Definition~\ref{defn-FellActor} follows since
			\begin{align*}
				M^{BA}_{\gamma,tu}(\xi_\gamma,\zeta_t\zeta_u)&=M^B_{\gamma\cdot \rho_k(tu),tu}(M^A_{\gamma,\rho_k(tu)}(\xi_\gamma,1_{\rho_k(tu)}),\zeta_t\zeta_u)\\
				&=M^B_{\gamma,\rho_k(t)}(M^A_{\gamma,\rho_k(t)}(\xi_\gamma,1_{\rho_k(t)}),\zeta_t)\zeta_u\\
				&=M^{BA}_{\gamma,t}(\xi_\gamma,\zeta_t)\zeta_u,
			\end{align*}
			for $\zeta_t,\zeta_u\in\Gg$ composable.
		
			To see that the fourth axiom of Definition~\ref{defn-FellActor} is satisfied by $M^{BA}$, we note that
			\begin{align*}
				M^{BA}_{\gamma,t}(\xi_\gamma,\zeta_t)^*\zeta_u&=M^B_{\gamma\cdot\rho_k(t),t}(M^A_{\gamma,\rho_k(t)}(\xi_\gamma,1_{\rho_k(t)}),\zeta_t)^*\zeta_u\\
				&=\zeta_t^*M^B_{(\gamma\cdot\rho_k(t))^{-1},u}(M^A_{\gamma,\rho_k(t)}(\xi_\gamma,1_{\rho_k(t)})^*,\zeta_u)\\
				&=\zeta_t^*M^B_{(\gamma\cdot\rho_k(t))^{-1},u}(M^A_{\gamma^{-1},r(\gamma\cdot\rho_k(t))}(\xi_\gamma^*,1_{r(\gamma\cdot\rho_k(t))}),\zeta_u)\\
				&=\zeta_t^*M^{BA}_{\gamma^{-1},u}(\xi_\gamma^*,\zeta_u).
			\end{align*}
			
			Lastly, we must show that this composition associates.
			Fix saturated Fell bundles $\Ee_i$ over groupoids $G_i$ for $i=1,2,3,4$.
			Let $A_i:\Ee_i\curvearrowright\Ee_{i+1}$ be saturated Fell actors with underlying actors $h_i:G_i\curvearrowright G_{i+1}$ and anchors $\rho_i:G_{i+1}\to G_i^{(0)}$ for $i=1,2,3$.
			Suppose that $\Ee_2$ and $\Ee_3$ are categorical.
			For $(\gamma,t)\in \Ee_1\baltimes{s}{\rho_{321}}\Ee_4$, $\xi_\gamma\in(\Ee_1)_\gamma$ and $\zeta_t\in(\Ee_4)_t$ we have
			\begin{align*}
				M^{A_3(A_2 A_1)}_{\gamma,t}(\xi_\gamma,\zeta_t)&=M^{A_3}_{\gamma\cdot \rho_{3}(t),t}(M^{A_2A_1}_{\gamma,\rho_{3}(t)}(\xi_\gamma,1_{\rho_{3}(t)}),\zeta_t)\\
				&=M^{A_3}_{\gamma\cdot \rho_{3}(t),t}(M^{A_2}_{\gamma\cdot\rho_{32}(t),\rho_{3}(t)}(M^{A_1}_{\gamma,\rho_{32}(t)}(\xi_\gamma,1_{\rho_{32}(t)}),1_{\rho_3(t)}),\zeta_t)\\
				&=M^{(A_3A_2)}_{\gamma\cdot\rho_{32}(t),t}(M^{A_1}_{\gamma,\rho_{32}(t)}(\xi_{\gamma},1_{\rho_{32}(t)}),\zeta_t)\\
				&=M^{(A_3A_2)A_1}_{\gamma,t}(\xi_\gamma,\zeta_t),
			\end{align*}
			hence composition of such actors is associative.
			The identity actor of a categorical Fell bundle $\Ee$ over $G$ on itself is given by the left multiplication $\mu=\mu_\Ee$ of the Fell bundle, with underlying actor of the left multiplication of the groupoid on itself.
			Concretely, if $A:\Ff\curvearrowright\Ee$ is an actor with multiplication map $M^A$, then the composition $\mu A$ has multiplication map $M^{\mu A}(f,\xi)=\mu(M^{A}(f,1_{r(\xi)}),\xi)=M^{A}(f,1_{r(\xi)}\xi)=M^{A}(f,\xi)$.
			If, reversely, $A:\Ee\curvearrowright\Ff$, then $M^{A\mu}(\xi,f)=M^A(\mu(\xi,1_{s(\xi)}),f)=M^A(\xi,f)$.
			Thus categorical Fell bundles over \'etale groupoids with Fell actors form a category.
		\end{proof}
	\end{proposition}		
	
	\begin{remark}
		Since the Fell bundles and actors $A:\Ee\curvearrowright\Ff$ and $B:\Ff\curvearrowright\Gg$ considered above are all saturated, one could apply Lemma~\ref{lem-saturatedActorsMultnSurjective} to write each element of $\Gg$ as $\nabla^B(a)\zeta$ for some $a\in\Ff$ and $\zeta\in\Gg$.
		One could then attempt to do away with the ad hoc assumption that the Fell bundles be categorical and ask if $M^{BA}(\xi,\nabla^B(a)\zeta):=M^B(M^A(\xi,a),\zeta)$ defines the right notion of composite Fell actor.
		Indeed, if $\Ff$ is categorical, then it is clear that this definition agrees with the one above.
		However, two obstacles arise.
		Firstly, when $\Ff$ is categorical, we may pick a local choice of continuous factorisation of $f\in\Cc_c(U\cdot V,\Gg)$ into $\Phi^0(f_1)f_2$ for some $f_1\in\Cc_c(U,\Ff)$ and $f_2\in\Cc_c(V,\Gg)$ by picking $f_1$ to be a local section of units in $\Ff$.
		If $\Ff$ is not a topological category, there is no canonical (local) continuous choice of Cohen-Hewitt factorisation $\Gg\to\Ff\times\Gg$, and without this it is not immediately obvious how one may show the map $M^{BA}$ is continuous.
		Secondly, axiom 2 of Definition~\ref{defn-FellActor} uses that $M^A(\xi_1,1)M^A(\xi_2,2)=M^A(\xi_1\xi_2,1)$ (where the notation `1' here is overloaded to denote the unit in the appropriate fibre of $\Ff$).
		This holds since the units in fibres in $\Ff$ commute with other operators on the fibres, so this argument fails if these units are replaced with arbitrary elements of the fibres.
		The other axioms of Definition~\ref{defn-FellActor} should still hold under this alternative definition of composite actor (although, the computations showing this become much more cumbersome).
	\end{remark}	
	
	Throughout the remainder of this article we shall assume all Fell bundles are saturated and categorical, and that all Fell actors are saturated.
	We shall now show that the construction of ${}^*$-homomorphisms from Fell actors is functorial.
	
	\begin{proposition}[{cf. \cite[Theorem~4.9]{B1}}]\label{prop-fellActorHomFunctorial}
		Let $A:\Ee\curvearrowright\Ff$ and $B:\Ff\curvearrowright\Gg$ be proper Fell actors with underlying actors $h:G\curvearrowright H$ and $k:H\curvearrowright K$.
		Let $\Phi^0_A:\Cc_c(G,\Ee)\to\Cc_c(H,\Ff)$, $\Phi^0_B:\Cc_c(H,\Ff)\to\Cc_c(K,\Gg)$, and $\Phi_{BA}^0:\Cc_c(G,\Ee)\to\Cc_c(K,\Gg)$ be the ${}^*$-homomorphisms induced in Proposition~\textup{\ref{prop-properFellActorGivesDirectMap}} by $A$, $B$, and $BA$ respectively.
		Then $\Phi_{BA}^0=\Phi_B^0\circ\Phi_A^0$, and so $\Phi_B\circ\Phi_A=\Phi_{BA}$.
		In particular, there is a functor taking saturated categorical Fell bundles and saturated Fell actors to $C^*$-algebras and ${}^*$-homomorphisms.
		
		\begin{proof}
			It suffices to show that $\Phi_{BA}^0(f)=\Phi_B^0(\Phi_A^0(f))$ for a section $f\in\Cc_c(U,\Ee)$ supported on a bisection $U\subseteq G$, since such sections span $\Cc_c(G,\Ee)$.
			For such $f$ we have that both $\Phi_{BA}^0(f)$ and $\Phi_B^0(\Phi_A^0(f))$ are supported on $U\cdot K^{(0)}=(U\cdot H^{(0)})\cdot K^{(0)}$, so we fix $(\gamma,t)\in U\baltimes{s}{\rho_{kh}} G^{(0)}$ and compute $\Phi^{BA}_0(f)g$ for $g\in\Cc_c(K^{(0)},\Gg)$.
			Since $g$ has compact support, the image $\rho_k(\supp(g))$ is also compact, Lemma~\ref{lem-localUnitSections} gives a section $u\in\Cc_c(H^{(0)},\Ff)$ such that $u(\rho_k(t'))=1_{\rho_k(t')}$ for all $t'\in\supp(u)$.
			Hence $\Phi_B^0(u)g=g$.
			We then have
			\begin{align*}
				[\Phi^0_{BA}(f)g](\gamma\cdot t)&=M^{BA}_{\gamma,t}(f(\gamma),g(t))\\
				&=M^B_{\gamma\cdot\rho_k(t),t}(M^A_{\gamma,\rho_k(t)}(f(\gamma),1_{\rho_k(t)}),g(t))\\
				&=M^B_{\gamma\cdot\rho_k(t),t}(M^A_{\gamma,\rho_k(t)}(f(\gamma),u(\rho_k(t))),g(t))\\
				&=M^B_{\gamma\cdot\rho_k(t),t}([\Phi^0_A(f)u](\gamma\cdot\rho_k(t)),g(t))\\
				&=[\Phi_B^0(\Phi_A^0(f)u)g](\gamma\cdot t)\\
				&=[\Phi_B^0(\Phi_A^0(f))\Phi_B^0(u)g](\gamma\cdot t)\\
				&=[\Phi_B^0(\Phi_A^0(f))g](\gamma\cdot t),
			\end{align*}
			whereby $\Phi_{BA}^0=\Phi^0_B\circ\Phi_A^0$.
		\end{proof}
	\end{proposition}
	
	For a saturated Fell bundle $\Ee$ over an \'etale groupoid $G$, Kwa\'sniewski and Meyer define a local multiplier algebra-valued conditional expectation $EL:C^*(\Ee)\to\Mloc(C_0(G^{(0)},\Ee))$ by describing an inverse semigroup action of $\Bis(G)$ on $C_0(G^{(0)},\Ee)$ (see \cite[Section 7]{KM2}).
	The Hilbert $C_0(G^{(0)},\Ee)$-bimodule associated to an open bisection $U\subseteq G$ is the space $C_0(U,\Ee)$ of $C_0$-sections of the bundle $\Ee$ restricted to $U$ with structure inherited from $C^*(\Ee)$ (see \cite[Lemma~7.3]{KM2}).
	The conditional expectation $EL$ is defined on $C_0(U,\Ee)$ by multiplying the ideal $C_0(U\cap G^{(0)},\Ee)$ of $C_0(G^{(0)},\Ee)$, that is, for $f\in C_0(U,\Ee)$ and $g\in C_0(U\cap G^{(0)},\Ee)$
	$$EL(f)g:=fg,$$
	which defines an element $EL(f)\in M(C_0(U\cdot (U\cap G^{(0)}),\Ee))=M(C_0(U\cap G^{(0)},\Ee))$.
	Moreover, we observe that for such $f\in C_0(U,\Ee)$, the multiplier $EL(f)$ is completely determined by the restriction of $f$ to the unit space $G^{(0)}$.
	Since $EL$ is linear, the expectation $EL(f)$ for a function $f\in\Cc_c(G,\Ee)$ is determined by the restriction of $f$ to $G^{(0)}$.
	Using this fact, we can mirror the characterisation of $\Mloc(C_0(X))=\bB(X)/\mM(X)$ for commutative $C^*$-algebras by showing that sections in $\Cc_c(G,\Ee)$ belong to the kernel of $EL$ if their restriction to $G^{(0)}$ has meagre support.
	
	Kwa\'sniewski and Meyer define the \emph{essential Fell bundle {$C^*$}-algebra} \cite[Definition~7.12]{KM2} of a Fell bundle $\Ee$ as the quotient of $C^*(\Ee)$ by the ideal $N_{EL}$, which is equal to $\{f\in C^*(\Ee): EL(f^*f)=0\}$ by \cite[Theorem~4.11]{KM2}.
	This motivates understanding how to characterise when a function is zero under the conditional expectation $EL$.
	
	\begin{lemma}\label{lem-meagreUnitSupport0Exp}
		Let $f\in\Cc_c(G,\Ee)$ be a section with meagre support when restricted to $G^{(0)}$.
		Then $EL(f)=0$.
		\begin{proof}
			We note that for any dense open subset $V\subseteq G^{(0)}$, the algebra $C_0(V,\Ee)$ is an essential ideal in $C_0(G^{(0)},\Ee)$.
			Suppose $f\in\Cc_c(G,\Ee)$ has meagre support when restricted to $G^{(0)}$.
			We may express $f$ as $f=\sum_{j=1}^n f_n$ for some $f_j\in\Cc_c(U_j,\Ee)$, where $U_n\subseteq G$ are open bisections.
			We note that $f|_{G^{(0)}}$ must be discontinuous exactly on its support, since $f$ is zero on the comeagre (hence dense) complement of its support.
			For each $j\leq n$, the function $f_j$ restricted to $G^{(0)}$ can only have discontinuities at the boundary of $U_j\cap G^{(0)}$, since it is continuous on $U_j\cap G^{(0)}\subseteq U_j$.
			The boundary $C_j:=\partial U_j\cap G^{(0)}$ in $G^{(0)}$ is a closed set with empty interior for each $j\leq n$, hence the finite union $C:=\bigcup_{j=1}^n C_j$ is also a closed set with empty interior, and $C$ contains all points where $f$ may be continuous.
			Thus $C$ contains the support of $f$, so $EL(f)$ annihilates $C_0(G^{(0)}\setminus C,\Ee)$, which is an essential ideal of $C_0(G^{(0)},\Ee)$.
			Hence $EL(f)=0$.
		\end{proof}
	\end{lemma}	
	
	As in the case for actors between groupoids, we wish to see which ${}^*$-homomorphisms arising from proper Fell actors descend to essential $C^*$-algebras, and we shall again employ freeness.
	
	\begin{proposition}\label{prop-fellActorEntwines}
		Let $A:\Ee\curvearrowright\Ff$ be a proper Fell actor with underlying actor $h:G\curvearrowright H$ and associated anchor map $\rho:H\to G^{(0)}$.
		Let $\Phi_A:C^*(\Ee)\to C^*(\Ff)$ be the associated ${}^*$-homomorphism as in Corollary~\textup{\ref{cor-StarHomFellBundles}}.
		Suppose that $G$ is Hausdorff, or that $\rho$ is skeletal.
		If $h$ is free, then $\Phi_A$ descends to a ${}^*$-homomorphism $C^*_\ess(\Ee)\to C^*_\ess(\Ff)$.
		If, $G$ is Hausdorff, then $\Phi_A$ entwines conditional expectations.
		\begin{proof}
			Firstly suppose that $G$ is Hausdorff.
			Then all sections in $\Cc_c(G,\Ee)$ are continuous and so $EL_\Ee$ is a genuine conditional expectation taking values in $C_0(G^{(0)},\Ee)$.
			For any $f\in\Cc_c(G,\Ee)$ and $t\in H^{(0)}$ we have
			$$[\Phi_A(EL_\Ee(f))](t)=f(\rho(t))=\sum_{\gamma\cdot t=t} f(\gamma)=(\Phi_A(f)).$$
			Hence $\Phi_A\circ EL_\Ee= EL_\Ff\circ\Phi_A$, that is, $\Phi_A$ entwines conditional expectations.
			In particular $\Phi_A$ descends to a ${}^*$-homomorphism of the essential $C^*$-algebras.		
		
			Suppose now that $\rho$ is skeletal.
			Fix $f\in C^*(\Ee)$ with $EL_\Ee(f^*f)=0$.
			That is, $f$ is an element of the ideal $N_{EL_\Ee}$, the quotient by which defines the essential $C^*$-algebra of $\Ee$.
			Fix $\varepsilon>0$ and bisections $U_1,\dots,U_n\subseteq G$ and sections $f_i\in\Cc_c(U_i,\Ee)$ such that $||f^*f-\big(\sum_{i=1}^nf_i\big)\big(\sum_{j=1}^nf_j\big)||<\varepsilon$.
			Label $f_0:=\sum_{i=1}^n f_i$
			Then $||EL_\Ee(f_0^*f_0))||<\varepsilon$, since $EL_\Ee$ is contractive.
			For each $1\leq i,j\leq n$, define $\tilde{U}_{i,j}:=(U_i^{-1}U_j\cap G^{(0)})\cup (G^{(0)}\setminus U_i^{-1}U_j)^\circ$, so that $EL(f_i^*f_j)$ is a multiplier of $C_0(\tilde{U}_{i,j},\Ee)$ and each $\tilde{U}_{i,j}$ is dense and open in $G^{(0)}$.
			The intersection $\tilde{U}:=\bigcap_{i,j=1}^n\tilde{U}_{i,j}$ is then again dense and open in $G^{(0)}$, and $EL_\Ee(f_0^*f_0)$ is given by a multiplier on $C_0(\tilde{U},\Ee)$, which we may identify with multiplication by the restriction of $f_0^*f_0$ to $\tilde{U}$.
			
			Since $\tilde{U}$ is dense and open in $G^{(0)}$, the preimage $\tilde{V}:=\rho^{-1}(\tilde{U})$ is dense and open as $\rho$ is skeletal.
			Thus $C_0(\tilde{V},\Ff)$ is an essential ideal of $C_0(H^{(0)},\Ff)$, and we have
			$$\left|\left|\Phi_A(f_0)|_{\tilde{V}}\right|\right|=\sup_{t\in\tilde{V}}\left|\left|\nabla_t\left(\sum_{i=1}^nf(\rho(t))\right)\right|\right|\leq\sup_{u\in\tilde{U}}\left|\left|\sum_{i=1}^n f(u)\right|\right|=\left|\left|EL\left(\sum_{i=1}^nf_i\right)\right|\right|<\varepsilon.$$
			Hence $||EL_\Ff(\Phi(f^*f)||<\varepsilon$ for all $\varepsilon>0$, whereby $\Phi(f)\in N_{EL_\Ff}$ and $\Phi_A$ descends to the essential $C^*$-algebras.
		\end{proof}
	\end{proposition}
	
	In general we cannot ask whether a ${}^*$-homomorphism between Fell bundle $C^*$-algebras entwines expectations unless the underlying groupoid of the domain Fell bundle is Hausdorff.
	This is because there is no obvious reason why the criteria of Lemma~\ref{lem-liftToLocMultAlg} should be satisfied for these ${}^*$-homomorphisms.

	\begin{remark}
		A Fell actor could still induce a ${}^*$-homomorphism that entwines conditional expectations without having an underlying free actor of groupoids as in Proposition~\ref{prop-fellActorEntwines}, since any extra isotropy arising from the actor not being free could lie in the kernels of the connecting maps $\nabla_{\gamma,x}$.
		An extreme example of this would be an actor $\Ee\curvearrowright 0_H$, where $0_H$ is the zero-bundle over the groupoid $H$.
		In this case, $C^*(0_H)=\{0\}$ and so all ${}^*$-homomorphisms $C^*(\Ee)\to C^*(0_H)$ entwine conditional expectations, no matter what the underlying actor $G\curvearrowright H$ is.
	\end{remark}
	
	One may consider a saturated Fell bundle $\Ee$ over a groupoid $G$ as a saturated Fell bundle over the inverse semigroup $\Bis(G)$ of open bisections over $G$ (see \cite[Section~7]{KM2}).
	Kwa\'sniewski and Meyer show that the (full) Fell bundle $C^*$-algebras for both these Fell bundles are then isomorphic \cite[Proposition~7.6]{KM2}, and the conditional expectation $C^*_\ess(\Ee)\to\Mloc(C^*(\Ee^{(0)}))$ is faithful by \cite[Theorem~4.11]{KM2}.
	We note that in this case $C^*(\Ee^{(0)})=C^*_\ess(\Ee^{(0)})$, since the conditional expectation is the identity on $C^*(\Ee^{(0)})$.
	
	\subsection{Twists over \'etale groupoids}
	We focus now on a particular class of Fell bundles over \'etale groupoids; namely the one-dimensional bundles.
	
	\begin{definition}[compare, for example, {\cite{K1}}]
		Let $G$ be a topological groupoid.
		A \emph{twist} $(G,\Sigma)$ over $G$ consists of a topological groupoid $\Sigma$ and a continuous central extension 
		$$G^{(0)}\times\TT\to\Sigma\to G.$$
	\end{definition}

	A twist over a groupoid is itself a groupoid, but is never \'etale.
	We shall however not stray far from the \'etale setting, as the class of twists we shall consider are over \'etale groupoids.
	For such groupoids, the source and range fibres are no longer discrete.
	Rather, for a twist $\Sigma$ over an \'etale groupoid, the source fibre at a unit $x\in \Sigma^{(0)}=G^{(0)}$ is isomorphic to $G_x\times\TT$, which consists of disjoint copies of the circle.
	
	Given a twist $(G,\Sigma)$ over an \'etale groupoid $G$, we may construct a complex line bundle $L_\Sigma:=\frac{\Sigma\times\CC}{\TT}$, where the quotient is by the circle action $\lambda(\sigma,z)=(\lambda\sigma,\bar\lambda z)$ on $\Sigma\times\CC$.
	Note that this action is well defined on the $\Sigma$-component since a twist is a central extension by $\TT$.
	This endows $L_\Sigma$ with a Fell line bundle structure over $G$.
	The bundle map $L_\Sigma\to G$ is given by $(\sigma,z)\mapsto p(\sigma)$, where $p:\Sigma\to G$ is the quotient map in the twist.
	
	Given a Fell line bundle $L$ over $G$, one can construct a twist over $G$.
	We say $u\in L_\gamma$ is \emph{unitary} if $u^*u=1_{s(\gamma)}$ and $uu^*=1_{r(\gamma)}$.
	Define $U(L)$ as the space of unitary elements of $L$ with its subspace topology, and multiplication inherited from the Fell line bundle structure.
	One readily checks that this is a topological groupoid with unit space given by the unit elements of the fibres $1_x\in L_x\cong \CC$ for $x\in G^{(0)}$.
	Thus the unit space of $\Sigma_L$ is homeomorphic to $G^{(0)}$, and the restricted bundle map $\Sigma_L\to G$ is a surjective groupoid homomorphism.
	The space $\Aa:=U(L)\cap L^{(0)}$ is an open subgroupoid of $U(L)$ consisting of unitaries of the line bundle belonging to fibres over the unit space of $G$.
	The map $G^{(0)}\times\TT\ni (x,z)\mapsto z\cdot 1_{x}\in\Aa$ is then an isomorphism and the sequence
	$$G^{(0)}\times\TT\cong \Aa\to U(L)\to G,$$
	is a twist over $G$.
	
	Starting with a twist $(G,\Sigma)$, elements of the twist $\Sigma$ can be identified with unitaries of the associated line bundle $L_\Sigma$ via $\Sigma\ni\sigma\mapsto[1,\sigma]\in U(L_\Sigma)\subseteq L_\Sigma$, and this gives rise to an isomorphism of the twists $(G,\Sigma)$ and $(G,U(L_\Sigma))$.
	Similarly, from a Fell line bundle $L$ over $G$ there is an isomorphism $L_{U(L)}\to L$ of Fell bundles given by $[z,u]\mapsto zu$.
	Thus there is a one-to-one correspondence between Fell line bundles over $G$ and twists over $G$ (see \cite{DKR1} for more details).
	Under this correspondence, we identify sections of the Fell line bundle $L_\Sigma$ with functions $f:\Sigma\to\CC$ that are $\TT$-anti-equivarient, that is, $f(z\sigma)=\bar{z}f(\sigma)$ for all $\sigma\in\Sigma$ and $z\in\TT$.
	In particular, since the twist is trivial over the units, we see that the Fell line bundle is also trivial over $G^{(0)}$.
	Hence the canonical subalgebra $C^*(L^{(0)})=C_0(G^{(0)},L)$ is canonically isomorphic to $C_0(G^{(0)})$.
	
	\begin{proposition}\label{prop-oneDFellBundleActorsAndTwistsActorsSame}
		Let $(G,\Sigma)$ and $(H,\Omega)$ be twists over \'etale groupoids.
		Let $L_\Sigma$ and $L_\Omega$ be the respective Fell line bundles over $G$ and $H$ given by the twists.
		A saturated Fell actor $L_\Sigma\curvearrowright L_\Omega$ is given uniquely by an actor $\Sigma\curvearrowright\Omega$, and vice versa.
		\begin{proof}
			Let $A:L_\Sigma\curvearrowright L_\Omega$ be a saturated Fell actor with underlying actor $h:G\curvearrowright H$.
			Let $\rho:H\to G^{(0)}$ be the anchor map for $h$.
			Let $p_{\Sigma}:\Sigma\to G$ and $p_\Omega:\Omega\to H$ be the quotient maps in the extensions.
			Without loss of generality we may consider $\Sigma=U(L_\Sigma)$ and $\Omega=U(L_\Omega)$, the twists given by unitaries of the line bundles.
			The multiplication $M:L_\Sigma\baltimes{s}{\rho}L_\Omega\to L_\Omega$ maps pairs of unitaries to unitaries, since for $u\in U(L_\Sigma)$ and $v\in U(L_\Omega)$ with $s(u)=\rho(v)$ we have
			\begin{align*}
				M_{p_\Sigma(u),p_\Omega(v)}(u,v)^*M_{p_\Sigma(u),p_\Omega(v)}(u,v)&=v^*M_{p_\Sigma(u^*),p_\Omega(u\cdot v)}(u^*,M_{p_{\Sigma}(u),p_{\Omega}(v)}(u,v))\\
				&=v^*M_{p_\Sigma(u^*),p_\Omega(u\cdot v)}(u^*u,v)\\
				&=v^*v\\
				&=1_{s(v)},
			\end{align*}
			and similarly $(u\cdot v)(u\cdot v)^*=1_{r(v)}$.
			Saturation of the actor $A$ is needed, since otherwise $M_{p_\Sigma(u^*),p_\Omega(u\cdot v)}(u^*u,v)=0$ can occur.
			The multiplication given by restricting the multiplication $M$ of the Fell actor $A$ to the unitaries of the Fell line bundles $L_\Sigma$ and $L_\Omega$, together with the anchor map $\rho\circ r:\Omega\to\Sigma^{(0)}$ defines an actor $\beta_A:\Sigma\curvearrowright\Omega$.
			
			Now let $\beta:\Sigma\curvearrowright\Omega$ be an actor.
			We define a multiplication $L_\Sigma\baltimes{s}{\rho}L_\Omega\to L_\Omega$ by
			$$[\sigma,z]\cdot[\omega,w]:=[\sigma\cdot\omega,zw].$$
			This defines a Fell actor $A_k:L_\Sigma\curvearrowright L_\Omega$.
			This actor is saturated, since elements of the form $[\sigma,1]\in L_\Sigma$ act by non-zero linear maps, and all fibres are one-dimensional.
			This actor restricts to the actor $k$ on the unitaries of $\Sigma$ and $\Omega$, and so this construction is at least a one-sided inverse to the construction above taking a Fell actor to an actor of the twists.
			
			Conversely, starting with a Fell actor $A:L_\Sigma\curvearrowright L_\Omega$ and constructing the induced actor $k_A$ on twists as above, we note that the Fell actor $A_{k_A}:L_{U(L_\Sigma)}\curvearrowright L_{U(L_\Omega)}$ acts by
			$$[[\sigma,1],1]\cdot[[\omega,1],1]=[[\sigma\omega,1]1]$$
			which is exactly the image of the actor $A$ under the isomorphisms $L_\Sigma\cong L_{U(L_\Sigma)}$ and $L_\Omega\cong L_{U(L_\Omega)}$ induced by the identifications $\Sigma\xrightarrow{\sim} U(L_\Sigma)$, $\sigma\mapsto[\sigma,1]$, and $\Omega\xrightarrow{\sim} U(L_\Omega)$, $\omega\mapsto[\omega,1]$. 
		\end{proof}
	\end{proposition}
	
	\begin{lemma}\label{lem-entwiningIffLineBundles}
		Let $L_1$ and $L_2$ be Fell line bundles over \'etale groupoids $G_1$ and $G_2$ with locally compact Hausdorff unit space.
		Let $A:L_1\curvearrowright L_2$ be a proper Fell actor with underlying actor $h:G_1\curvearrowright G_2$.
		If $h$ is free, then the associated ${}^*$-homomorphism $\Phi_A$ entwines conditional expectations.
		\begin{proof}
			Using the fact that $C^*(L_i^{(0)})\cong C_0(G^{(0)}_i)$ for $i=1,2$, let $\tilde\Phi_A:\Mloc(C_0(G_1^{(0)}L_1))\to\Mloc(\Gamma_0(G_2^{(0)},L_2)$ be the extension of $\Phi_A|_{\Gamma_0(G^{(0)}_1,L_1)}$ from Corollary~\ref{cor-alwaysLiftCommutative}.
			For $f\in\Cc_c(G,L_1)$ and $t\in H^{(0)}$ we have
			$$[\tilde{\Phi}_A(f|_{G^{(0)}})](t)=\nabla_t(f(\rho(t))=\sum_{\gamma\cdot t=t}\nabla_{\gamma,t}(f(\gamma))=[\Phi_A(f)](t)$$
			since $h$ is free.
			Hence $\Phi$ entwines expectations on a dense subspace of $C^*(L_1)$, hence on all of $C^*(L_1)$ by continuity.
		\end{proof}
	\end{lemma}
	
	\section{Cartan morphisms of commutative Cartan pairs}
	
	Let $A\subseteq B$ is an inclusion of $C^*$-algebras with $A$ commutative.
	The identity map on $A$ extends to a completely positive contractive linear map $B\to I(A)$, where $I(A)$ is Hamana's injective hull of $A$ (see \cite{H1}).
	Such completely positive contractive linear maps that restrict to the identity on $A$ are called \emph{generalised expectations} by Kwa\'sniewski and Meyer \cite{KM2}, and called \emph{pseudoexpectations} if the codomain is the injective hull.
	In the case where $A$ is commutative, we have already recalled that $I(A)=\Mloc(A)$, so pseudoexpectations and $\Mloc(A)$-valued expectations coincide for inclusions of commutative $C^*$-algebras.
	We call a $\Mloc(A)$-valued expectation a \emph{local expectation}.
	
	A \emph{normaliser} for the inclusion $A\subseteq B$ is an element $n\in B$ that conjugates $A$ to itself, that is, $n^*An$ and $nAn^*$ are contained in $A$.
	Kumjian \cite[1.6]{K1} showed that such elements implement dynamics on the spectrum of $A$.
	We write $N(A,B)$ to denote the collection of all normalisers of $A$ in the inclusion $A\subseteq B$, and we say that $A$ is a \emph{regular} subalgebra of $B$ is $N(A,B)$ spans a dense subspace of $B$.
	
	We recall the definition of a(n essential) Cartan pair, which Renault originally defines for separable $C^*$-algebras with genuine conditional expectations (that is, conditional expectations taking values directly in the subalgebra).
	In \cite{R1}, Renault shows that such pairs are exactly reduced twisted groupoid $C^*$-algebras of twists over \'etale effective locally compact Hausdorff second countable groupoids.
	The requirement of second countability has since been removed in work of Kwa\'sniewski and Meyer \cite{KM1}, as well as Raad \cite{Raad1}.
	In the author's PhD thesis \cite{T1}, it is shown that such results still hold when the conditional expectation takes values in the local multiplier algebra of $A$.
	The construction of the underlying groupoid and twist remain the same, however the groupoid need no longer be Hausdorff (and will not be unless the conditional expectation is genuine).
	In this setting, one must consider the essential twisted groupoid $C^*$-algebra instead of the reduced, which coincide of the groupoid model is Hausdorff (equivalently, when the conditional expectation is genuine).
	
	We shall show in this section that all non-degenerate ${}^*$-homomorphisms between essential Cartan pairs that preserve the Cartan structure arise uniquely as actors between the associated Fell line bundles.
	With this we show that the category of essential Cartan pairs with such morphisms is equivalent to the category of Weyl twists of such pairs, with Fell actors as morphisms.
	
	\begin{definition}[{\cite[Definition~6.0.1]{T1}, \cite[Definition~5.1]{R1}}]
		Let $A\subseteq B$ be an inclusion of $C^*$-algebras with $A$ commutative.
		We say the pair $A\subseteq B$ is an \emph{essential Cartan pair} if the following conditions hold
		\begin{enumerate}
			\item $A$ is a regular subalgebra of $B$, that is, the normalisers $N(A,B):=\{n\in B:n^*An+nAn^*\subseteq B\}$ span a dense subspace of $B$;
			\item $A$ is a maximal abelian subalgebra of $B$;
			\item there exists a faithful local conditional expectation $E:B\to\Mloc(A)$.
		\end{enumerate}
		If $E$ is a genuine conditional expectation, that is, takes values in $A$, we say the pair $A\subseteq B$ is a \emph{Cartan pair}.
	\end{definition}

	Renault's definition \cite[Defintion~5.1]{R1} originally also requires the inclusion $A\subseteq B$ to also be non-degenerate.
	This was later shown to be redundant by Pitts \cite[Theorem~2.6]{P1}.
	
	We recall the main result of Renault's work on Cartan pairs \cite{R1} and a result of the author's PhD thesis \cite{T1}.	
	
	\begin{theorem}[{\cite[Theorem~6.2.11]{T1}, \cite[Theorem~5.6]{R1}}]\label{thm-renaultTaylor}
		Let $A\subseteq B$ be an essential Cartan pair.
		There is an effective \'etale groupoid $G$ with locally compact Hausdorff unit space, and a twist $\Sigma$ over $G$ such that $B$ is isomorphic to $C^*_\ess(G,\Sigma)$ via an isomorphism that entwines conditional expectations.
		The groupoid and twist $(G,\Sigma)$ are unique up to isomorphism among twists over \'etale effective groupoids with locally compact Hausdorff unit space.
	\end{theorem}
	
	The groupoid $G$ and twist $\Sigma$ in Theorem~\ref{thm-renaultTaylor} are constructed as follows (see \cite[Section~4]{R1} for more details).
	Let $C_0(X)\subseteq B$ be a non-degenerate inclusion of $C^*$-algebras and let $n\in N(C_0(X),B)$ be a normaliser.
	It follows that $n^*n$ is an element of $C_0(X)$, and so $\dom(n):=\supp^\circ(n^*n)$ is an open subset of $X$.
	Kumjian \cite[1.7]{K1} shows there is a unique homeomorphism $\alpha_n:\dom(n)\to\dom(n^*)$ satisfying
	$$(n^*fn)(x)=(n^*n)(x)f(\alpha_n(x)),$$
	for all $x\in\dom(n)$ and $f\in C_0(X)$.
	The partial homeomorphisms $\alpha_n$ form an inverse subsemigroup of the semigroup of partial homeomorphisms on $X$, which then acts tautologically on $X$.
	The \emph{Weyl groupoid} is defined as the transformation groupoid associated to this action.
	
	Define $D:=\{(n,x):n\in N(C_0(X),B), x\in\dom(n)\}\subseteq N(C_0(X),B)\times X$ and equip $D$ with the subspace topology of the product topology.
	We define an equivalence relation on $D$ by $(n,x)\sim(m,y)$ if and only if $x=y$ and there exist functions $f,g\in C_0(X)$ such that $f(x),g(x)>0$ and $nf=mg$.
	The quotient $\Sigma:=D/\sim$ equipped with the multiplication
	$$[n,x][m,y]:=[nm,y],$$
	defined whenever $x=\alpha_m(y)$ gives $\Sigma$ the structure of a topological groupoid.
	Let $\Aa:=\{[f,x]:x\in\dom(f), f\in C_0(X), x\in\dom(f)\}\subseteq \Sigma$.
	Then $\Aa\hookrightarrow \Sigma\rightarrow G$ is a twist over $G$, where the surjection $\Sigma\to G$ is given by $[n,x]\mapsto [\alpha_n,x]$.
	
	The isomorphism in Theorem~\ref{thm-renaultTaylor} preserves all the structure of a Cartan pair: it preserves the subalgebra, entwines conditional expectations, and maps normalisers to normalisers.
	This motivates the following definition.
	
	\begin{definition}[{cf. \cite[Definition~6.4.1]{T1}}]\label{defn-cartanMorphism}
		Let $A_i\subseteq B_i$ be a regular non-degenerate inclusions of $C^*$-algebras with pseudoexpectations $E_i:B_i\to I(A_i)$.
		A ${}^*$-homomorphism $\Phi:B_1\to B_2$ is a \emph{Cartan morphism} if
		\begin{enumerate}
			\item $\Phi(A_1)\subseteq A_2$;
			\item $\Phi(N(A_1,B_1))\subseteq N(A_2,B_2)$;
			\item $\Phi(E_1(b))=E_2(\Phi(b))$ for all $b\in B_1$ with $E_1(b)\in A_1$.
		\end{enumerate}
		We say a Cartan morphism $\Phi$ is \emph{non-degenerate} if $\Phi$ is a non-degenerate ${}^*$-homomorphism.
		That is, $\Phi(B_1)B_2=B_2$.
	\end{definition}
	
	If $A_i\subseteq B_i$ are (commutative) Cartan pairs, Cartan morphisms are exactly those considered by Li in \cite{L1} when considering inductive limits of Cartan pairs.
	The above definition is more general in two ways: the subalgebras $A_i$ need not be commutative, and the conditional expectations $E_i$ need not be genuine.
	Given an essential Cartan pair $A\subseteq B$, the Weyl twist $(G,\Sigma)$ associated to $A\subseteq B$ gives rise to a Fell line bundle $L_\Sigma$ over $G$.
	We call this the \emph{Weyl line bundle} of the pair $(A,B)$.
	As noted earlier, the canonical inclusions $C_0(G^{(0)})\subseteq C^*_\ess(G,\Sigma)$ and $C_0(G^{(0)})\subseteq C^*_\ess(L_\Sigma)$ are isomorphic via a Cartan isomorphism.
	
	We shall show that free proper actors between Fell bundles over groupoids arising from (essential) Cartan pairs always lift to non-degenerate Cartan morphisms of the corresponding $C^*$-inclusions, but this is not a true statement in general.
	
	\begin{example}\label{eg-actorHomNotCartan}
		Consider $X=[0,1]\times\{-1,1\}$ and the homeomorphism $\alpha:X\to X$ given by $\alpha(x,n)=(x,-n)$.
		This homeomorphism induces an action of $G_1:=\ZZ/2\ZZ$ on $X$, and the transformation groupoid $G_2:=G_1\ltimes X$ is principal since $\alpha$ has no fixed points.
		Viewing $G_1=\{1,-1\}$ with multiplication, there is then an actor $h:G_1\curvearrowright G_2$ given by $m\cdot (n_1,(n_2,x))=(mn_1,(n_2,x))$.
		This actor is free and proper, but we claim the induced ${}^*$-homomorphism $\varphi_h:C^*_\ess(G_1)\to C^*_\ess(G_2)$ does not map normalisers to normalisers.
		Since $G_1$ is a group, its canonical subalgebra $C_0(\{1\})$ is isomorphic to the complex numbers.
		Viewing $C^*_\ess(G_1)$ as having two generators $\delta_1$ and $\delta_{-1}$ with relation $\delta_{n}\delta_{m}=\delta_{nm}$ for $n,m\in G_1$, we see $C_0(\{1\})$ embeds as the scalar multiples of $\delta_1$.
		One readily computes that the element $\nu:=\delta_1-i\delta_{-1}$ is a normaliser of $C_0(\{1\})\subseteq C^*_\ess(G_1)$.
		The image of $\nu$ under $\varphi_h$ is the function
		$$[\varphi_h(\nu)](m,(x,n))=\nu(m)=\begin{cases}1,& m=1\\
		-i,&m=-1.\end{cases}$$
		Since $G_2$ is principal and Hausdorff, normalisers of the inclusion $C_0(X)\subseteq C^*_\ess(G_2)$ are exactly functions supported on open bisections of $G_2$ by \cite[Proposition~4.7]{R1}.
		Hence $\varphi_h(\nu)$ cannot be a normaliser, since it has support equal to $G_2$.
	\end{example}
	
	The failure of the actor in Example~\ref{eg-actorHomNotCartan} to induce a Cartan morphism arises from two properties: the inclusion $\CC\subseteq C^*(\ZZ/2\ZZ)$ has normalisers not arising from functions supported on bisections, and there are functions in $C_0(G_2^{(0)})$ that do not lie in the image of $\CC$ under the induced ${}^*$-homomorphism.
	If the groupoid $G_1$ comes from a(n essential) Cartan pair, then the first problem goes away.
	If a ${}^*$-homomorphism of $C^*$-inclusions restricts to a surjective map between the subalgebras, this second problem goes away.
	
	\begin{lemma}
		Let $C_0(G^{(0)}_1)\subseteq C^*_\ess(G_1,\Sigma_1)$ be an essential Cartan pair with Weyl twist $(G_1,\Sigma_1)$, and let $L_1$ be the corresponding Fell line bundle associated to $(G_1,\Sigma_1)$.
		Let $L_2$ be a Fell line bundle over $G_2$, and let $A:L_1\curvearrowright L_2$ be a free and proper Fell actor.
		The ${}^*$-homomorphism $\Phi_A:C^*_\ess(L_1)\to C^*_\ess(L_2)$ induced by $A$ in Corollary~\textup{\ref{cor-StarHomFellBundles}} is a Cartan morphism of the pairs $C_0(G_1^{(0)})\subseteq C^*_\ess(L_1)$ and $C_0(G_2^{(0)})\subseteq C^*_\ess(L_2)$.
		\begin{proof}
			That $\Phi_A$ preserves subalgebras is Corollary~\ref{cor-FellActorHomPreservesSubalg}, and entwining of expectations follows from Lemma~\ref{lem-entwiningIffLineBundles}.
			Let $n\in N(C_0(G_1^{(0)}),C^*_\ess(L_1))$ be a normaliser.
			Using \cite[Lemma~6.2.7]{T1} (or \cite[Corollary~5.3]{R1} if $G_1$ is Hausdorff and second countable) we see that for any $f\in C_c(G_1^{(0)})$ the element $nf$ is represented in $C^*_\ess(L_1)$ by a section supported on a bisection $U\subseteq G_1$.
			The image $\Phi_A(nf)$ is then supported in the bisection $U\cdot G_2$ by Corollary~\ref{cor-FellActorHomPreservesSubalg}.
			Taking an approximate unit for $C_0(G_1^{(0)})$ consisting of such compactly supported functions, we see $\Phi_A(n)$ may be represented as a section with support on a bisection of $G_2$, and hence is a normaliser.
		\end{proof}
	\end{lemma}	
	
	\begin{lemma}
		Let $A_i\subseteq B_i$ be regular non-degenerate inclusions of $C^*$-algebras and let $\Phi:B_1\to B_2$ be a ${}^*$-homomorphism such that $\Phi(A_1)=A_2$.
		Then $\Phi$ preserves normalisers.
		\begin{proof}
			Let $n\in N(A_1,B_1)$ be a normaliser and fix $a_2\in A_2$.
			Since $\Phi(A_1)=A_2$, there exists $a_1\in A_1$ such that $\Phi(a_1)=a_2$ and we then have $\Phi(n)^*a_2\Phi(n)=\Phi(n^*a_1n)\in\Phi(A_1)=A_2$, whereby $\Phi$ preserves normalisers.
		\end{proof}
	\end{lemma}
	
	We shall now show that all non-degenerate Cartan morphisms between essential Cartan pairs are induced by Fell actors between the associated Weyl line bundles.	
	
	\begin{proposition}\label{prop-cartanMapsComeFromActors}
		Let $A_i\subseteq B_i$ be essential Cartan pairs.
		Let $L_1,L_2$ be the associated Weyl line bundles, so that the pairs $(A_i,B_i)$ and $(C_0(G_i^{(0)}),C^*_\ess(L_i))$ are isomorphic as Cartan pairs by \cite[Theorem~6.2.11]{T1} for each $i$.
		Every ${}^*$-homomorphism arising from a proper Fell actor $L_1\curvearrowright L_2$ with free underlying actor $G_1\curvearrowright G_2$ is a non-degenerate Cartan homomorphism, and all non-degenerate Cartan homomorphisms $(C_0(G_1^{(0)}),C^*_\ess(L_1))\to(C_0(G_2^{(0)}),C^*_\ess(L_2))$ are of this form.
		\begin{proof}
			Let $A:L_1\curvearrowright L_2$ be a proper Fell actor with free underlying actor $G_1\curvearrowright G_2$.
			The ${}^*$-homomorphism $\Phi_A$ entwines conditional expectations by Lemma~\ref{lem-entwiningIffLineBundles}, and maps $C_0(G_1^{(0)})$ to $C_0(G_2^{(0)})$ by Corollary~\ref{cor-FellActorHomPreservesSubalg}.
			To see that $\Phi_A$ maps normalisers to normalisers, we make use of the fact that $C_0(G_1^{(0)})\subseteq C^*_\ess(L_1)$ is an essential Cartan inclusion.
			The definition of the Weyl groupoid and twist (see \cite[Section~4]{R1}) gives that a normaliser for the inclusion $C_0(G_1^{(0)})\subseteq C^*_\ess(L_1)$ can be approximated by elements of $\Cc_c(U,L_1)$ for some open bisection $U\subseteq G$.
			The subspace $\Cc_c(U,L_1)$ is mapped to $\Cc_c(U\cdot_h G_2^{(0)},L_2)$ by Corollary~\ref{cor-FellActorHomPreservesSubalg}, which is contained in the normalisers of the range inclusion.
			Since the space of normalisers for an inclusion of $C^*$-algebras is closed, the map $\Phi_A$ is a Cartan morphism.
			To see that $\Phi_A$ is non-degenerate, note that the restriction of $\Phi_A$ to the subalgebra $C_0(G_1^{(0)})$ is given by the pullback along $\rho$, and is therefore non-degenerate since $\rho$ is globally defined on $G_2^{(0)}$.
			
			Let $\Phi:C^*(L_1)\to C^*(L_2)$ be a Cartan morphism.
			Without loss of generality, we take $G_i$ and $L_i$ to be the Weyl groupoids and Weyl line bundles for the inclusions $C_0(G_i^{(0)})\subseteq C^*_\ess(L_i)$ for each $i=1,2$.
			The restriction of $\Phi$ to a ${}^*$-homomorphism $C_0(G_1^{(0)})\to C_0(G_2^{(0)})$ is non-degenerate, so is the pullback of a proper continuous map $\rho_0:G_2^{(0)}\to G_1^{(0)}$.
			Define $\rho:G_2\to G_1^{(0)}$ as $\rho:=\rho_0\circ r$, where $r$ is the range map in $G_2$.
			For elements $[\alpha_n,t_1]\in G_1$ and $[\alpha_m,t_2]\in G_2$ with $s[\alpha_n,t_1]=t_1=\rho_0(\alpha_m(t_2)=\rho[\alpha_m,t_2]$, define $[\alpha_n,t_1]\cdot[\alpha_m,t_2]:=[\alpha_{\Phi(n)m},t_2]$.
			Similarly, for elements $[n,t_1]\in\Sigma_1$ and $[m,t_2]\in\Sigma_2$ with $s[n,t_1]=t_1=\rho[m,t_2]$ define $[n,t_1]\cdot[m,t_2]:=[\Phi(n)m,t_2]$.
			These two compostitions define actors $h:G_1\curvearrowright G_2$ and $A:\Sigma_1\curvearrowright\Sigma_2$.
			We identify the actor $A$ with the corresponding Fell actor $L_1\curvearrowright L_2$ via Proposition~\ref{prop-oneDFellBundleActorsAndTwistsActorsSame}.
			We must show now that $\Phi=\Phi_A$.
			Let $U\subseteq G_1$ be an open bisection and fix $f\in\Cc_c(U,\Sigma_1)$.
			By Corollary~\ref{cor-FellActorHomPreservesSubalg} the support of $\Phi_A(f)$ is contained in $U\cdot G_2^{(0)}$, so fix $[\alpha_n,x]\cdot[\alpha_m,t]\in U\cdot G_2^{(0)}$, where $[\alpha_n,x]\in U$ and $[\alpha_m,t]\in G_2^{(0)}$.
			We then have
			$$f[\alpha_n,x]=\Phi_A(f)([\alpha_n,x]\cdot[\alpha_m,t])=\Phi_A(f)[\alpha_{\Phi(n)m},t].$$ 
			If $f[\alpha_n,x]=0$ then $f^*f(x)=0$ and we have
			$$|\Phi(f)[\alpha_{\Phi(n)m},t]|^2=|\Phi(f^*f)(t)|=|f^*f(\rho(t))|.$$
			and since $[\alpha_m,t]\in G_2^{(0)}$ we have $\rho(t)=\rho(s[\alpha_m,t])=\rho(r[\alpha_m,t])=s[\alpha_n,x]=x$, which in turn implies $|\Phi(f)[\alpha_{\Phi(n)m},t]|^2=|f^*f(x)|=0$.
			Otherwise if $f[\alpha_n,x]\neq 0$ then $\alpha_{\Phi(n)m}$ and $\alpha_{\Phi_A(f)}$ have the same germ at $t$, that is, $[\alpha_{\Phi(n)m},t]=[\alpha_{\Phi_A(f)},t]$.
			Then
			$$\Phi_A(f)[\alpha_{\Phi_A(f)},t]=\Phi_A(f^*f)(t)=(f^*f)\circ\rho(t)=\Phi(f^*f)(t)=\Phi(f)[\alpha_{\Phi(f)},t].$$
			Lastly for $t'$ in a neighbourhood of $t$ we have 
			$$\alpha_{\Phi_A(f)}(t')=r[\alpha_{\Phi_A(f)},t]=r([\alpha_f,t]\cdot_h[\alpha_g,t])=r[\alpha_{\Phi(f)g},t]=\alpha_{\Phi(f)}(t)$$ 
			for any $g\in C_0(G_2^{(0)})$ satisfying $g(t)=1$.
			In particular $\Phi_A(f)[\alpha_{\Phi(f)},t]=\Phi(f)[\alpha_{\Phi(f)},t]$ for all $[\alpha_{\Phi(f)},t]\in U\cdot_h G_2^{(0)}$ where $\Phi_A(f)$ is non-zero.
			
			Thus $\Phi$ and $\Phi_A$ coincide on normalisers, which span a dense subspace of $C^*_\ess(L_1)$, and hence $\Phi=\Phi_A$.
		\end{proof}
	\end{proposition}
	
	\begin{corollary}\label{cor-mutualInverseFell1dBundle}
		The constructions in Proposition~\textup{\ref{prop-cartanMapsComeFromActors}} are mutually inverse.
		\begin{proof}
			The proof of Proposition~\ref{prop-cartanMapsComeFromActors} shows that if $\Phi$ is a non-degenerate Cartan morphism between $(C_0(G_1^{(0)}),C^*_\ess(L_1))$ and $(C_0(G_2^{(0)}),C^*_\ess(L_2)$), and $A_\Phi:L_1\curvearrowright L_2$ is the associated Fell actor, then $\Phi_{A_\Phi}=\Phi$, where $\Phi_{A_\Phi}$ is the Cartan homomorphism induced by $A_\Phi$.
			The assignment taking a non-degenerate Cartan morphism $\Phi$ to the actor $A_\Phi$ is thus injective with one-sided inverse $A\mapsto\Phi_A$.
			It is also surjective by Proposition~\ref{prop-cartanMapsComeFromActors}, as every such proper Fell actor is of the form $A_\Phi$ for some non-degenerate Cartan morphism $\Phi$.
			Thus the assignments $A\mapsto\Phi_A$ and $\Phi\mapsto A_\Phi$ are mutually inverse.
		\end{proof}
	\end{corollary}
	
	\begin{corollary}\label{cor-faithfulFullSucat}
		Let $\EssCart$ be the category of essential Cartan pairs with non-degenerate Cartan morphisms as arrows, and let $\Cart$ be the full subcategory of Cartan pairs.
		Let Let $\Fell_1$ be the category of Fell line bundles over \'etale groupoids with free proper Fell actors as morphisms, and let $\Fell_{1,\Haus}$ be the full subcategory of Fell line bundles over locally compact Hausdorff effective \'etale groupoids.
		The functor $\Ww:\EssCart\to\Fell_1$ taking an essential Cartan pair to the Fell line bundle associated to its Weyl twist, and taking a Cartan morphism $\Phi$ to the actor $\Ww\Phi:=A_\Phi$ is faithful, faithful on objects up to isomorphism, surjects to a full subcategory, and restricts to an equivalence of categories between $\Cart$ and $\Fell_{1,\Haus}$.
		\begin{proof}
			That $\Ww$ is faithful and surjects onto a full subcategory is Corollary~\ref{cor-mutualInverseFell1dBundle}.
			That $\Ww$ is faithful on objects up to isomorphism follows from \cite[Corollary~6.3.6]{T1}.
			That $\Ww$ restricts to an equivalence of categories follows from \cite[Proposition~4.11]{R1}, \cite[Proposition~3.4]{Raad1}.
		\end{proof}
	\end{corollary}
	
	\section{Inductive systems of actors}
	
	In this section we shall define inductive systems of free and proper groupoid actors, and show that they have colimits in the category of groupoids with actors.
	Moreover, we shall show that the $C^*$-algebra of such a colimit is exactly the inductive limit $C^*$-algebra induced by the inductive system of ${}^*$-homomorphisms arising from the actors.
	
	In \cite{L1}, Li constructs inductive limit groupoids for systems of Cartan morphisms between $C^*$-algebras with Cartan subalgebras.
	Li demonstrates that the injective Cartan morphisms between Cartan pairs induce a zig-zag on the underlying Weyl twists.
	These zig-zags consist of injective continuous open groupoid homomorphisms and proper continuous surjective fibrewise bijective groupoid homomorphisms coming from intermediate ``stepping stone'' groupoids and twists, constructed from the Cartan morphism.
	Inductive systems of Cartan morphisms between Cartan pairs then give rise to inductive systems of proper surjections in one direction, then open embeddings in the other.
	Barlak and Li \cite{BL1} construct an inductive limit groupoid for such a system, and show that $C^*$-algebra of this limit groupoid agrees with the inductive limit of the corresponding $C^*$-algebras.
	Combining these, Li is able to prove his groundbreaking result showing classifiable $C^*$-algebras have Cartan subalgebras (hence twisted groupoid models) by showing that the inductive limit of groupoids and twists gives a model for the inductive limit of $C^*$-algebras.
	
	In this section, we show that the stepping stone groupoid that Li describes is a special case of the transformation groupoid for free and proper actors, and can be constructed for homomorphisms between $C^*$-algebras arising from actors.
	We mimic Li's construction of the inductive limit groupoid for non-degenerate Cartan morphisms with this ersatz, and show that the corresponding groupoid $C^*$-algebra is the inductive limit $C^*$-algebra of the system, while also removing the requirement that the considered groupoids be globally Hausdorff.
	
	\begin{definition}
		A groupoid homomorphism $p:H\to G$ is \emph{fibrewise bijective} if $p$ restricts to a bijection $H_t\to G_{p(t)}$ on each source fibre $H_t$ for $t\in H^{(0)}$.
	\end{definition}	
	
	\begin{lemma}
		Let $h:G_1\curvearrowright G_2$ be a free and proper actor with surjective anchor map $\rho_h:G_2\to G_1^{(0)}$.
		There exists an \'etale groupoid $H$ with locally compact Hausdorff unit space, a continuous open injective groupoid homomorphism $i:H\to G_2$, and a continuous surjective, proper, and fibrewise bijective groupoid homomorphism $p:H\to G_1$ such that $\varphi_h$ is the continuous extension of $i_*\circ p^*:\Cc_c(G_1)\to\Cc_c(G_2)$.
		\begin{proof}
			Define $H\subseteq G_2$ by $H:=G_1\cdot_h G_2^{(0)}$.
			This is a subgroupoid since for $\gamma\cdot t\in H$ we have $(\gamma\cdot t)^{-1}=\gamma^{-1}\cdot r(\gamma\cdot t)$, and $(\eta\cdot r(\gamma\cdot t))(\gamma\cdot t)=\eta\gamma\cdot t$.
			It is open since $H$ is the union of $U\cdot G_2^{(0)}$, which are all open bisections by Lemma~\ref{lem-actorBisectionProp}.
			Define $i:H\to G_2$ as the inclusion map.
			
			We define $p:H\to G_1$ by $p(\gamma\cdot t)=\gamma$, which is surjective since $\rho_h$ is defined on all $G_2$.
			By definition, each element of $H$ is of the form $\gamma\cdot t$ for unique $\gamma\in G_1$ and $t\in G_2^{(0)}$, since the actor is free.
			From this definition, $p$ restricts to the anchor map $\rho_h$ of the actor $h$ since $t=\rho(t)\cdot t$ for all $t\in G_2^{(0)}$.
			The map $G_{\rho(t)}\to H_t$ given by $\gamma\mapsto\gamma\cdot t$ is a fibrewise inverse function for each $t\in G_2^{(0)}$, whereby $p$ is fibrewise bijective.
			To see that $p$ is proper, let $C\subseteq G_1$ be a compact subset.
			Recall that $\rho_h^{-1}(s(C))$ is compact as $\rho_h$ is proper.
			The preimage is then 
			$$p^{-1}(C)=\{x\in G_2: x=\gamma\cdot t\text{ for some }\gamma\in C, t\in G_2^{(0)}\}=C\cdot \rho_h^{-1}(C).$$
			This is precisely the image of $C\times\rho_h^{-1}(s(C))$ under the multiplication of the actor $h$, which is then compact since $\rho_h$ is proper.
			To see that $p$ is continuous, fix an open bisection $U\subseteq G_1$.
			The preimage $p^{-1}(U)=U\cdot\rho_h^{-1}(s(U))$ is then open since the multiplcation of the actor is a local homeomorphism by Corollary~\ref{cor-multRestrictsToHomeo}.
			
			Lastly, we show that $\varphi_h$ is the continuous extension of $i_*\circ p^*$.
			For an open bisection $U\subseteq G_1$ and $f\in\Cc_c(U)$ and $\gamma\cdot t\in U\cdot G_2^{(0)}$, we have
			$$i_*\circ p^*(f)(\gamma\cdot t)=f(p(\gamma\cdot t)=f(\gamma)=\varphi_h(f)(\gamma\cdot t).$$
			If $x\in G_2$ does not lie in $U\cdot G_2^{(0)}$, we then have either $x\notin H$ or $p(x)\notin U$.
			In both cases $i_*\circ p^*(f)(x)=0=\varphi_h(f)(x)$.
			Pullbacks and pushforwards are linear, and so $i_*\circ p^*$ agrees with $\varphi_h$ on $\Cc_c(G)$.
		\end{proof}
	\end{lemma}
	
	Following \cite{BL1}, we shall construct a groupoid model for the inductive limit $C^*$-algebra arising from the actors $h_n:G_n\curvearrowright G_{n+1}$.
	
	\begin{lemma}[{cf. \cite[Lemma~3.5]{BL1}}]\label{lem-projLimGrpdIsGrpd}
		Let $G_n$ be \'etale groupoids with locally compact Hausdorff unit spaces.
		Let $p_n:G_{n+1}\to G_n$ be proper continuous groupoid homomorphisms.
		Let $\bar{G}:=\varprojlim_n G_n$ be the projective limit as topological spaces and let $p^\infty_n:\bar{G}\to G_n$ be the canonical projections associated to the limit.
		Then $\bar{G}$ becomes an \'etale groupoid with locally compact Hausdorff unit space such that the $p^\infty_n$ are proper groupoid homomorphisms and $\bigcup_n(p^\infty_n)^*(\Cc_c(G_n))$ is dense in $\Cc_c(\bar{G})$ in the inductive limit topology.
		\begin{proof}
			The proof by Barlak and Li in \cite[Lemma~3.5]{BL1} works in this context, simply by removing the second countability condition and not requiring the projective limit groupoid to be Hausdorff.
			The projective limit groupoid $\bar{G}$ does have Hausdorff unit space, as the unit space of $\bar{G}$ is the projective limit of the unit spaces of the $G_n$, which are all Hausdorff.
		\end{proof}
	\end{lemma}
	
	Barlak and Li define inductive systems of (twists over) effective \'etale Hausdorff groupoids using their characterisation in the $C^*$-algebraic setting as Cartan pairs.
	For two such twists $(G_1,\Sigma_1)$ and $(G_2,\Sigma_2)$ and an injective Cartan morphism $\varphi:C^*_\red (G_1,\Sigma_1)\to C^*_\red(G_2,\Sigma_2)$, Barlak and Li \cite{BL1} define a pair $C\subseteq A$ of $C^*$-algebras as follows: $C$ is the ideal generated by the image of $\varphi(C_0(G_1^{(0)})$ in $C_0(G_2^{(0)})$, and $A:=C^*(C,\varphi(C^*_\red(G_1,\Sigma_1))\subseteq C^*_\red(G_2,\Sigma_2)$.
	One quickly verifies that $C\subseteq A$ is a Cartan pair if the inclusions $C_0(G_i^{(0)})\subseteq C^*_\red(G_i,\Sigma_i)$ are, and so \cite[Theorem~5.6]{R1} shows the pair $(C,A)$ is isomorphic to $(C_0(H^{(0)}), C^*_\red(H,\Omega))$ where $(H,\Omega)$ is the Weyl twist associated to $C\subseteq A$.
	There is a canonical embedding $i:(H,\Omega)\to(G_2,\Sigma_2)$ as the restriction of the twist $\Sigma_2$ to an open subgroupoid of $G_2$. 
	There is also a continuous, surjective, proper, and fibrewise bijective homomorphism of twists $p:(H,\Omega)\to (G,\Sigma)$, and these two maps together give rise to a ${}^*$-homomorphism $i_*\circ p^*:C^*_\red(G_1,\Sigma_1)\to C^*_\red(G_2,\Sigma_2)$.
	Li \cite[Proposition~5.4]{L1} shows that this recovers the original ${}^*$-homomorphism $\varphi$.
	
	We shall define (yet another) transformation groupoid for groupoid actions, which shall serve as our ersatz for the intermediate groupoid arising from the Cartan pair $C\subseteq A$.
	
	\begin{definition}{{\cite[Lemma~4.9]{MZ1}}}
		Let $h:G\curvearrowright X$ be a left action of $G$ on a topological space $X$ with anchor map $\rho:X\to G^{(0)}$.
		The transformation groupoid $G\ltimes_h H^{(0)}$ associated to $h$ is given by the following data:
		\begin{itemize}
			\item $G\ltimes_h X=\{(\gamma,x):\gamma\in G, x\in X, s(\gamma)=\rho(x)\}=G\baltimes{s}{\rho}X$, where the topology is the subspace topology of the product topology on $G\times X$;
			\item a pair $((\gamma,x),(\eta,y))$ is composable if and only if $x=\eta\cdot_h y$, and the multiplication is given by $(\gamma,x)(\eta,y)=(\gamma\eta,y)$.
		\end{itemize}
	\end{definition}
	
	The transformation groupoid is a groupoid by \cite[Lemma~4.9]{MZ1}.
	It is \'etale since it has a basis of open bisections of the form $U\baltimes{s}{\rho} V$, where $U\subseteq G$ is a bisection and $V\subseteq X$ is open.
	We also note that this definition agrees with Definition~\ref{defn-grpActionTrsfmGrpd} for actions of groups.
	
	\begin{lemma}[{cf. \cite[Proposition~4.16]{MZ1}}]\label{lem-transfGrpdHasMaps}
		Let $h:G\curvearrowright H$ be a free proper actor with surjective anchor $\rho$.
		The map $i:G\ltimes H^{(0)}\to H$ given by $i(\gamma,t)=\gamma\cdot t$ is a continuous injective groupoid homomorphism with open image that maps the unit space of $H^{(0)}$ identically to itself, and so $G\ltimes H^{(0)}$ embeds as an open subgroupoid of $H$.
		The map $p:G\ltimes H^{(0)}\to G$ given by $p(\gamma,t)=\gamma$ is a continuous, surjective, proper, and fibrewise bijective groupoid homomorphism.
		\begin{proof}
			If $i(\gamma,t)=i(\gamma',t')$ for some $(\gamma,t),(\gamma',t')\in G\ltimes H^{(0)}$ then $t=s(\gamma\cdot t)=s(\gamma'\cdot t')=t'$, and so $t=(\gamma\cdot t)^{-1}(\gamma'\cdot t')=\gamma^{-1}\gamma'\cdot t$.
			Since $h$ is free, this implies that $\gamma^{-1}\gamma'$ is a unit in $G$, whereby $\gamma=\gamma'$.
			Thus $i$ is injective.
			The map $i$ is continuous and open as it is the restriction of the multiplication of the actor $h$ to the open subset $G\baltimes{s}{\rho}H^{(0)}\subseteq G\baltimes{s}{\rho}H$, and the multiplication is a local homeomorphism by Corollary~\ref{cor-multRestrictsToHomeo}.
			It follows from the axioms of an actor that $i$ is a groupoid homomorphism.
			
			For $\gamma\in G$, there exists $t\in H^{(0)}$ such that $\rho(t)=s(\gamma)$, and so $(\gamma,t)\in G\ltimes H^{(0)}$ and $p$ is surjective.
			Properness of $p$ follows from properness of $\rho$, and $p$ is a groupoid homomorphism by definition of the multiplication in $G\ltimes H^{(0)}$.
			The map $p$ is continuous as it is the coordinate projection in $G\ltimes H^{(0)}$, which carries the subspace topology of the product topology on $G\times H^{(0)}$.
			Fix $t\in H^{(0)}$.
			For $(\gamma,t),(\eta,t)\in (G\ltimes H^{(0)})_t$, we have $p(\gamma,t)=p(\eta,t)$ if and only if $\gamma=\eta$, since $h$ is free.
			Thus $p$ is fibrewise bijective.
		\end{proof}
	\end{lemma}

	\begin{remark}
		The argument in Lemma~\ref{lem-transfGrpdHasMaps} combined with Example~\ref{eg-fibrewiseBijectiveMapsExample} show that continuous fibrewise bijective maps and free actors are in bijective correspondence with one another.
	\end{remark}
	
	We shall now show that if $\varphi_n:C^*_\ess(G_n)\to C^*_\ess(G_{n+1})$ are injective Cartan morphisms for groupoids giving rise to (essential) Cartan pairs, then the transformation groupoid $G_n\ltimes_{h_n}G_{n+1}$ associated to the actor $h_n:G_n\curvearrowright G_{n+1}$ inducing $\varphi_n$ is exactly the intermediate groupoid described by Bartak and Li \cite{BL1}.
	
	\begin{lemma}
		Suppose $G$ and $H$ are groupoids giving rise to essential Cartan pairs and let $\varphi:C^*_\ess(G)\to C^*_\ess(H)$ be an injective non-degenerate morphism of Cartan pairs.
		Let $A\subseteq C^*_\ess(H)$ be the $C^*$-algebra generated by $C_0(H^{(0)})$ and the image of $C^*_\ess(G)$ in $C^*_\ess(H)$.
		Then the pair $C_0(H^{(0)})\subseteq A$ is then isomorphic to the pair $C_0(H^{(0)})\subseteq C^*_\ess(G\ltimes_h H^{(0)})$, where $h:G\curvearrowright H$ is the actor inducing $\varphi$ from Proposition~\textup{\ref{prop-cartanMapsComeFromActors}}.
		\begin{proof}
			The pair $C_0(H^{(0)})\subseteq A$ is an essential Cartan pair by the same argument as in the proof of \cite[Proposition~5.4]{L1}: $C_0(H^{(0)})$ is clearly a regular and maximal abelian subalgebra of $A$.
			The local expectation $C^*_\ess(H)\to\Mloc(C_0(H^{(0)}))$ restricts to a faithful local expectation on the pair $C_0(H^{(0)})\subseteq A$.
			Since $H$ induces an essential Cartan pair, we may assume it is the Weyl groupoid of the pair $C_0(H^{(0)})\subseteq C^*_\ess(H)$ by \cite[Corollary~6.3.6]{T1}.
			The Weyl groupoid $K$ of $C_0(H^{(0)})\subseteq A$ embeds into $H$ directly as a subset.
			Specifically, $K=\{[\alpha_m,t]:m\in N(C_0(H^{(0)}),A),t\in\dom(n)\}\subseteq H$.
			
			Let $h:G\curvearrowright H$ be the unique actor inducing $\varphi$ by Proposition~\ref{prop-cartanMapsComeFromActors}.
			By Lemma~\ref{lem-transfGrpdHasMaps} the map $i:G\ltimes_h H^{(0)}\to H$ is an injective groupoid homomorphism with open image and maps the unit space $H^{(0)}$ identically to itself.
			We claim that $K$ is the image of $i$.
			The proof of \cite[Proposition~5.4]{L1} shows that each element in the Weyl groupoid of $C\subseteq A$ is of the form $[\alpha_{\varphi(n)},t]$ for some $n\in N(C_0(G^{(0)}),C^*_\ess(G))$.
			The point $\rho(t)=\varphi|_{C_0(G^{(0)})}^*(t)$ then lies in the open support of $n^*n$ since $t\in\dom(\varphi(n))$, and so the pair $([\alpha_n,\rho(t)],t)$ is an element of the transformation groupoid $G\ltimes_h H^{(0)}$.
			It follows that $[\alpha_{\varphi(n)},t]=[\alpha_n,\rho(t)]\cdot_h t=i([\alpha_n,\rho(t)],t)$ lies in the image of $i$, whereby $i$ is an isomorphism $G\ltimes_h H^{(0)}\to K$ of groupoids.
		\end{proof}
	\end{lemma}
	
	We shall now restrict our attention to free and proper actors with surjective anchor maps.
	Such actors induces injective ${}^*$-homomorphisms between the associated essential groupoid $C^*$-algebras by Proposition~\ref{prop-cartanMapInjIffInjOnSubalg}.
	By functoriality, an inductive system of such actors will induce inductive systems of $C^*$-algebras, and we shall then show that the construction of the essential groupoid $C^*$-algebra preserves these inductive limits.
	
	Let $(G_n,h_n)_{n\in\NN}$ be an inductive system of free and proper actors with surjective anchor maps.
	For each $n\in\NN$, the map $G_n\ltimes_{h_n}G_{n+1}^{(0)}\to G_{n+1}$ implementing the actor is an open embedding of groupoids by Lemma~\ref{lem-transfGrpdHasMaps}.
	Let $H_n$ be the image of $G_n\ltimes_{h_n}G_{n+1}^{(0)}$ in $G_{n+1}$, and let $p_n:H_n\to G_n$ be the continuous surjective proper fibrewise bijective groupoid homomorphism in Lemma~\ref{lem-transfGrpdHasMaps}.
	For $n,m\in\NN$ with $m\geq n$, let $h_n^m$ be the composition $h_{m-1}h_{m-2}\cdots h_n:G_n\curvearrowright G_m$ (where the empty composition is the identity actor on $G_n$).
	Let $H_n^m\subseteq G_m$ be the image of $G_n\ltimes_{h_n^m}G_m^{(0)}$ in $G_m$ under the open embedding in Lemma~\ref{lem-transfGrpdHasMaps}, and note that $H_n^m$ is always contained in $H_m$.
	Let $p_n^m:H_n^m\to G_n$ be the the continuous surjective proper fibrewise bijective groupoid homomorphism in Lemma~\ref{lem-transfGrpdHasMaps}.
	Define $G_{n,m}:=H_n^{n+m}$ for $n,m\in\NN$.
	The maps $p_{n+m}:G_{n,m+1}\to G_{n,m}$ define projective systems for each $n\in\NN$, and the limits $\bar{G}_n:=\varprojlim_m (G_{n,m},p_{n+m})$ are each \'etale groupoids with locally compact Hausdorff unit spaces by Lemma~\ref{lem-projLimGrpdIsGrpd}.
	
	Each projective limit $\bar{G}_n$ can be expressed as the subspace of $\prod_{m\in\NN}G_{n,m}$ compatible with the connecting maps $p_{n,m}=p_{n+m}$.
	That is, $(\gamma_{n,m})_{m\in\NN}\in \bar{G}_n$ if and only if $p_{m+n}(\gamma_{n,m+1})=\gamma_{n,m}$ for all $m\in\NN$.
	For each $n\in\NN$ there is a connecting map $\bar{\imath}_n:\bar{G}_n\to\bar{G}_{n+1}$ given by
	$$\bar\imath_n(\gamma_{n,m})_{m\in\NN}=(\gamma_{n,m+1})_{m\in\NN},$$
	i.e. deleting the first entry in the tuple and shifting all other entries down one index.
	Each map $\bar\imath_n$ is injective since the deleted entry in the shift can be recovered uniquely as $\gamma_{n,0}=p_{n+m}(\gamma_{n,1})$.
	
	Let $(\gamma_{n,m})_m\in\bar{G}_n$ be a generic element. 
	Since $G_{n,m}$ is the image of $G_n\ltimes_{h_{n}^{n+m}}G_{n+m}^{(0)}$ under the multiplication map induced by the actor, we may write $\gamma_{n,m}=\eta_{n,m}\cdot t_{n+m}$ for some $\eta_{n,m}\in G_n$ and $t_{n+m}\in G_{n+m}^{(0)}$ which are composable under $h_{n}^{n+m}$.
	The element $t_{n+m}$ is determined uniquely as the source of $\gamma_{n,m}$, and $\eta_{n,m}$ is determined as $p_n\circ p_{n+1}\circ\cdots\circ p_{n+m}(\gamma_{n,m})$ for any $m\in\NN$.
	In particular, $\gamma_{n,m}=\eta_n\cdot s(\gamma_{n,m})$ for $\eta_n=\eta_{n,m}$, and this does not depend on choice of $m\in\NN$.
	Thus a generic element of $\bar{G}_n$ takes the form $(\eta_n\cdot t_{n+m})_m$, where each $t_{n+m}\in G_{n+m}^{(0)}$ for $m\in\NN$ has anchor equal to the source of $\eta_n$.
	
	\begin{remark}
		Barlak and Li \cite{BL1} inductively define $G_{n,0}=G_n$ and $G_{n,m+1}:=p_{n+m}^{-1}(G_n)$.
		This definition gives the same groupoids as ours, since the map $p_n^m$ is exactly the composition $p_{n+m}\circ p_{n+m-1}\circ\cdots\circ p_n$.
	\end{remark}	
	
	\begin{lemma}\label{lem-inductOnProjGrpdLem}
		Let $(G_n,h_n)_{n\in\NN}$ be an inductive system of free and proper actors with surjective anchor maps.
		For $n\in\NN$ and $m\in\NN$, let $G_{n,m}$ be defined as above.
		Then $G_{n,m+k}\subseteq G_{n+k,m}$ for all such $n,m,k\in\NN$, and the maps $\bar\imath_n^{n+k}:\bar{G}_n\to\bar{G}_{n+k}$ given by $(\gamma_{n,m})_m\mapsto (\gamma_{n,m+k})_m$ are continuous injective groupoid homomorphisms with open image.
		\begin{proof}
			To show that $G_{n,m+1}\subseteq G_{n+1,m}$ we proceed by induction in $m$.
			For $m=0$ we have $G_{n,1}=p^{-1}_{n+m}(G_n)\subseteq G_{n+1}=G_{n+1,0}$.
			Assume now that $G_{n,m+1}\subseteq G_{n+1,m}$ for some $m\in\NN$.
			Then $G_{n,m+2}=p^{-1}_{n+m+1}(G_{n,m+1})\subseteq p^{-1}_{n+1+m}(G_{n+1,m})\subseteq G_{n+1,m+1}$, as required.
			Thus $G_{n,m+1}\subseteq G_{n+1,m}$ for all $n\in\NN$ and $m\in\NN$.
			That $G_{n,m+k}$ is contained in $G_{n+k,m}$ follows by induction in $k$.
			
			The shift maps $\bar\imath_n$ are continuous since their composition with the canonical projections $\bar{G}_{n+1}\to G_{n+1,m}$ are exactly the projections $\bar{G}_n\to G_{n,m+1}$.
			Each $\bar\imath_n$ is injective, since the zeroth entry of a tuple $(\gamma_{n,m})_m$ is uniquely determined as $\gamma_{n,0}=p_n(\gamma_{n,1})$.
			Since the inclusions $G_{n,m+1}\subseteq G_{n+1,m}$ are all groupoid homomorphisms, the induced maps $\bar\imath_n$ are too.
			The image of $\bar\imath_n$ is the preimage of $G_{n+1}$ under the canonical projection $\bar{G}_{n+1}\to G_{n+1,0}=G_{n+1}$, and so is open.
			The proof for $\bar\imath_n^k$ for $k>1$ follows inductively, as $\bar\imath_n^{n+k}=\bar\imath_{n+k-1}\circ\dots\circ\bar\imath_{n+1}\circ\bar\imath_n$.
		\end{proof}		 
	\end{lemma}

	\begin{lemma}
		Let $i_n:G_n\to G_{n+1}$ for $n\in\NN$ be injective open continuous groupoid homomorphisms between \'etale groupoids with locally compact Hausdorff unit spaces.
		The topological inductive limit $\varinjlim_{n}G_n$ equipped with pointwise multiplication is an \'etale groupoid with locally compact Hausdorff unit space.
		\begin{proof}
			For $k\geq n$ let $i_n^k:G_n\to G_k$ be the composition of maps $i_{k-1}\circ i_{k-2}\circ\dots\circ i_n$, where $i_n^n:=\Id_{G_n}$. 
			The topological inductive limit is the quotient of the disjoint union $\bigsqcup_{n}G_n$ by the relation $g_n\sim g_m$ for $g_n\in G_n$ and $g_m\in G_m$ if and only if there exists $k\geq n,m$ with $i_n^k(g_n)=i_m^k(g_m)$.
			
			The groupoids $G_m$ each embed into $G:=\varinjlim_n G_n$ by the composition $\iota^m:G_m\to\bigsqcup_{n}G_n\to G$, and the images of the $G_n$ form an open cover of $G$.
			Thus any element of $G$ is of the form $[g_n]$ for some $n\in\NN$ and some $g_n\in G_n$.
			Define 
			$$G^{(2)}:=\{([g_n],[g_m])\in G\times G: (i_n^k(g_n),i_m^k(g_m))\in G_k^{(2)}\text{ for some }k\geq n,m\}.$$
			Note that if $(i_n^k(g_n),i_m^k(g_m))\in G_k^{(2)}\text{ for some }k\geq n,m$ then $(i_n^k(g_n),i_m^k(g_m))\in G_k^{(2)}$ for all $k\geq n,m$.
			Define a multiplication $G^{(2)}\to G$ via $([g_n],[g_m])\mapsto [i_n^k(g_n)i_m^k(g_m)]$, where $k\in\NN$ is greater than $n$ and $m$.
			This is well defined, since for $k'\geq k$ we have $i_k^{k'}(i_n^k(g_n)i_m^k(g_m))=i_k^{k'}(i_n^k(g_n))i_k^{k'}(i_{m}^k(g_m))=i_n^{k'}(g_n)i_m^{k'}(g_m)$, as the connecting maps in the inductive limit are groupoid homomorphisms.
			The inclusions $\iota^n:G_n\to G$ are all then groupoid homomorphisms, forcing the multiplication in $G$ to be associative.
			
			The units in $G$ are exactly the images of units of building blocks $G_n$ under the canonical inclusions, and so the unit space of $G$ is the inductive limit of unit spaces of the $G_n$.
			In particular, $G^{(0)}$ is locally compact and Hausdorff.
		
			The source map $s_n$ in $G_n$ and the inclusion $\iota^n$ are both local homeomorphisms, and their composition locally defines the source map in $G$.
			Thus the source map in $G$ is a local homeomorphism.
			The inverse map is a homeomorphism on each $G_n$ and it commutes with the maps $i_n$, and hence it commutes with $\iota^n$.
			This describes the inverse map on $G$, and it follows that inversion is a homeomorphism on $G$.
			Hence $G$ is \'etale.
		\end{proof}
	\end{lemma}
	
	We are now able to show that inductive systems in the category of groupoid actors have colimits, as well as give a concrete description.
	The construction is originally due to Li \cite{L1}, where only effective Hausdorff groupoids are considered as they arise from Cartan pairs.
	Here we may consider non-effective groupoids, by taking the actors as part of our initial data.
	
	\begin{theorem}\label{thm-indLimGrpd}
		Let $(G_n,h_n)_{n\in\NN}$ be an inductive system of free and proper actors with surjective anchor maps.
		The groupoid $\bar{G}:=\varinjlim_{i_n}\varprojlim_{p_m}G_{n,m}$ is the inductive limit of the system $(G_n,h_n)_{n\in\NN}$, and the universal cocone is consists of actors $j_k:G_k\curvearrowright\bar{G}$ with anchor maps $\rho_k:\bar{G}\to G_k^{(0)}$ given by $\rho_k[(\eta_n\cdot t_{n+m})_m]=\rho_k^n(\eta_n)$ for $n\geq k$, and  the multiplication is given by
		$$\gamma_k\cdot[(\eta_n\cdot t_{n+m})_m]=[((\gamma_k\cdot\eta_n)\cdot t_{n+m})_m].$$
		\begin{proof}
			Each element of $\bar{G}$ takes the form $[(\gamma_{n,m})_m]$, where $(\gamma_{n,m})_m\in\bar{G}_n$ is some element of the projective limit groupoid $\bar{G}_n=\varprojlim_m (G_{n,m},p_{n+m})$.
			We define an actor $j_k:G_k\curvearrowright \bar{G}$ as follows: for $[(\gamma_{n,m})_m]\in \bar{G}$, the entries for $n<k$ are determined by the later entries in the sequence, so we may assume without loss of generality that $n\geq k$.
			Recall that $G_{n,m}\subseteq G_{n+m,0}=G_{n+m}$ by Lemma~\ref{lem-inductOnProjGrpdLem}.
			The actor $h_k^{n+m}:=h_{n+m}h_{n+m-1}\cdots h_k:G_k\curvearrowright G_{n+m}$ restricts to an actor on $G_{n,m}$, since for any $\gamma_k\in G_k$ and $\eta_n\cdot_{h_{n}^{n+m}} t_{n+m}\in G_{n,m}$ we have
			$$\gamma_k\cdot_{h_{k}^{n+m}}(\eta_n\cdot_{h_{n}^{n+m}} t_{n+m})=(\gamma_k\cdot_{h_k^n}\eta_n)\cdot_{h_n^{n+m}}t_{n+m}\in G_{n,m},$$
			whenever $\eta_n\in G_n$ and $t_{n+m}\in G_{n+m}^{(0)}$ are composable for the actor $h_{n}^{n+m}$.
			Define $\rho_k:\bar{G}\to G_k^{(0)}$ by $\rho_k[(\eta_n\cdot t_{n+m})_m]=\rho_{k}^{n}(\eta_n)$, where $\eta_n\in G_{n,0}$ is the first element in the sequence $(\eta_n\cdot t_{n+m})_{m}\in\bar{G}_n$.
			Let $\rho_k^n:G_n\to G_k$ be the anchor map associated to $h_k^n$ (and similarly for $\rho_k^{n+m}$).
			Since $\rho_{k}^{n+m}(\eta_n\cdot t_{n+m})=\rho_k(\eta_n)$ for all $m\in\NN$, we see that the anchor $\rho_k$ is compatible with the inductive limit structure induced by the maps $\bar\imath_n$.
			The actor $j_k$ is then defined as
			$$\gamma_k\cdot[(\eta_n\cdot t_{n+m})_m]=[((\gamma_k\cdot\eta_n)\cdot t_{n+m})_m],$$
			and this well defined by the above arguments.
			
			Let $\ell_n:G_n\curvearrowright H$ be a family of actors for $n\in\NN$ such that $\ell_{n+1} h_n=\ell_n$ with anchor maps $\rho_{\ell_n}:H\to G_n^{(0)}$.
			Define $\rho_\ell:H\to\bar{G}^{(0)}$ by $\rho_\ell(x)=[(\rho_{\ell_{n+m}}(x))_m]$.
			This is well defined for each $\bar{G}_n$ since the connecting maps $p_{n+m}$ in the projective limits $\bar{G}_n$ are given by the respective anchor maps on units in each $G_{n,m}$, and the actors $\ell_n$ respect the inductive limit structure.
			The map $\rho_\ell$ is well defined in the inductive limit $\bar{G}$ since two elements are equivalent in the inductive limit if and only if the tails of their sequences agree, and these are uniquely determined by the element $x\in H$.
			We claim the actor multiplication
			$$[(\eta_n\cdot t_{n+m})_m]\cdot_\ell x=\eta_n\cdot_{\ell_n}x,$$
			with anchor map $\rho_\ell$ is well-defined.
			For $(\eta_n\cdot t_{n+m})_{n}\in \bar{G}_{n}$, we have $s(\eta_n\cdot t_{n+m})=t_{n+m}=\rho_{\ell_n}(x)$ for all $m\in\NN$.
			Since the actors $\ell_n$ form a family compatible with the inductive structure maps $h_n$, we see
			$$\eta_n\cdot_{\ell_n}x=(\eta_n\cdot_{h_n^{n+m}}\rho_{\ell_{n+m}}(x))\cdot_{\ell_{n+m}} x=(\eta_n\cdot_{h_n}t_{n+m})\cdot_{\ell_{n+m}} x.$$
			Hence the multiplication $\cdot_\ell$ is well-defined.
			That $\ell$ is an actor follows pointwise from each $\ell_n$ being an actor.
			
			The actor $\ell$ is uniquely determined by the actors $\ell_n$, whereby $\bar{G}$ is the inductive limit.
			Since the groupoids $\bar{G}_n\subseteq\bar{G}$ form an open cover of $\bar{G}$, it is clear tha tevery actor $\bar{G}\curvearrowright H$ arises in this way.
		\end{proof}
	\end{theorem}
	
	An inductive system of actors induces an inductive system of groupoid $C^*$-algebras, and so a natural question is whether the functor from groupoids to $C^*$-algebras preserves these colimits.
	The next result proves this is the case.
	
	\begin{theorem}
		Let $(G_n,h_n)_{n\in\NN}$ be an inductive system of free and proper actors with surjective anchor maps and let $\varphi_n:C^*_\ess(G_n)\to C^*_\ess(G_{n+1})$ be the induced injective ${}^*$-homomorphisms.
		The inductive limit $C^*$-algebra associated to the system $(C^*_\ess(G_n),\varphi_n)_{n\in\NN}$ is isomorphic to the essential groupoid $C^*$-algebra of $\bar{G}$, the inductive limit of the system $(G_n,h_n)_{n\in\NN}$ from Theorem~\textup{\ref{thm-indLimGrpd}}
		\begin{proof}
			The inductive limit structure actors $j_n:G_n\curvearrowright\bar{G}$ are each free since the actors $h_n$ are, and the anchor maps $\rho_{j_n}:\bar{G}\to G_n^{(0)}$ are each proper and surjective since $\rho_{h_n}$ each are.
			The induced ${}^*$-homomorphisms $i_n:C^*_\ess(G_n)\to C^*_\ess(\bar{G})$ are thus injective by Proposition~\ref{prop-cartanMapInjIffInjOnSubalg}.
			The maps $i_n$ then form a cocone under the inductive system by Lemma~\ref{lem-actorToStarHomFunctorial}.
			Let $A$ be the inductive limit $C^*$-algebra of the system $(C^*_\ess(G_n),\varphi_n)_{n\in\NN}$ and let $\lambda_n:C^*_\ess(G_n)\to A$ be the canonical maps in the  universal cocone, and let $\varphi:A\to C^*_\ess(\bar{G})$ be the canonical map induced by the universal property of $A$ as a colimit.
			Then $i_n=\varphi\circ \lambda_n$ for each $n\in\NN$, whereby each $\lambda_n$ must be injective. 
			The groupoid $\bar{G}$ has a basis of bisections of the form
			$$\left[\prod_{m\in\NN}(U_n\cdot V_{n+m})\right]=\{[(\eta_n\cdot t_{n+m})_m]\in\bar{G}:\eta_n\in U_n,t_{n+m}\in V_{n+m}\text{ for all }m\in\NN\},$$
			where $U_n\subseteq G_n$ is an open bisection and $V_{n+m}\subseteq G_{n+m}^{(0)}$ are open subsets, and $V_{n+m}=G_{n+m}^{(0)}$ for all but finitely many $m\in\NN$.
			Any continuous function with compact support on such a bisection will be the image of a continuous function in $\Cc_c(U_n)\subseteq C^*_\ess(G_n)$ under the inclusion $i_n$.
			Hence the images of the maps $i_n$ span a dense subspace of $C^*_\ess(G)$, whereby $\varphi$ is surjective.
			
			To see that $\varphi$ is injective, we note that on each subspace $\lambda_n(C^*_\ess(G_n))$ for each $n\in\NN$, the map $\varphi$ is given by $\varphi(\lambda_n(f_n))=i_n(f_n)$.
			Thus $\varphi$ is isometric on the subspace $\lambda_n(C^*_\ess(G_n))$ as $i_n$ and $\lambda_n$ are.
			The union of the subspaces $\lambda_n(C^*_\ess(G_n))$ over all $n\in\NN$ is dense in $A$, whereby $\varphi$ is isometric.
		\end{proof}
	\end{theorem}
	
	The construction above lends itself to the more general setting of Fell actors.
	We shall define inductive systems of Fell actors such that the underlying actors of groupoids form an inductive system as defined above, and we may construct a Fell bundle over the inductive limit groupoid representing the colimit of the diagram in the category of Fell bundles and Fell actors.
	As before, we consider saturated categorical Fell bundles and saturated Fell actors throughout.
	
	\begin{definition}
		Let $A:\Ee\curvearrowright\Ff$ be a Fell actor with underlying actor $h:G\curvearrowright H$.
		Let $\rho:H\to G^{(0)}$ denote the anchor map of $h$.
		We say $A$ is \emph{free} if $h$ is, and the ${}^*$-homomorphisms $\nabla_t:\Ee_{\rho(t)}\to M(\Ff_t)$ for each $t\in H^{(0)}$ are injective. 
	\end{definition}
	
	\begin{lemma}\label{lem-freeFellActorEmbedding}
		Let $A:\Ee\curvearrowright\Ff$ be a free proper Fell actor with underlying actor $h:G\curvearrowright H$.
		For each $(\gamma,t)\in G\baltimes{s}{\rho}H^{(0)}$, the map $\nabla_{\gamma,t}:\Ee_\gamma\to\Kk(\Ff_t,\Ff_{\gamma\cdot t})$ is a $\nabla_{r(\gamma\cdot t)}$--$\nabla_t$-equivariant isometric embedding of Hilbert bimodules.
		\begin{proof}
			The map $\nabla_{\gamma,t}$ is left $\nabla_{r(\gamma\cdot t)}$-equivariant by Lemma~\ref{lem-connectingMapsOnFibres}.
			It is right $\nabla_t$-equivariant by Lemma~\ref{lem-connectingMapsKindaMultiplicative}.
			
			To see that $\nabla_{\gamma,t}$ entwines right inner products through $\nabla_t$, we compute an inner product for $\xi_\gamma,\eta_\gamma\in\Ee_\gamma$.
			Recall that by Lemma~\ref{lem-connectingMapsOnFibres} we have $\nabla_{\gamma,t}(\xi_\gamma)^*=\nabla_{\gamma^{-1},\gamma\cdot t}(\xi_\gamma^*)$. We compute:
			\begin{align*}
				\nabla_{\gamma,t}(\xi_\gamma)^*\nabla_{\gamma,t}(\eta_\gamma)&=\nabla_{\gamma^{-1},\gamma\cdot t}(\xi_\gamma^*)\nabla_{\gamma,t}(\eta_\gamma)\\
				&=\nabla_{\gamma^{-1}\gamma,t}(\xi_\gamma^*\eta_\gamma),&\text{by Lemma~\ref{lem-connectingMapsKindaMultiplicative}},\\
				&=\nabla_t(\xi_\gamma^*\eta_\gamma).
			\end{align*}
			Similarly, $\nabla_{\gamma,t}$ entwines left inner products through $\nabla_{r(\gamma\cdot t)}$.
		\end{proof}
	\end{lemma}
	
	\begin{lemma}\label{lem-HDorSkelGivesEssMap}
		Let $A:\Ee\curvearrowright\Ff$ be a free and proper Fell actor with underlying actor $h:G\curvearrowright H$.
		Suppose that either $G$ is Hausdorff or that the anchor map $\rho$ is skeletal.
		If the anchor map $\rho:H\to G^{(0)}$ is surjective, then the induced ${}^*$-homomorphism $\Phi_A:C^*_\ess(\Ee)\to C^*_\ess(\Ff)$ on essential Fell bundle $C^*$-algebras is injective.
		\begin{proof}
			If $G$ is Hausdorff, then $EL_\Ee$ is a genuine conditional expectation and $\Phi_A$ entwines conditional expectations by Proposition~\ref{prop-fellActorEntwines}.
			Since $EL_\Ff$ is faithful, an element $f\in C^*_\ess(\Ee)$ belongs to the kernel of $\Phi_A$ if and only if $0=EL_\Ff(\Phi_A(f^*f))=\Phi_A(EL_\Ee(f^*f))$.
			Hence $EL_\Ee(f^*f)$ lies in the kernel of $\Phi_A$.
			For any $t\in H^{(0)}$ we have
			$$||EL_\Ee(f^*f)(\rho(t))||=||\nabla_{t}(f^*f(\rho(t)))||=||\Phi_A(f^*f)(t)||=0.$$
			Since $\rho$ is surjective, we see that $EL_\Ee(f^*f)=0$ whereby $f=0$ as $EL_\Ee$ is faithful.
			
			Now suppose that $\rho$ is skeletal.
			For $f\in\ker(\Phi_A)$ fix $\varepsilon>0$ and $f_\varepsilon\in\Cc_c(G,\Ee)$ with $||f^*-f_\eps^*f_\eps||<\varepsilon$.
			There exist open bisections $U_1,\dots,U_n\subseteq G$ and sections $f^j_\eps\in\Cc_c(U_j,\Ee)$ such that $f_\eps=\sum_{j=1}^nf_j$.
			Let $V_j:=(U_j\cap G^{(0)})\cup G^{(0)}\setminus\overline{V}_j$, so that $V_j$ is a dense open subset of $G^{(0)}$ and the expectation $EL_\Ee(f_j)$ is given by the restriction of $f_j$ to $V_j$, acting on $C_0(V_j,\Ee)$.
			Let $\tilde{U}:=\bigcap_{j=1}^n V_j$, so that $\tilde{U}$ is a dense open subset of $G^{(0)}$ and the expectation $EL_\Ee(f_\eps)$ of $f_\eps$ is the restriction of $f_\eps$ to $\tilde{U}$.
			Since $\rho$ is skeletal, the preimage $\rho^{-1}(\tilde{U})$ is dense and open in $H^{(0)}$, and $\rho(\rho^{-1}(\tilde{U})) = \tilde{U}$ since $\rho$ is surjective.
			We then have
			$$||EL(f^*_\eps f_\eps)||=\sup_{x\in\tilde{U}}||f^*_\eps f_\eps(x)||\leq \sup_{t\in\rho^{-1}(\tilde{U})}||\nabla_t(f^*_\eps f_\eps(\rho(t)))||=||EL_\Ff(\Phi(f^*_\eps f_\eps)||<\eps,$$ 
			as each $\nabla_t$ is isometric and $C_0(\rho^{-1}(\tilde{U}),\Ff)$ is an essential ideal of $C_0(H^{(0)},\Ff)$.
			Taking $\eps\to 0$ we see that $EL_\Ee(f^*f)=0$, whereby $f=0$ as $EL_\Ee$ is faithful.
		\end{proof}
	\end{lemma}
	
	We now define the particular class of inductive systems of Fell actors we wish to consider.
	
	\begin{definition}
		Let $\Ee_n\to G_n$ be saturated categorical Fell bundles over \'etale groupoids $G_n$ with locally compact Hausdorff unit spaces, and let $(A_n:\Ee_n\curvearrowright\Ee_{n+1})_{n\in\NN}$ be proper saturated Fell actors.
		We say the collection $(\Ee_n,A_n)_{n\in\NN}$ is an \emph{inductive system of Fell actors} if each $A_n$ is free, and the anchor maps $\rho_n:G_{n+1}\to G_n^{(0)}$ are surjective.
		
		We say an inductive system of Fell actors $(\Ee_n,A_n)_{n\in\NN}$ is \emph{essentially compatible} if the induced ${}^*$-homomorphisms $\Phi_{A_n}$ arising from $A_n$ descend to injective maps on the essential Fell bundle $C^*$-algebras.
	\end{definition}
	
	The term `essentially compatible' is used so that we may consider essential $C^*$-algebras throughout.
	If for each $n\in\NN$ either $G_n$ is Hausdorff or the anchor map associated to the actor $A_n$ in an inductive system is skeletal, the system is essentially compatible by Lemma~\ref{lem-HDorSkelGivesEssMap}.
	
	Let $\Ee_n\to G_n$ be saturated categorical Fell bundles over \'etale groupoids $G_n$ with locally compact Hausdorff unit spaces, and let $(A_n:\Ee_n\curvearrowright\Ee_{n+1})_{n\in\NN}$ be an essentially compatible inductive system of Fell actors, and let $\bar{G}$ be the inductive limit groupoid for the underlying actors $h_n:G_n\curvearrowright G_{n+1}$ in Theorem~\textup{\ref{thm-indLimGrpd}}.
	Over each projective limit groupoid $\bar{G}_n=\varprojlim\limits_m G_{n,m}$ we can define the bundle $\bar\Ee_n$ over $\bar{G}_n$ fibrewise by declaring the fibre over a generic element $\bar{\gamma}_n=(\gamma_n,\gamma_n\cdot t_{n+1},\gamma_n\cdot t_{n+2},\dots)\in\bar{G}_n$ as $(\Ee_n)_{\gamma_n}$.
	The Fell bundle structure of $\bar\Ee_n$ is then inherited from $\Ee_n$ in this way, since the groupoid operations in $\bar{G}_n$ are pointwise.
	
	Recall the inductive system $(\bar{G}_n,\bar\imath_n)_{n\in\NN}$ from Lemma~\ref{lem-inductOnProjGrpdLem}.
	For each $n\in\NN$ we also define injective Fell bundle maps $I_n:\bar\Ee_n\to\bar\Ee_{n+1}$ lifting the maps $\bar\imath_n:\bar{G}_n\to\bar{G}_{n+1}$ as follows: for $\bar\gamma_n=(\gamma_n,\gamma_n\cdot t_{n+1},\dots)\in\bar{G}_n$ we note that $(\bar\Ee_n)_{\bar\gamma_n}=(\Ee_n)_{\gamma_n}$ and $(\bar\Ee_{n+1})_{\bar\imath_n(\bar\gamma_n)}=(\Ee_{n+1})_{\gamma_n\cdot t_{n+1}}$.
	We identify $(\Ee_{n+1})_{\gamma_n\cdot t_{n+1}}$ with $\Kk((\Ee_n)_{t_{n+1}},(\Ee_n)_{\gamma_n\cdot t_{n+1}})$ via the canonical map $\xi\mapsto |\xi\rangle$ sending vectors in $(\Ee_{n+1})_{\gamma_n\cdot t_{n+1}}$ to their left creation operators.
	Define $(I_n)_{\bar\gamma_n}:=\nabla_{\gamma_n,t_{n+1}}:(\Ee_n)_{\gamma_n}\to(\Ee_{n+1})_{\gamma_n\cdot t_{n+1}}$.
	The $I_n$ are Fell bundle morphisms by Lemma~\ref{lem-connectingMapsKindaMultiplicative}.
	These form an inductive system of Banach spaces, and so we may define the inductive limit by
	$$\bar\Ee_{[\bar\gamma_n]}:=\varinjlim\limits_{m}((\Ee_{n+m})_{\gamma_n\cdot t_{n+m}},\nabla_{\gamma_n\cdot t_{n+m},t_{n+m+1}}).$$
	
	\begin{proposition}
		Equipped with pointwise operations, these form the fibres of a Fell bundle over $\bar{G}$.
		The topology on $\bar\Ee$ is given by specifying the local continuous sections: a local section $f:U\to\bar\Ee$ is continuous if for all $\eps>0$ there is $m\in\NN$ and a continuous section $f':U\to\Ee_{n+m}\subseteq\bar\Ee$ such that $||f[\gamma_n]-f_\eps[\gamma_n]||<\varepsilon$ for all $[\gamma_n]\in U$.
			
		Moreover $\bar\Ee$ is saturated and categorical.
		\begin{proof}
			Indeed, if $f$ is such a section, then for all $\alpha>0$ the set
			$$\{x\in U: ||f(x)||<\alpha\}=\bigcup_{\varepsilon>0}\{x\in U:||f(x)||<\alpha-\varepsilon\}.$$
			For each $\varepsilon>0$ we may pick a local section $f_\eps$ such that $||f(x)-f_\eps(x)||<\eps$ for all $x\in U$, giving $||f(x)||\leq ||f(x)-f_\eps(x)||+||f_\eps(x)||<\varepsilon+f_\eps(x)$.
			
			Moreover, we may without loss of generality pick $f_\eps$ with $||f_\eps(x)||<||f(x)||$ by approximating $(1-\delta)f$ for $\delta>0$.
			Fix such a function $f_\eps$ for each $\eps >0$.
			Let $V\subseteq U$ be the set $V:=\{x\in U: ||f_\eps(x)||<\alpha-\eps$ for $\eps>0\}$.
			Fix $x\in U$ with $||f(x)||<\alpha$.
			Then for any such $f_\eps$ as above we have $||f_\eps(x)||\leq ||f(x)||<\alpha$, so $x\in V$ and $U\subseteq V$.
			Now suppose $x\in V$, set $\eps<\alpha-||f(x)||$.
			Picking $f_\eps$ as above, we have
			$$||f(x)||\leq ||f(x)-f_\eps(x)||+||f_\eps(x)||<\alpha-\eps+\eps=\alpha,$$
			so $x\in U$.
			Hence $U=V$ is open, as it is a union of open sets, and \cite[Proposition~2.4]{BE1} gives a unique topology on $\bar\Ee$ turning it into a Banach bundle such that the local sections $f$ described above are continuous.
			There is a basis for the topology on $\bar\Ee$ given by open `sleeves'
			$$\Omega(U,f,\eps):=\{v\in\Ee: p(v)\in U, ||v-f(p(v))||<\eps\},$$
			where $f$ is a local section as described above, and $\eps>0$.
			
			We shall use \cite[Proposition~2.9]{BE1} to show that $\bar\Ee$ is a Fell bundle.
			Equipping $\bar\Ee$ with the (limits of the) pointwise products in $\Ee_n$, the algebraic axioms follow from the pointwise operations together with continuity of the multiplication and involution (which we shall now show).
			Fix two continuous local sections $f:U\to\Ee$ and $g:V\to\Ee$.
			Fix $\eps>0$, and pick $f_\eps$ and $g_\eps$ as above.
			Then the pointwise product map $(\gamma,\eta)\mapsto f_\eps(\gamma)g_\eps(\eta)$ is continuous on $G^{(2)}\cap U\times V$, since it is a product of sections of $\Ee_n$ for some $n\in\NN$ large enough.
			Taking uniform limits converging on $f$ and $g$ then preserves continuity, so $fg$ is continuous.
			Similarly, the map $U\to\Ee$, $\gamma\mapsto f(\gamma)^*$ is continuous since it can be uniformly approximated by continuous sections of $\Ee_n$.
			Thus \cite[Proposition~2.9]{BE1} implies that the Fell bundle structure is 
			
			The resulting Fell bundle will be saturated since the inductive system Fell bundles all are, and the union of the images of these bundles in the inductive limit forms a dense subbundle of $\bar\Ee$.
			
			We must now show that $\bar\Ee$ is categorical.
			Since $\Ee_n$ is categorical, the fibres of $\Ee_n$ over units of $\bar{G}_n$ are unital $C^*$-algebras, and so each $\nabla_{t_{n+m}}$ is a unital ${}^*$-homomorphism by Lemma~\ref{lem-multDefinesStarHomsOnUnitSpaceFibres}.
			Hence the inductive limit fibres over units are also unital, and the unit is the image of the unit of any of the fibres in the inductive limit building blocks.
			The source map on $\bar\Ee$ is given by the bundle map composed with source map in the groupoid, together with the section of units in $\Ee_n$ which exists by Lemma~\ref{lem-localUnitSections}.
			In particular, it is composition of continuous functions, so is itself continuous.
			A similar argument shows the range map is continuous.
		\end{proof}
	\end{proposition}
	
	\begin{theorem}\label{thm-indLimFellBundle}
		There are free proper saturated Fell actors $J_k:\Ee_k\curvearrowright\bar\Ee$ with underlying actors $j_k:G_k\curvearrowright\bar{G}$ from Theorem~\textup{\ref{thm-indLimGrpd}} turning $\bar\Ee$ into the inductive limit of the system $(A_n:\Ee_n\curvearrowright\Ee_{n+1})_{n\in\NN}$.
		\begin{proof}		
			For $n\in\NN$, we specify the underlying actor of $J_k:\Ee_k\curvearrowright\bar\Ee$ to be the $j_k:G_k\curvearrowright\bar{G}$ as in the proof of Theorem~\ref{thm-indLimGrpd}.
			For $\gamma_k\in G_k$ and $[\bar\gamma_n]=[(\gamma_{n,m})_{m\in\NN}]\in\bar{G}$ with $\rho_k[\bar\gamma_n]=s(\gamma_k)$ we may without loss of generality assume that $n\geq k$.
			We define $M^{J_k}_{\gamma_k,[\bar\gamma_n]}:(\Ee_k)_{\gamma_k}\times\bar\Ee_{[\bar\gamma_n]}\to\bar\Ee_{[\gamma_k\cdot\bar\gamma_n]}$ by
			$$M^{J_k}_{\gamma_k,[\bar\gamma_n]}(\xi_k,[(\zeta_{n,m})_m]):=[(M^{A^{n+m}_k}_{\gamma_k,\gamma_{n+m}}(\xi_k,\zeta_{n,m}))_m],$$
			where $M^{A_k^{n+m}}_{\gamma_k,\gamma_n}$ is the multiplication in the actor $A_k^{n+m}:\Ee_k\curvearrowright\Ee_{n+m}$ in the inductive system.
			This is well defined since if $[(\zeta_{n,m})_m]=[(\zeta'_{n,m})_m]$ in $\bar\Ee_{[\bar\gamma_n]}$ then for large enough $n\in\NN$ we have $(\zeta_{n,m})_m=(\zeta_{n,m}')_m$, and so for such $n$ we have $(M_{\gamma_k,\gamma_{n+m}}(\xi_k,\zeta_{n,m}))_m]=[(M_{\gamma_k,\gamma_{n+m}}(\xi_k,\zeta'_{n,m}))_m$.
			This multiplication satisfies axioms (1)--(4) of Definition~\ref{defn-FellActor} since each $A_k^{n+m}$ does pointwise.
			The actors $J_k$ are free, proper, and saturated since $A_k^{n+m}$ are.
			
			By definition of composite actors in Proposition~\ref{prop-defineFellActComp} we have 
			\begin{align*}
				M^{J_{k+1}A_k}(\xi_k,[(\zeta_{n,m})_m])&=M^{J_{k+1}}(M^{A_k}(\xi_k,1_{\rho_{k+1}([\bar\eta_n])}),[(\zeta_{n,m})_m])\\
				&=[(M^{A^{n+m}_{k+1}}(M^{A_k}(\xi_k,1_{\rho_{k+1}([\bar\eta_n])}),\zeta_{n,m}))_m]\\
				&=[(M^{A^{n+m}_k}(\xi_k,\zeta_{n,m}))_m]\\
				&=M^{J_k}(\xi_k,[(\zeta_{n,m})_m]),
			\end{align*}
			whereby $J_k=J_{k+1}A_k$, so the actors $J_k$ form a cocone of the inductive system $(A_n:\Ee_n\curvearrowright\Ee_{n+1})_{n\in\NN}$.
			
			Given a Fell bundle $\Ff$ over an \'etale groupoid $H$ and a system of proper actors $(B_n:\Ee_n\curvearrowright\Ff)_{n\in\NN}$ with underlying actors $\ell_n:G_n\curvearrowright H$ satisfying $B_n=B_{n+1}A_n$, we define $B:\bar\Ee\curvearrowright\Ff$ as follows.
			Let $\ell:\bar{G}\curvearrowright H$ be the actor induced by the inductive limit structure of $\bar{G}$: for any leg $j_k:G_k\curvearrowright\bar{G}$ as described in Theorem~\ref{thm-indLimGrpd} we have $\ell_k=\ell j_k$.
			This shall serve as the underlying actor for $B$.
			Let $\rho:H\to\bar{G}^{(0)}$ be the anchor map for $\ell$.
			For $([\bar\gamma_n],x)\in\bar{G}\baltimes{s}{\rho}H$ define $M^B_{[\bar\gamma_n],x}:\bar\Ee_{[\bar\gamma_n]}\times\Ff_x\to\Ff_{[\bar\gamma_n]\cdot x}$ by
			$$M^B_{[\bar\gamma_n],x}([(\zeta_{n,m}],f_x):=M^{B_n}_{\gamma_n,x}(\zeta_{n,0},f_x),$$
			where $M^{B_n}_{\gamma_n,x}$ is the multiplication arising from the actor $B_n$.
			This is well-defined, since for any $m\in\NN$ we have 
			$$M^{B_n}_{\gamma_n\cdot t_{n+m},x}(\zeta_{n,0},f_x)=M^{B_{n+m}}_{\gamma_n\cdot t_{n+m},x}\left(M^{A_n^{n+m}}_{\gamma_n,t_{n+m}}(\zeta_{n,0}, 1_{t_{n+m}}),f_x\right),$$
			since $\zeta_{n,0}\in(\Ee_n)_{\gamma_n}$ satisfies $\zeta_{n,m}=M^{A_n^{n+m}}_{\gamma_n,t_{n+m}}(\zeta_{n,0},1_{t_{n+m}})=\zeta_{n,m}$, where $1_{t_{n+m}}\in(\Ee_t){t_{n+m}}$ is the unit.
			It follows from construction that $B=B_nJ_n$ for any $n\in\NN$, and that these $B_n$ are the unique such actors.
			Since $\bar{G}_n\subseteq\bar{G}$ form an open cover of $\bar{G}$ and $\Ee_n$ span a dense subbundle of $\bar\Ee$, every Fell actor $\bar\Ee\curvearrowright\Ff$ arises this way.
		\end{proof}
	\end{theorem}
	
	\begin{lemma}\label{lem-HDorSkelFellIndSyst}
		Let $(\Ee_n,A_n)_{n\in\NN}$ be an inductive system of Fell actors where for each $n\in\NN$ either $G_n$ is Hausdorff or $\rho_n$ is skeletal, and let $J_k:\Ee_k\curvearrowright\bar\Ee$ be the Fell actors from Theorem~\textup{\ref{thm-indLimFellBundle}}.
		Then the induced ${}^*$-homomorphisms $I_k:C^*(\Ee_k)\to C^*(\bar\Ee)$ descend to a map of the essential quotients $C^*_\ess(\Ee_k)\to C^*_\ess(\bar\Ee)$.
		\begin{proof}
			We divide into two cases: either there are infinitely many $n\in\NN$ such that $G_n$ is Hausdorff, or eventually all anchor maps $\rho_n$ are skeletal.
			
			Suppose firstly that there are infinitely many $n\in\NN$ such that $G_n$ is Hausdorff.
			Since each $A_n$ is free, proper, and saturated, the induced ${}^*$-homomorphisms $\Phi_n$ descend to maps $C^*_\ess(\Ee_n)\to C^*_\ess(\Ee_{n+1})$ of the essential $C^*$-algebras by Lemma~\ref{lem-HDorSkelGivesEssMap}.
			For a Fell bundle $\Ee_n$ in the system, let $k\in\NN$ be the next largest natural number bigger than $n$ such that $G_k$ is Hausdorff.
			The induced ${}^*$-homomorphism associated to $J_k$ descends to essential $C^*$-algebras by Lemma~\ref{lem-HDorSkelGivesEssMap} since $G_k$ is Hausdorff, and we have already shown that $A_{k-1}\dots A_n$ has the same property.
			Hence the ${}^*$-homomorphism $\Phi_{J_kA_{k-1}\dots A_n}$ descends to essential $C^*$-algebras.
			By Proposition~\ref{prop-fellActorHomFunctorial} we have $\Phi_{J_kA_{k-1}\dots A_n}=\Phi_{J_k}\circ\Phi_{A_{k-1}\dots A_n}$, which descends to essential $C^*$-algebras since both $\Phi_{J_k}$ and $\Phi_{A_{k-1}\dots A_n}$ do.
			Lastly, by Theorem~\ref{thm-indLimFellBundle} we have $J_kA_{k-1}\cdots A_n=J_n$, whereby $\Phi_{J_n}=\Phi_{J_kA_{k-1}\dots A_n}$ descends to essential groupoids $C^*$-algebras.
			
			Now suppose that there are only finitely many Hausdorff $G_n$ in the inductive system.
			For $n\in\NN$, if there exists $k>n$ with $G_k$ Hausdorff, then we repeat the process above.
			Otherwise, for all $k>n$ the anchor $\rho_k:G_{k+1}\to G_k^{(0)}$ is skeletal.
			Without loss of generality we may assume all anchor maps $\rho_n$ are skeletal for all $n\in\NN$ (by truncating the inductive system).
			Recall the anchor $\rho_{J_k}:\bar{G}\to G_k^{(0)}$ is given by
			$$\rho_{j_k}([(\gamma_{\ell,m})_m])=\rho_k^\ell(\gamma_{\ell,0}),$$
			where $\gamma_{\ell,0}$ is the first element in the sequence $(\gamma_{\ell,m})_m$.
			Note we can always take $\ell>k$ since a sequence in the inductive limit $\bar{G}=\lim\limits_{\ell\to\infty}G_\ell$ is equivalent to all its (finite) truncations from the left.
			The map $\rho_k^\ell=\rho_\ell\circ\rho_{\ell-1}\circ\dots\circ\rho_k$ is skeletal, since each of the $\rho_k,\dots,\rho_\ell$ is.
			Fix $C_1\subseteq G_1^{(0)}$ nowhere dense, and define recursively $C_{\ell+1}=\rho_{\ell}^{-1}(C_\ell)$ for $\ell\geq 1$.
			Let $C=\prod_{\ell\in\NN}C_\ell\subseteq \bar{G}_1=\varprojlim_{\ell}G_\ell$.
			Then $$\overline{C}^\circ=\overline{\left(\prod_{\ell\in\NN}C_\ell\right)}^{\circ}=\left(\prod_{\ell\in\NN}\overline{C_\ell}\right)^\circ=\emptyset,$$
			since the product of closures of sets is the closure of the product in the product topology, $\bar{C_\ell}\neq G_{\ell}$ for any $\ell$ since each $G_\ell$ is Baire.
			The map $\bar{G}_1\to\bar{G}$ induced by the inductive limit structure is a homeomorphism onto its image, and the image of $C$ under this inclusion is exactly
			$$\{[(\gamma_{1,\ell})_\ell]:\gamma_{1,\ell}\in C_\ell\}=\{[(\gamma_{1,\ell})_\ell]:\gamma_{1,\ell}\in (\rho_1^\ell)^{-1}(C_1)\}=\rho_{J_k}^{-1}(C_1).$$
			Hence $\rho_{J_k}^{-1}(C_1)\cong C$ is nowhere dense, so $\rho_{J_k}$ is skeletal.
			Proposition~\ref{prop-fellActorEntwines} then shows that $I_k$ descends to a map $C^*_\ess(\Ee_k)\to C^*_\ess(\bar\Ee)$.
		\end{proof}
	\end{lemma}

	\begin{lemma}\label{lem-indLimBundleDenseSpanned}
		Let $I_k:C^*_\ess(\Ee_k)\to C^*_\ess(\bar\Ee)$ be the ${}^*$-homomorphisms induced by the actors $J_k:\Ee_k\curvearrowright\bar\Ee$ in Theorem~\textup{\ref{thm-indLimFellBundle}}.
		The images of $C^*_\ess(\Ee_k)$ under the maps $I_k$ span a dense subspace of $C^*_\ess(\bar\Ee)$.
		\begin{proof}
			By Lemma~\ref{lem-HDorSkelFellIndSyst} the maps $I_k$ induced by the actors $J_k$ all descend to maps $C^*_\ess(\Ee_k)\to C^*_\ess(\bar\Ee)$ of the essential $C^*$-algebras.
			
			Fix an open bisection $\tilde{U}\subseteq \bar{G}$ and $f\in C_c(\tilde{U},\bar\Ee)$.
			Note that $\tilde{U}$ takes the form 
			$$\tilde{U}=\left[\prod_{m\in\NN}U_n\cdot V_{n+m}\right]$$
			for some bisections $U_n\subseteq G_n$ and $V_{n+m}\subseteq G_{n+m}^{(0)}$, where $n\in\NN$ is large enough and $V_{n+m}=G_{n+m}^{(0)}$ for all but finitely many $m\in\NN$.
			Since the fibres $(\Ee_k)_{\gamma_n\cdot t_{k}}$ span a dense subspace of $\bar\Ee_{[\bar\gamma]}$ for each $[\bar\gamma]\in\tilde{U}$, for $k\in\NN$ large enough, we can pick a function $f_k^{[\bar{\gamma}]}\in \Cc_c(U_n\cdot V_k,\Ee_k)$ with $||I_k(f_k)[\bar{\gamma}']-f[\bar{\gamma}']||<\varepsilon$ for all $[\bar{\gamma}']$ in a precompact neighbourhood $\tilde{U}_{[\bar{\gamma}]}$ of $[\bar{\gamma}]$.
			The open sets $\tilde{U}_{[\bar{\gamma}]}$ form an open cover of the compact support of $f$, so we may take a finite subcover $\tilde{U}_{[\bar{\gamma}_1]},\dots,\tilde{U}_{[\bar{\gamma}_q]}$.
			Taking the maximum such $k$ for the points $[\bar{\gamma}_1],\dots,[\bar{\gamma}_q]$, we may assume the functions $f_k^{[\bar{\gamma}_1]},\dots,f_k^{[\bar{\gamma}_q]}$ are sections of the same bundle $\Ee_k$.
			Writing $U^i=\tilde{U}_{[\bar{\gamma}_i]}$ and $f_k^i=f_k^{[\bar{\gamma}_i]}$, the open sets $U^i$ take the form
			\[U^i=\left[\prod_{m\in\NN}U^i_n\cdot V^i_{m+n}\right],\label{eq1}\tag{$*$}\]
			where each $U^i_n\subseteq G_n$ is an open bisection and each $V^i_{m+n}$ is an open subset of $G_{n,m}^{(0)}\subseteq G_{n+m}$ for some $n\in\NN$.
			We may without loss of generality assume $k\geq n$.
			Since the equivalence relation in $\bar{G}$ is given by asymptotic equality under the unilateral shift, the set $U^i$ uniquely determines $U^i_k=U^i_n\cdot V^i_{n,k-n}\subseteq G_{n,k-n}$ in the expression (\ref{eq1}).
			Let $\varsigma_i:G_k^{(0)}\to[0,1]$ be a continuous partition of unity suboordinate to the family $(s(U^i_k))_{i=1}^q$, and note that since $U^i$ is precompact, so is $U^i_k$, so each $\varsigma_i$ has compact support.
			We also note that since the sections $f$ and $I_k(f^i_k)$ are all supported on the bisection $\tilde{U}$, their norms are given by the supremum norm over the bisection.
			We then have
			\begin{align*}
				\left|\left|\sum_{i=1}^q I_k(f_k^i\varsigma_i)-f\right|\right|&=\sup\limits_{[\bar{\gamma}]\in\tilde{U}}\left|\left|\sum_{i=1}^q f_k^i(\gamma_k)\varsigma_i(s(\gamma_k))-\sum_{i=1}^q\varsigma_i(s(\gamma_k))f[\bar{\gamma}]\right|\right|\\
				&\leq\sup\limits_{[\bar{\gamma}]\in\tilde{U}}\sum_{i=1}^q||\varsigma_i(s(\gamma_k))\left(f_k^i(\gamma_k)-f[\bar{\gamma_k}]\right)||\\
				&\leq \sum_{i=1}^q\sup\limits_{[\bar{\gamma}]\in U^i}\varsigma_i(s(\gamma_k))||f_k^i(\gamma_k)-f[\bar{\gamma}]||\\
				&< \sum_{i=1}^q\sup\limits_{[\bar{t}]\in V_{[\bar{t}_i]}}\varsigma_i(s(\gamma_k))\varepsilon\\
				&\leq\varepsilon.
			\end{align*}
			Thus the union of the images of the algebras $\Cc_c(G_k,\Ee_k)$ is dense in each $\Cc_c(\tilde{U},\bar\Ee)$ for any open bisection $\tilde{U}\subseteq\bar{G}$, which span a dense subspace of $C^*_\ess(\bar\Ee)$.
		\end{proof}
	\end{lemma}
	
	We now give a version of Theorem~\ref{thm-indLimGrpd} for essentially compatible inductive systems of Fell bundles.
	
	\begin{theorem}
		Let $\Phi_n:C^*_\ess(\Ee_n)\to C^*_\ess(\Ee_{n+1})$ be the injective ${}^*$-homomorphisms induced by the free proper Fell actors $A_n:\Ee_n\curvearrowright\Ee_{n+1}$ in the essentially compatible inductive system above.
		Let $\bar\Ee\to\bar{G}$ be the inductive limit bundle as in Theorem~\textup{\ref{thm-indLimFellBundle}}.
		The essential $C^*$-algebra of $\bar\Ee$ is (naturally isomorphic to) the inductive limit $C^*$-algebra of the system $(C^*_\ess(\Ee_n),\Phi_n)_{n\in\NN}$.
		\begin{proof}
			The Fell actors $J_k:\Ee_k\curvearrowright\bar\Ee$ from Theorem~\ref{thm-indLimFellBundle} are free, proper, and saturated, so give rise to injective ${}^*$-homomorphisms $I_k:C^*_\ess(\Ee_k)\to C^*_\ess(\bar\Ee)$ satisfying $I_k=I_{k+1}\circ\Phi_k$ by Lemma~\ref{lem-HDorSkelFellIndSyst}.
			
			Let $B$ be the inductive limit of $(C^*_\ess(\Ee_n),\Phi_{A_n})$ in the category of $C^*$-algebras, and let $\Lambda_n:C^*_\ess(\Ee_n)\to B$ be the universal ${}^*$-homomorphisms embedding the system into $A$.
			$B$ can be identified with the subalgebra of $\frac{\prod_{k\in\NN} B_k}{\bigoplus_{k\in\NN}B_k}$ consisting of equivalence classes $[b_k]_{k\in\NN}$ of sequences that are zero up to a point $k_0$, and then $b_{k+1}=\Phi_{A_k}(b_k)$ for all $k\geq k_0$.
			The maps $\Lambda_n$ then map an element $f_n\in C^*_\ess(\Ee_n)$ to the equivalence class of the sequence $(0,\dots,0,f_n,\Phi_n^{n+1}(f_n),\dots)$, where the first $f_n$ occurs at the $n$-th entry of the sequence.
			Such a sequence is in the equivalence class of zero if and only if the sequence of norms $||\Phi_n^{n+\ell}(f_n)||$ converges to zero, which can only occur if $f_n=0$ since each $\Phi_n$ is isometric.
			Hence $\Lambda_n(f_n)=0$ only if $f_n=0$ whereby $\Lambda_n$ is injective, hence isometric.
			There is a universal homomorphism $\Phi:B\to C^*_\ess(\bar\Ee)$ induced by the inclusion maps $I_n$, so that $I_n=\Phi\circ\Lambda_n$ for all $n\in\NN$.
			The map $\Phi$ is then surjective, since the images of the $C^*$-algebras $C^*_\ess(\Ee_n)$ under the maps $I_n$ span a dense subalgebra of $C^*_\ess(\bar\Ee)$ by Lemma~\ref{lem-indLimBundleDenseSpanned}.
			
			To see that $\Phi$ is injective, we note that on each subspace $\Lambda_n(C^*_\ess(\Ee_n))$ for each $n\in\NN$, the map $\Phi$ is given by $\Phi(\Lambda_n(f_n))=I_n(f_n)$.
			Thus $\Phi$ is isometric on the subspace $\Lambda_n(C^*_\ess(\Ee_n))$ as $I_n$ and $\Lambda_n$ are.
			The union of the subspaces $\Lambda_n(C^*_\ess(\Ee_n))$ over all $n\in\NN$ is dense in $A$, whereby $\Phi$ is isometric.
		\end{proof}
	\end{theorem}

\end{document}